\documentclass[11pt,a4paper]{amsart}

\usepackage[T1]{fontenc}
\usepackage[utf8]{inputenc}
\usepackage{lmodern}
\usepackage{amsmath,amssymb,amsfonts,amsthm,mathtools}
\usepackage{amsthm}

\usepackage{mathrsfs}
\usepackage{enumitem}
\usepackage{booktabs,tabularx}
\usepackage{xltabular}
\usepackage{hyperref}
\usepackage[margin=1.15in]{geometry}
\usepackage{microtype}
\allowdisplaybreaks[2]

\hypersetup{
  colorlinks=true,
  linkcolor=blue,
  citecolor=blue,
  urlcolor=blue,
  pdftitle={Quantum Generalized Cluster Algebras from Scattering Diagrams: Positivity and Mutation},
  pdfauthor={Huihui Ye},
  pdfsubject={Positivity and mutation of generalized quantum cluster scattering diagrams},
  pdfkeywords={generalized quantum cluster algebra, scattering diagram, cylinder projection, quantum mutation, F-polynomial}
}

\newtheorem{theorem}{Theorem}[section]
\newtheorem{proposition}[theorem]{Proposition}
\newtheorem{lemma}[theorem]{Lemma}
\newtheorem{corollary}[theorem]{Corollary}

\theoremstyle{definition}
\newtheorem{definition}[theorem]{Definition}
\newtheorem*{definition*}{Definition}

\newtheorem{example}[theorem]{Example}

\theoremstyle{remark}
\newtheorem{remark}[theorem]{Remark}

\newcommand{\ZZ}{\mathbb Z}
\newcommand{\QQ}{\mathbb Q}
\newcommand{\RR}{\mathbb R}

\newcommand{\g}{\mathfrak g}
\newcommand{\Ghat}{\widehat G}
\newcommand{\ghat}{\widehat{\mathfrak g}}
\newcommand{\Scat}{\operatorname{Scat}}

\newcommand{\Ad}{\operatorname{Ad}}
\newcommand{\ord}{\operatorname{ord}}
\newcommand{\pr}{\operatorname{pr}}
\newcommand{\supp}{\operatorname{Supp}}
\newcommand{\Joints}{\operatorname{Joints}}

\newcommand{\id}{\operatorname{id}}
\newcommand{\sgn}{\operatorname{sgn}}

\newcommand{\eps}{\varepsilon}
\newcommand{\Db}{\overline{\mathfrak D}}
\newcommand{\Ddag}{\mathfrak D^\dagger}
\newcommand{\Dstar}{\mathfrak D^*}
\newcommand{\Din}{\mathfrak D_{\mathrm{in}}}

\newcommand{\uf}{\mathrm{uf}}

\newcommand{\E}{\mathbb E}

\newcommand{\PathProd}{\mathfrak p}
\newcommand{\WC}{\Theta}
\newcommand{\EssSupp}{\operatorname{EssSupp}}
\newcommand{\codim}{\mathrm{codim}}

\title[Quantum generalized cluster scattering diagrams]{Quantum Generalized Cluster Algebras from Scattering Diagrams}
\author{Meiling Yang}
\address{School of Mathematics and Statistics, Jiangsu Normal University, Xuzhou 221116, China}
\author{Huihui Ye}
\address{School of Mathematics and Statistics, Jiangsu Normal University, Xuzhou 221116, China}
\email{yehuihui@jsnu.edu.cn}


\begin{document}
\raggedbottom

\begin{abstract}

We construct generalized quantum cluster scattering diagrams over a ring of formal reciprocal exchange coefficients, for compatible seeds satisfying \(\widetilde B^T\Lambda=[I_n\mid0].\) A central lattice extension records the roots of the exchange polynomials as independent degrees. Applying Davison--Mandel positivity in this extension and proving a positive cylinder projection theorem yields positive wall factorizations and theta coefficients in the root variables. Their regular classical limits recover generalized cluster scattering diagrams. Nonvanishing under classical specialization, together with control of the mutated positive monoids, gives quantum mutation compatibility. The resulting rational scattering atlas has the generalized cluster fan, and every reachable cluster monomial is a theta function. Classical detection of theta coefficients gives a geometric proof of the Laurent phenomenon, with positivity after adjoining the polynomial roots. Evaluation of the exchange coefficients recovers the generalized quantum cluster algebras with this normalization. The Fock--Goncharov mutation product becomes a wall-crossing action, identifying extended \(g\)-vectors with leading theta exponents and \(F\)-polynomials with the coefficient factors on their right. Constructing these polynomials before evaluation gives the separation formula for arbitrary prescribed reciprocal Laurent polynomial coefficients, without requiring positivity at \(t=1.\) Positivity in the roots is distinguished throughout from positivity in the exchange coefficients.
\end{abstract}

\maketitle
 \enlargethispage{4pt}

\section{Introduction}

Cluster algebras were introduced by Fomin and Zelevinsky to study total positivity and canonical bases through a recursive operation called mutation \cite{FZ02}. Their Laurent phenomenon makes the expansion of a cluster variable in any cluster a finite algebraic object. Two further questions concern the signs of its coefficients and the data that determine the expansion. The separation formulas of \cite{fomin_zelevinsky_2007} address the latter question through \(g\)-vectors and
\(F\)-polynomials. Berenstein and Zelevinsky introduced quantum cluster algebras and proved their Laurent phenomenon \cite{BZ05}; Tran subsequently constructed quantum \(F\)-polynomials and the corresponding separation formulas \cite{Thao2011}.

Generalized cluster algebras allow exchange polynomials of arbitrary positive degree. They arose in the work of Chekhov and Shapiro on Teichm\"uller spaces of surfaces with orbifold points
\cite{Chekhov_Shapiro14}. Nakanishi developed their seed structure and separation formulas \cite{Nakanishi15}, and studied quantized coefficient mutations using quantum dilogarithms of higher degrees
\cite{Nakanishi-qgca}. Bai, Chen, Ding, and Xu introduced generalized quantum cluster algebras with prescribed reciprocal Laurent polynomial coefficients \cite{BCDX18} and established their Laurent phenomenon in arbitrary rank \cite{BCDX23}. These results give an algebraic framework for higher-degree quantum exchange relations.

Scattering diagrams provide a geometric approach to mutation and
positivity. Originating in mirror symmetry \cite{Kontsevich_Soibelman06,Gross_Siebert11}, they attach automorphisms to walls and construct theta functions by summing contributions of broken lines. Gross, Hacking, Keel, and Kontsevich used this construction to realize cluster monomials as theta functions and prove classical positivity \cite{GHKK}. Davison and Mandel developed quantum theta bases and proved strong positivity using quantum scattering diagrams and Donaldson--Thomas theory \cite{Davison_Mandel21}. For generalized cluster algebras, Cheung, Kelley, and Musiker constructed scattering diagrams and theta functions in the reciprocal case \cite{CKM}, while Mou developed a construction with factorized exchange polynomials \cite{Mou24}. Burcroff, Lee, and Mou proved positivity of generalized cluster scattering diagrams and the resulting classical Laurent expansions using tight gradings \cite{BLM25,BLM25a}.

For generalized quantum cluster algebras, the relation between these geometric constructions and the recursive exchange rules requires additional work. A quantum wall corresponding to a higher-degree exchange polynomial contains several root factors. Applying quantum positivity while retaining their separate degrees must be compatible with projection, mutation, and coefficient evaluation. A scattering realization must also identify its theta functions with the variables of the recursive algebra. The cited results provide an algebraic mutation pattern and a classical scattering model; comparing their quantum wall and cluster coordinates requires the construction developed here.

A related question concerns polynomial separation. Fu, Peng, and Ye
established the Fock--Goncharov decomposition of generalized quantum mutations and obtained a quantum Gupta product for the associated rational \(F\)-elements \cite{FPY24}, extending the ordinary formulas of Gupta and Lin, Musiker, and Nakanishi \cite{gupta2019}. Their polynomiality theorem and separation formula assume
\[
 h_{i,s}(1)>0 \qquad (1\le i\le n,\ 1\le s<r_i).
\]
The calculation in \cite[Example~3.17]{FPY24} suggests that this assumption can be removed. Thus polynomiality for arbitrary prescribed reciprocal coefficients remains a question beyond that theorem. We resolve this question for the normalization specified below by constructing the \(F\)-polynomials over formal exchange coefficients before evaluating them.

\subsection*{Setting and main results}

Fix positive integers \(r_1,\ldots,r_n.\) Put
\[
 R_{\mathbf z}^{\mathrm{cl}}
 :=\frac{\QQ[z_{i,s}:1\le i\le n,\ 1\le s<r_i]}
          {(z_{i,s}-z_{i,r_i-s})},
 \qquad
 R_{\mathbf z}:=\QQ[t^{\pm1}]\otimes_{\QQ}
                         R_{\mathbf z}^{\mathrm{cl}},
\]
and set \(P_i(Y)=\sum_{s=0}^{r_i}z_{i,s}Y^s,\) with \(z_{i,0}=z_{i,r_i}=1\). In the seed-theoretic results we assume
\begin{equation}\label{eq:intro-normalized-compatibility}
 \widetilde B^T\Lambda=[I_n\mid0].
\end{equation}
Here \(\widetilde B\) is an integral \(m\times n\) exchange matrix and \(\Lambda\) is an integral skew form on the exponent lattice \(M.\) Thus all compatibility integers are \(1\); the exchange degrees \(r_i\) remain arbitrary. All seeds use fixed lattices \(N\) and \(M=N^*,\) with mutation changing the bases and their coordinate matrices. The following statements summarize the main results.

\begin{enumerate}[label=\textup{(\arabic*)},leftmargin=2.2em]
\item \emph{Positive scattering diagrams and mutation.} We construct a consistent generalized quantum scattering diagram \(\mathfrak D_S\) over \(R_{\mathbf z}.\) A central lattice extension records the roots of the exchange polynomials as independent degrees. The positive cylinder projection theorem produces a positive cylindrical representative in the lifted space, whose projection preserves truncations, path products, and consistency (Theorem~\ref{thm:cylinder-projection}). After splitting \(P_i(Y)=\prod_j(1+q_{ij}Y),\) this gives positive plethystic wall factorizations (Theorem~\ref{thm:seed-positive}).

The regular classical specialization of the wall-crossing actions
recovers the generalized cluster scattering diagram. Every essential
oriented wall survives this specialization. We use this fact and
control of the mutated degree monoids to prove
\[
 T_{k,S}(\mathfrak D_S^+)
 \sim\mathfrak D_{\mu_k(S)}
\]
(Theorem~\ref{thm:mutation}).

\item \emph{Cluster theta functions and root positivity.}
The wall-crossing actions define a rational scattering atlas and an algebra \(\mathcal A_q^{\mathrm{sc}}(S)\) over \(R_{\mathbf z}.\) Its reachable chambers recover the generalized cluster fan, with the frozen directions factored out. For a reachable seed \(u,\) let \(g_{i;u}\in M\) be the corresponding chamber generator. Then
\[
 X_{i;u}^{\mathrm{sc}}
 =\vartheta^{t_0}_{g_{i;u},Q_0},
\]
where \(t_0\) is the initial seed and \(Q_0\) is a general point of its positive chamber. More generally, every normalized reachable cluster monomial is a theta function (Theorem~\ref{thm:universal-cluster-theta}).

These functions are finite Laurent polynomials in every reachable
chart. After adjoining the exchange-polynomial roots, each Laurent
coefficient admits an expression in \(\ZZ_{\ge0}[t^{\pm1},q_{ij}]\) (Theorem~\ref{thm:universal-Laurent-positivity}). This gives a geometric
proof of Laurentness together with positivity in the root variables.

\item \emph{Universal polynomial separation and coefficient evaluation.}
The initial coordinates \(\widetilde{\mathbf g}_{i;u}\in\ZZ^m\) of \(g_{i;u}\) agree with the extended \(g\)-vectors of \cite{FPY24}. There is a universal quantum polynomial
\(F_{i;u}^{\mathrm{univ}},\) with coefficients in \(R_{\mathbf z}\) and constant term \(1,\) such that
\[
 X_{i;u}^{\mathrm{sc}}
 =X_{t_0}(\widetilde{\mathbf g}_{i;u})
       F_{i;u}^{\mathrm{univ}}(\widehat Y),
 \qquad
 \widehat Y^\alpha=X_{t_0}(\widetilde B_{t_0}\alpha).
\]
Polynomiality means that only finitely many based monomials \(\widehat Y^\alpha,\) \(\alpha\in\ZZ_{\ge0}^n,\) occur (Theorem~\ref{cor:F-polynomial-wall-crossing}).

For prescribed reciprocal \(h_{i,s}(t)\in\ZZ[t^{\pm1}]\) with \(h_{i,0}=h_{i,r_i}=1\), evaluation \(\sigma_h:z_{i,s}\mapsto h_{i,s}(t)\) recovers the generalized quantum cluster algebra of \cite{BCDX18,BCDX23}, over \(\QQ[t^{\pm1}]\) with frozen variables inverted (Corollary~\ref{cor:BCDX-after-evaluation}). It also gives
\[
 X_{i;u}^{(h)}
 =X_{t_0,h}(\widetilde{\mathbf g}_{i;u})
       F_{i;u}^{(h)}(\widehat Y_h),
 \qquad
 F_{i;u}^{(h)}=\sigma_h(F_{i;u}^{\mathrm{univ}}).
\]
Proposition~\ref{prop:FPY-F-identification} identifies this polynomial with the rational \(F\)-element defined in \cite{FPY24}.
Consequently, under \eqref{eq:intro-normalized-compatibility},
polynomiality and separation hold for arbitrary prescribed reciprocal coefficients, without the condition \(h_{i,s}(1)>0\) (Corollary~\ref{cor:BCDX-separation-after-evaluation}).
\end{enumerate}

\subsection*{Method and relation to previous results}

The main construction is the descent of positivity from a central
lattice extension. Davison--Mandel's theorem \cite[Theorem~2.15]{Davison_Mandel21} already permits positive one-parameter BPS polynomials and a degenerate skew form. We apply it to an enlarged lattice in which each root factor has its own degree. To descend the resulting positive walls, we construct a representative that is simultaneously positive and cylindrical. At each order there are finitely many relevant walls. Rational subdivisions of their base hyperplanes and uniqueness of the local BPS coefficients then give finitely many positive cylinders. Retaining the full factors inserted at earlier stages makes these constructions compatible in the inverse limit. This projection step supplies the geometric control needed for independent central coefficients.

Mutation follows the classical-limit strategy of \cite[Appendix~A]{Davison_Mandel21}. Reciprocal root relations require a separate noncancellation argument, and the source and target positive
completions must be compared at finite order. We specialize the regular adjoint actions, prove that essential oriented walls remain nontrivial, and use classical generalized mutation to locate the transformed quantum degrees. A quantum path-product identity proves consistency. The identification of cluster variables with theta functions further uses a comparison between formal wall crossing and rational chart pullbacks. Together these steps extend the quantum scattering realization to the higher-degree exchange relations considered here.

The algebraic consequences complement the results of \cite{BCDX18,BCDX23,FPY24}. The Laurent phenomenon for prescribed coefficients was proved in \cite{BCDX23}; our construction derives it
from theta functions over formal coefficients and also gives root
positivity. The Fock--Goncharov decomposition and the quantum Gupta product are established in \cite{FPY24}. Here their non-tropical factors are realized as wall-crossing actions, and polynomiality is proved before evaluation. This order of construction removes the classical positivity assumption in \cite[Theorem~3.11(2) and
Corollary~3.12]{FPY24} within \eqref{eq:intro-normalized-compatibility}.

The positivity conclusion concerns expressions in the roots of the
exchange polynomials. It is weaker than coefficientwise positivity in
the elementary exchange coefficients \(z_{i,s}.\) Arbitrary prescribed
coefficients may cause cancellation on evaluation; the finite
polynomial separation formula nevertheless survives. The comparison with BCDX uses \(t=q^{1/2},\) exchange degrees \(d_i^{\mathrm{BCDX}}=r_i,\) and exchange matrix \(\widetilde B^{\mathrm{BCDX}}=\widetilde B R,\) where \(R=\operatorname{diag}(r_1,\ldots,r_n).\)

\paragraph{Organization.}
Sections~\ref{sec:prelim}--\ref{sec:positivity} develop the scattering
construction and positive cylinder projection.
Sections~\ref{sec:seed-diagrams} and \ref{sec:mutation} treat seed
diagrams, classical specialization, and mutation.
Section~\ref{sec:scattering-algebra} constructs the scattering atlas
and proves the cluster theta and Laurent theorems.
Section~\ref{sec:coefficient-evaluation} compares the evaluated algebra
with the prescribed-coefficient construction.
Section~\ref{sec:fg-wall} establishes polynomial separation and its
wall-crossing interpretation, and gives a rank-two example.
Appendix~\ref{app:hypothesis-audit} records the proof dependencies.

\section{Quantum scattering diagrams}\label{sec:prelim}
In this section, we recall quantum scattering diagrams and the positivity theorem of Davison--Mandel \cite{Davison_Mandel21}. We also fix the coefficient substitutions used in the central lift.

\subsection{Quantum torus Lie algebras}

Let $L$ be a lattice of finite rank, let $L_0\subset L$ be a sublattice, and let $\omega:L\times L\to\ZZ$ be a skew-symmetric bilinear form. Fix a strongly convex rational polyhedral cone $\sigma\subset L_\RR$ and put $P=L_0\cap\sigma$ and $L_0^+=P\setminus\{0\}.$ The monoid $P$ is sharp and finitely generated. Choose an additive function
\begin{equation}\label{eq:background-order-function}
  \ord:P\longrightarrow\ZZ_{\geq0},
  \qquad \ord^{-1}(0)=\{0\}.
\end{equation}
Put $\QQ_t=\QQ[t^{\pm1}]$ and $(a)_t=t^a-t^{-a},$ following the Laurent coefficient convention of \cite[Section~2.1.1]{Davison_Mandel21}. The quantum torus $\QQ_t^\omega[L]$ has basis $\{x^m:m\in L\}$ and multiplication
\[
 x^a x^b=t^{\omega(a,b)}x^{a+b},
 \qquad [x^a,x^b]=(\omega(a,b))_t x^{a+b}.
\]

For $v\in L_0^+,$ let $|v|_0$ be its primitive index in $L_0.$ Following the normalized construction in \cite[(10)]{Davison_Mandel21}, define
\begin{equation}\label{eq:normalized-quantum-torus-Lie}
 \widehat x^v:=\frac{x^v}{(|v|_0)_t},
 \qquad
 \g:=\operatorname{Lie}_{\QQ_t}
       \langle\widehat x^v:v\in L_0^+\rangle
       \subseteq\bigoplus_{v\in L_0^+}\QQ(t)x^v.
\end{equation}
Thus $\g$ is the $\QQ_t$-Lie subalgebra generated by the displayed elements. Its homogeneous components and order filtration are
\[
 \g_v:=\g\cap \QQ(t)x^v,
 \qquad
 \g^{\geq d}:=\bigoplus_{\ord(v)\geq d}\g_v.
\]
They satisfy $[\g_u,\g_v]\subseteq\g_{u+v}$ and $[\g^{\geq a},\g^{\geq b}]\subseteq \g^{\geq a+b}.$ Define the completed Lie algebra and its exponential group by
\begin{equation}\label{eq:background-completed-Lie-group}
 \ghat:=\varprojlim_N\g/\g^{\geq N+1},
 \qquad
 \Ghat:=\varprojlim_N
 \exp\bigl(\g/\g^{\geq N+1}\bigr),
\end{equation}
with multiplication given by the Baker--Campbell--Hausdorff formula.

The normalized adjoint action preserves Laurent coefficients:
\begin{equation}\label{eq:normalized-Lie-adjoint}
 [\widehat x^v,x^m]
 =\frac{(\omega(v,m))_t}{(|v|_0)_t}x^{v+m}
 \in \QQ_t^\omega[L].
\end{equation}
Indeed, $v/|v|_0\in L_0$ implies $|v|_0\mid\omega(v,m).$ Put $\QQ_t [P]=\bigoplus_{p\in P}\QQ_tx^p$ and
\[
 \QQ_t [P]^{\geq d}
 :=\bigoplus_{\ord(p)\geq d}\QQ_tx^p,
 \qquad
 \QQ_t [\![P]\!]
 :=\varprojlim_N\QQ_t [P]/\QQ_t [P]^{\geq N+1}.
\]
The group $\Ghat$ acts continuously on each translated module \(\prod_{a\in P}\QQ_t x^{m+a},\) filtered by \(\ord(a).\) For rational comparisons we use the ring of based formal sums whose support is contained in a finite union of translates of \(P.\) Multiplication is coefficientwise finite because \(P\) has finite order sublevel sets. The finite quantum torus embeds in this ring. This is not an ideal-adic completion of the Laurent torus, in which the positive monomials are already units.

For a commutative $\QQ_t$-algebra $A,$ write $\g_A=\g\otimes_{\QQ_t}A.$ Its completed Lie algebra and group are formed from $(\g/\g^{\geq N+1})\otimes_{\QQ_t}A$ before taking the inverse limit.
In particular,
\begin{equation}\label{eq:normalized-Lie-scalar-extension}
 \g_{\QQ(t)}=\g\otimes_{\QQ_t}\QQ(t)
 =\bigoplus_{v\in L_0^+}\QQ(t)x^v.
\end{equation}
The last equality follows because every $(|v|_0)_t$ is invertible in $\QQ(t).$ We retain the normalized form when applying positivity or classical specialization and indicate coefficient extensions explicitly.

\begin{remark}[Primitive indices]
\label{rem:DM-normalized-Lie-form}
Our $|v|_0$ is distinct from $\ord(v)$ and is taken in the degree lattice $L_0.$ If $|v|_L$ is the ambient index used in \cite{Davison_Mandel21}, then
\[
 \frac{x^v}{(|v|_0)_t}
 =\frac{(|v|_L)_t}{(|v|_0)_t}\,
   \frac{x^v}{(|v|_L)_t},
 \qquad
 \frac{(|v|_L)_t}{(|v|_0)_t}\in \QQ_t.
\]
Thus our normalized form is contained in the ambient-index form;
they agree when the indices coincide, in particular when $L_0$ is
saturated in $L.$ Both have the scalar extension \eqref{eq:normalized-Lie-scalar-extension}.
\end{remark}

\subsection{Walls and scattering diagrams}

For a $v\in L_0^+$ with $\omega(v,-)\ne0,$ write
\[
  v^{\omega\perp}:=\{m\in L_\RR:\omega(v,m)=0\}.
\]
A \emph{wall} is a pair $(\mathfrak d,f_\mathfrak d),$ where $v$ is primitive in $L_0,$ $\mathfrak d$ is a closed convex rational polyhedral affine cone of codimension one parallel to $v^{\omega\perp},$ and
\[
 f_\mathfrak d\in\exp\!\left(\prod_{k\ge1}\g_{kv}\right).
\]
The degree $v$ is also called BPS degree. Its \emph{degree ray} is $\RR_{\ge0}v;$ the geometric direction in \cite{Davison_Mandel21} is $-v$. It is \emph{incoming} if $\mathfrak d+\RR_{\ge0}v=\mathfrak d,$ and \emph{outgoing} otherwise. For diagrams over a coefficient algebra $A$, replace $\g$ and $\Ghat$ by their indicated scalar extensions throughout.

\begin{definition}[Coincident wall packets in the degenerate case]
\label{def:coincident-wall-packet}
Write $\omega_1(v)=\omega(v,-).$ Fix an oriented ray $\rho\subset\omega_1(L_\RR)$ and a codimension-one rational polyhedral affine cone $\mathfrak d$ parallel to $\ker\rho.$ A \emph{coincident wall packet} on $(\mathfrak d,\rho)$ is an orderwise finite family
\[
 \bigl\{(\mathfrak d,f_{\mathfrak d,v})\bigr\}_v,
 \qquad
 f_{\mathfrak d,v}\in
 \exp\!\left(\prod_{\ell\geq1}\g_{\ell v}\right),
 \qquad
 \RR_{>0}\omega_1(v)=\rho .
\]
It may be displayed as one wall carrying $\prod_vf_{\mathfrak d,v},$ since $\omega(v,v')=0$ for its constituent degrees. The actual degrees are retained, and packets on opposite rays are kept separate.
\end{definition}

\begin{remark}
Packets abbreviate ordinary walls whose degrees may differ by radical vectors. BPS uniqueness is always indexed by the actual degree before coefficient specialization.
\end{remark}

A scattering diagram $\mathfrak D$ is a collection of walls such that, modulo $G^{\geq N+1}$ for every $N,$ only finitely many walls remain nontrivial.  
Denote $\supp(\mathfrak D):=\cup_{\mathfrak d\in \mathfrak D}\mathfrak d,$ and $$\Joints(\mathfrak D):=\bigcup\limits_{\mathfrak d\in\mathfrak D}\partial \mathfrak d\cup\bigcup\limits_{\substack{\mathfrak d_1,\mathfrak d_2\in \mathfrak D\\ \codim(\mathfrak d_1\cap\mathfrak d_2)=2}}(\mathfrak d_1\cap\mathfrak d_2).$$
A path $\gamma$ is admissible if its endpoints lie off the walls, it avoids the joints, and it crosses each wall of $\mathfrak D$ truncation transversely and at distinct times.

Suppose that, modulo $G^{\ge N+1},$ $\gamma$ crosses the nontrivial walls $(\mathfrak d_i,f_i)$ at increasing times $t_1\le\cdots\le t_s.$ If $-v_i$ is the direction of that wall, put
\[
 \epsilon_i:=\sgn\omega(v_i,-\dot\gamma(t_i)).
\]
With the convention that later crossings are multiplied on the left, the order-$N$ path product is
\begin{equation}\label{eq:background-truncated-path-product}
 \PathProd^{\leq N}_{\gamma,\mathfrak D}
 :=f_s^{\epsilon_s}\cdots f_1^{\epsilon_1}
 \pmod {G^{\geq N+1}}.
\end{equation}
Taking the inverse limit defines the path product $\PathProd_{\gamma,\mathfrak D}=\lim\limits_{\overset{\leftarrow}{N}}\PathProd^{\leq N}_{\gamma,\mathfrak D}\in\Ghat,$ its \emph{wall-crossing automorphism} is
\[
  \WC_{\gamma,\mathfrak D}:=\Ad_{\PathProd_{\gamma,\mathfrak D}}.
\]
Scattering diagrams are \emph{equivalent} if their finite-order path products
agree on common admissible paths. A scattering diagram is \emph{consistent} if
every admissible closed loop has trivial path product.

\begin{definition}[Reduction]
Let $\mathfrak D$ be a scattering diagram. For each $N,$ truncate $\mathfrak D$ modulo terms of order $>N,$ and take a finite common rational subdivision of its wall supports. Replace all walls
\[
  (\mathfrak d,f_1),\ldots,(\mathfrak d,f_s)
\]
with the same support and degree ray by the single wall
\[
  (\mathfrak d,f_1\cdots f_s),
\]
and delete this wall if its wall element is the identity modulo
terms of order $>N.$

The resulting diagram is called a reduction of
$\mathfrak D^{\le N}.$ Reduction of $\mathfrak D,$ denoted by $\mathrm{red}(\mathfrak D)$ is understood orderwise, with subdivisions compared by common refinement.
\end{definition}
Subdivision, merging and deletion preserve path products on common admissible paths. Thus every reduction is equivent to original scattering diagram. We write $\operatorname{red}(\mathfrak D)$ for this orderwise wall data, or for a representative when one is chose. The term reduced always refers to these orderwise data. Reducing and truncating give the same wall elements after common refinement and deletion of identities; no unique subdivision is asserted. The essential oriented support records the general points of the surviving wall pieces together with their rays $R_{>0}\omega_1(v).$ We use oriented wall germs only to compare these local data. 

The following theorem is a fundamental result in the theory of scattering diagrams. The case of quantum torus Lie algebras was proved in \cite{Kontsevich_Soibelman14} (see \cite{GHKK} for a review of this argument in the cluster setup). 

\begin{theorem}[{\cite[Theorem~2.13]{Davison_Mandel21}}]
\label{thm:KS-uniqueness}
Every finite incoming scattering diagram $\mathfrak D_{\mathrm{in}}$ admits a consistent completion $\Scat(\mathfrak D_{\mathrm{in}})$ obtained by adding only outgoing walls. This completion is unique up to equivalence.
\end{theorem}
All subsequent incoming data are supported on linear hyperplanes. The homogeneous construction then gives representatives with conical supports at the origin; for these, incoming is equivalent to $v\in\mathfrak d.$
\subsection{Broken lines and quantum theta functions}
\label{subsec:broken-lines}
Let $\mathfrak{D}$ be a consistent scattering diagram on $L_{\mathbb{R}}$ with respect to $\mathfrak{g}$ over the base field $\mathbb{K}.$ Fix $p\in L$ and a point $Q\in L_{\mathbb{R}}\setminus \operatorname{Supp}(\mathfrak{D}).$ A broken line $\gamma$ with respect to $\mathfrak{D}$ ending at $(p,Q)$ consists of the following data: a continuous map $\gamma\colon (-\infty,0]\to L_{\mathbb{R}}\setminus \operatorname{Joints}(\mathfrak{D}),$ a partition $-\infty < \tau_0 \leq \tau_1 \leq \ldots \leq \tau_{\ell} = 0,$ and for each $i=0,\ldots,\ell$ an associated monomial $c_ix^{v_i} \in \mathbb{K}_t[L]$ (with $c_i\in \mathbb{K}[t^{\pm 1}],$ $v_i\in L$), satisfying:
\begin{enumerate}[label=(\roman*), noitemsep]
\item $\gamma(0)=Q.$
\item For $i=1\ldots, \ell,$ $\gamma'(\tau)=-v_i$ for all $\tau\in (\tau_{i-1},\tau_{i}).$ Similarly, $\gamma'(\tau)=-v_0$ for all $\tau\in (-\infty,\tau_0).$
\item $c_0=1$ and $v_0=p.$
\item For $i=0,\ldots,\ell-1,$ $\gamma(\tau_i)\in \operatorname{Supp}(\mathfrak{D}).$ Let
\begin{align}\label{fi}
f_i\coloneqq \prod_{\substack{(\mathfrak{d},f_{\mathfrak{d}})\in \mathfrak{D} \\ \mathfrak{d}\ni \gamma(\tau_i)}} f_{\mathfrak{d}}^{\operatorname{sign}\,\omega(v_{\mathfrak{d}},v_i)} \in \hat{G}.
\end{align}
That is, $f_i$ is the limit as $\epsilon\to 0$ of $\theta_{\gamma|_{(\tau_i-\epsilon,\tau_i+\epsilon)},\mathfrak{D}}$ (using a smoothing of $\gamma$). Then $c_{i+1}x^{v_{i+1}}$ is a homogeneous term (with respect to the $L$-grading) of $\operatorname{Ad}_{f_i}(c_ix^{v_i})$ and is not equal to $c_ix^{v_i}$ (if the moment $\tau_i$ were omitted from the data, one would obtain the term $c_ix^{v_i}$).
\end{enumerate}

In the sequel we assume that every wall of $\mathfrak{D}$ has the origin as its vertex, and we only consider endpoints $Q$ that do not lie on any rational hyperplane, i.e., $Q\notin v^{\perp}$ for any $v\in L^*.$

For \(p\in L,\) the theta function associated with \(p\) and \(Q\) is
\[
  \vartheta_{p,Q}
  =
  \sum_{\gamma}c_\gamma x^{m_\gamma},
\]
where the sum runs over broken lines with initial exponent \(p\) and endpoint \(Q,\) and \(c_\gamma x^{m_\gamma}\) is the final monomial of \(\gamma.\) The theta function is well-defined and has form $\vartheta_{p,Q}=x^p(1+\sum\limits_{v\in L_0^+}a_{v}x^{v})$ by \cite[Proposition 3.1]{Davison_Mandel21}. Given two point $Q_1,Q_2\in L_{\RR}/\supp(\mathfrak D)$ and $p\in L_{\RR}$ the theta function $\vartheta_{p,Q_1},\vartheta_{p,Q_2}$ satisfy 
\[\vartheta_{p,Q_2}=\Ad_{\PathProd_{\gamma,\mathfrak D}}( \vartheta_{p,Q_1}),\]
where $\gamma$ is a smooth admissible path in $L_{\RR}\setminus \Joints(\mathfrak D)$ form $Q_1$ to $Q_2$ by \cite[Theorem 2.14]{Mandel21}. 
\subsection{Quantum dilogarithms and plethystic substitutions}
For $v\in L_0^+,$ the quantum dilogarithm is
\[
  \Psi_t(x^v):=\exp\left(\sum_{k\ge1}\frac{(-1)^{k-1}}{k(t^k-t^{-k})}x^{kv}\right).
\]
The standard product formula is \(\Psi_{t}(x^v)=\prod_{k\ge1}\frac{1}{1+t^{2k-1}x^{v}}.\)
We also use the plethystic notation 
\(
  \E(-x^v)=\Psi_t(x^v),
\) following \cite{Davison_Mandel21}.
To include central variables, put 
\[A_s:= \mathbb Z[t^{\pm 1},y_1,\dots,y_s]\]
We use the following explicit substitutions, with distinguished variable \(-t,\) as in \cite[Section~2.4.2]{Davison_Mandel21}.
For each $\ell\ge1,$ define the $\ZZ$-algebra endomorphism $\psi_\ell:A_s\to A_s$ by
\begin{equation}\label{eq:explicit-plethystic-substitution}
 \psi_\ell(t)=(-1)^{\ell-1}t^\ell,
 \qquad \psi_\ell(y_a)=y_a^\ell.
\end{equation}
Thus, for $p=\sum_{a,\alpha}c_{a,\alpha}t^a\mathbf y^\alpha,$
\begin{equation}\label{eq:explicit-plethystic-monomials}
 \psi_\ell(p)=\sum_{a,\alpha}c_{a,\alpha}
       (-1)^{a(\ell-1)}t^{a\ell}\mathbf y^{\ell\alpha}.
\end{equation}
Define
\begin{equation}\label{eq:plethystic-E-with-coefficients}
 \E_{A_s}(-p x^v) :=\exp\left(\sum_{\ell\ge1}
       \frac{(-1)^{\ell-1}\psi_\ell(p)x^{\ell v}}
            {\ell(t^\ell-t^{-\ell})}\right).
\end{equation}
This is the $E$-operator of \cite[(24)--(26)]{Davison_Mandel21}, with distinguished
variable $-t.$ Put $C_s:=\QQ(t)[y_1,\ldots,y_s].$ The series is defined over $C_s:$ each finite-order quotient involves only finitely many terms, and monomials on one degree ray commute.  Its logarithm belongs to the normalized Lie form extended by $\mathbf y,$ since, for $d=|v|_0,$
\[
 \frac{x^{\ell v}}{\ell(t^\ell-t^{-\ell})}
 =\frac{[d]_{t^\ell}}{\ell}\widehat x^{\ell v},
 \qquad
 [d]_{t^\ell}:=\frac{t^{d\ell}-t^{-d\ell}}
                         {t^\ell-t^{-\ell}}\in\ZZ[t^{\pm\ell}].
\]
Additivity gives

\[
 \E_{A_s}(-(p+p')x^v)
 =\E_{A_s}(-p x^v)\E_{A_s}(-p'x^v).
\]
The substitution $\psi_\ell(p)$ is not the ordinary power $p^\ell;$
for example, $\psi_2(y_1+y_2)=y_1^2+y_2^2.$ If $c=\omega(v,m),$ quantum-torus multiplication gives
\[
 \Ad_{\Psi_t(x^v)}(x^m)
 =x^m\Psi_t(t^{2c}x^v)\Psi_t(x^v)^{-1}.
\]
Hence, for $c\ne0$ and $\eta\in\{\pm1\},$
\[
 \Ad_{\Psi_t(x^v)^\eta}(x^m)
 =x^m\prod_{j=1}^{|c|}
   (1+t^{\sgn(c)(2j-1)}x^v)^{\eta\sgn(c)}.
\]
If $c=0,$ the action fixes $x^m.$ In a path product, $\eta$ is the
crossing sign in \eqref{eq:background-truncated-path-product}.
For $X=\mathbf y^\alpha x^v$ and $a\in\ZZ,$ the same definition gives
\begin{equation}\label{eq:odd-plethystic-factor}
 \E_{A_s}(-t^aX)
 =\Psi_t\bigl((-1)^a t^aX\bigr)^{(-1)^a}.
\end{equation}
In particular,
\[
 \E_{A_s}(-t^{2b}X)=\Psi_t(t^{2b}X),\qquad
 \E_{A_s}(-t^{2b+1}X)=\Psi_t(-t^{2b+1}X)^{-1}.
\]
The latter also has form $\prod_{r \ge 0} \Psi_{t^{2^r}} \bigl( t^{2^r(2b+1)} X^{2^r} \bigr)$ by \cite{Davison_Mandel21}. The signed substitution in
\eqref{eq:explicit-plethystic-substitution} is essential.

\begin{lemma}[Plethystic base change]
\label{lem:plethystic-base-change}
Let $B$ be a commutative $\QQ$-algebra and let $\rho:C_s\to B$ be a
unital homomorphism. Put $t_B=\rho(t).$  Coefficientwise evaluation
induces a continuous map between the completed quantum tori over
$C_s$ and $B$, and satisfies
\[
 \widehat\rho\bigl(\E_{A_s}(-p x^v)\bigr)
 =\exp\left(\sum_{\ell\ge1}
   \frac{(-1)^{\ell-1}\rho(\psi_\ell(p))x^{\ell v}}
        {\ell(t_B^\ell-t_B^{-\ell})}\right).
\]
It commutes with products, inverses, adjoint actions and path products, and preserves consistency.
\end{lemma}

\begin{proof}
On each finite-order quotient, these operations are finite algebraic
expressions over $C_s.$ Apply $\rho$ and pass to the inverse limit.
\end{proof}

No substitutions on $B$ are needed. The classical limit $t=1$ is
taken separately on regular adjoint actions.

For the positivity theorem, we also allow a rational skew form
$\omega:L\times L\to\frac1D\ZZ$ satisfying $\omega(L_0,L)\subseteq\ZZ.$  The normalized Lie form is defined over $\QQ[t^{\pm1}],$ \eqref{eq:normalized-Lie-adjoint} still has Laurent coefficients. The ambient torus is formed over $\QQ[t^{\pm1/D}];$ field-valued path groups use $\g\otimes_{\QQ_t}\QQ(t^{1/D}).$
\begin{theorem}[Davison--Mandel positivity]
\label{thm:DM}Let the finite incoming diagram consist of walls
\[
 \left(v_i^{\omega\perp},\E(-p_i(t)x^{v_i})\right),
 \qquad
 v_i\in L_0^+\setminus\ker\omega_1,\quad
 p_i(t)\in\ZZ_{\ge0}[t^{\pm1}].
\]
Then its consistent completion has a representative with precisely
these incoming walls and with every added wall element of the form
\[
 \E(-p_{\mathfrak d}(t)x^{v_{\mathfrak d}}),
 \qquad
 v_{\mathfrak d}\in L_0^+\setminus\ker\omega_1,\quad
 p_{\mathfrak d}(t)\in\ZZ_{\ge0}[t^{\pm1}].
\]
\end{theorem}
\begin{remark}[Relative integrality]\label{rem:DM-relative-integrality}
We also call $p_{\mathfrak d}(t)$ the BPS polynomial. Only $\omega(L_0,L)\subseteq\ZZ$ is required; the ambient form may be rational and degenerate.  Since $v/|v|\in L_0,$ this condition gives $|v|\mid\omega(v,m)$ for $m\in L$ and ensures Laurent coefficients in the normalized adjoint action. For the normalized seeds of Theorem~\ref{thm:seed-positive}, the integral form is simply \(\Lambda\), and compatibility gives \(\Lambda(v_i,m)=\langle e_i,m\rangle.\)
\end{remark}
Coefficientwise positivity in the independent central variables $\mathbf y$ is proved by the central lift in
Corollary~\ref{cor:plethystic-coefficient-positivity}.
\section{Generalized quantum cluster scattering diagrams}\label{sec:generalized-diagram}
Fix positive integers
\[
  \mathbf r=(r_1,\ldots,r_n)
\]
and formal reciprocal coefficients $\mathbf z:=\{z_{i,s}\}_{1\le i\le n,0\le s\le r_i}$ with $z_{i,0}=z_{i,r_i}=1$ for all $1\le i\le n.$ Keeping the distinction between coefficient rings and field extensions used in
\cite[Section~2.1.1]{Davison_Mandel21}, set
\[
 R_{\mathbf z}^{\mathrm{cl}}
 :=
 \frac{\QQ[ z_{i,s}:1\le i\le n,\ 1\le s<r_i]}
 {\bigl( z_{i,s}- z_{i,r_i-s}:
          1\le i\le n,\ 1\le s<r_i\bigr)},
 \qquad
 R_{\mathbf z}:=\QQ[t^{\pm1}]\otimes_{\QQ}R_{\mathbf z}^{\mathrm{cl}},
\]
so \(R_{\mathbf z}^{\mathrm{cl}}\) is equivalently the polynomial ring in the independent variables \( z_{i,s}\) with \(1\leq s\leq\lfloor r_i/2\rfloor.\) Put
\[
 \mathbb F_{\mathbf z}^{\mathrm{cl}}
 :=\operatorname{Frac}(R_{\mathbf z}^{\mathrm{cl}}),
 \qquad
 \mathbb F_{\mathbf z}:=\operatorname{Frac}(R_{\mathbf z})
 =\mathbb F_{\mathbf z}^{\mathrm{cl}}(t), 
\]
The coefficient specialization $R_{\mathbf z}\to R_{\mathbf z}^{\mathrm{cl}}$ is $t\mapsto1,$ with the $z_{i,s}$ fixed.
The coefficient fraction field is used only after an explicitly indicated scalar extension.
Define
\[
 P_i(y):=\sum_{s=0}^{r_i} z_{i,s}y^s\in R_{\mathbf z}[y].
\]
For \(\ell\geq1,\) let \(s_{i,\ell}\in
R_{\mathbf z}^{\mathrm{cl}}\) be the Newton power sum determined by
\begin{equation}\label{eq:formal-Newton-power-sums}
 s_{i,\ell}
 +\sum_{a=1}^{\ell-1}
   (-1)^a z_{i,a}s_{i,\ell-a}
 +(-1)^\ell\ell z_{i,\ell}=0,
\end{equation}
where \( z_{i,a}=0\) for \(a>r_i.\) Let \(K_{\mathbf z}/\mathbb F_{\mathbf z}\) be a splitting field for the \(P_i,\) and write
\[
 P_i(y)=\prod_{j=1}^{r_i}(1+q_{ij}y)
 \qquad\text{in }K_{\mathbf z}[y].
\]
Then Newton's identities give
\begin{equation}\label{eq:Newton-power-sums-after-splitting}
 s_{i,\ell}=\sum_{j=1}^{r_i}q_{ij}^{\ell}.
\end{equation}

Let $L$ be a lattice with basis
\[
  e_1,\ldots,e_m,
  \qquad m\ge n,
\]
and put
\[
 P:=\sum_{i=1}^n\ZZ_{\geq0}e_i,
 \qquad
 L_0^+:=P\setminus\{0\},
 \qquad
 \ord\!\left(\sum_i a_ie_i\right):=\sum_i a_i.
\]
Assume $\omega:L\times L\to\ZZ$ is skew-symmetric, and set $L_0=\bigoplus_{i=1}^n\ZZ e_i.$ For the normalized Lie algebra $\g$ of \eqref{eq:normalized-quantum-torus-Lie}, put
\[
 \g_{R_{\mathbf z}}
 :=\g\otimes_{\QQ_t}R_{\mathbf z},
 \qquad
 \g_{\mathbb F_{\mathbf z}}
 :=\g_{R_{\mathbf z}}\otimes_{R_{\mathbf z}}\mathbb F_{\mathbf z}.
\]
The first algebra retains the normalized generators of \eqref{eq:normalized-quantum-torus-Lie}; only after extending scalars does one have $\g_{\mathbb F_{\mathbf z}}=\bigoplus_{v\in L_0^+}\mathbb F_{\mathbf z}z^v.$
All Lie algebras, wall groups, and quantum torus modules in Sections~\ref{sec:generalized-diagram}--\ref{sec:positivity} are
completed for this coordinate-sum filtration, exactly as in
\eqref{eq:background-completed-Lie-group}. Thus every completed base change below means the inverse limit of the base changes of the finite order quotients.
Put
\[
 I_{\mathrm{act}}
 :=\{i\in\{1,\ldots,n\}:\omega(e_i,-)\ne0\}.
\]
The \emph{active-index convention} means that the scattering diagram is
formed only from the indices in \(I_{\mathrm{act}};\) all sums, positive
monoids, lifts, and incoming-wall lists in Sections \ref{sec:generalized-diagram}--\ref{sec:positivity} use this set. The set $I_{\mathrm{act}}$ indexes the incoming walls. The scattering Lie algebra remains graded by $L_0^+.$ For the given incoming data, the Kontsevich--Soibelman construction takes place in the graded Lie subalgebra supported on \[ P_{\mathrm{act}}\setminus\{0\}, \qquad P_{\mathrm{act}} :=\sum_{i\in I_{\mathrm{act}}}\ZZ_{\ge0}e_i. \] In particular, outgoing walls may have mixed degrees. To avoid
repeated restrictions of notation, throughout those sections we relabel the active indices and assume
\begin{equation}\label{eq:active-scattering-directions}
 \omega(e_i,-)\ne0\qquad(1\leq i\leq n).
\end{equation}
If an index violates \eqref{eq:active-scattering-directions}, its
homogeneous quantum-torus element is central, its adjoint wall action is the identity, and \(e_i^{\omega\perp}=L_\RR\) is not a codimension-one wall.  Such an index is retained, if desired, only as a central or frozen coefficient and is omitted from the active diagram. This operation does not delete a mutable direction from a compatible generalized seed: Lemma~\ref{lem:compatible-seed-indices-active} below
shows that every unfrozen seed index is automatically active.

\begin{definition}[Higher-degree quantum wall element]
For $1\le i\le n$, define
\begin{align}
  \Phi_i(x^{e_i})
  &:=\exp\left(\sum_{k\ge1}\frac{(-1)^{k-1}}{k(t^k-t^{-k})}
      s_{i,k}x^{k e_i}\right) \notag\\
  &=\prod_{j=1}^{r_i}\Psi_t(q_{ij}x^{e_i})
  \qquad\text{after base change to }K_{\mathbf z}.
  \label{eq:Phi-factor}
\end{align}
Since $|k e_i|_0=k,$ its logarithm is
\[
 \log\Phi_i(x^{e_i})
 =\sum_{k\ge1}\frac{(-1)^{k-1}}{k}s_{i,k}\widehat x^{k e_i}
 \in\widehat{\g}_{R_{\mathbf z}}.
\]
Thus $\Phi_i$ is a wall group element over $R_{\mathbf z}$ in the normalized Lie sense. Its adjoint action has Laurent coefficients, as the next lemma shows; its ordinary power-series coefficients may have denominators in $t.$
\end{definition}
Let $a\in L$ and put $\eps=\sgn\omega(e_i,a).$ Then
\[
  \Ad_{\Phi_i(x^{e_i})^\eps}(x^a)
  =x^a\prod_{s=1}^{|\omega(e_i,a)|}
  \left(1+\sum_{j=1}^{r_i}z_{i,j}t^{\eps j(2s-1)}x^{j e_i}\right).
\]

\begin{definition}[Generalized quantum scattering diagram]
\label{def:GQSD}
Under the active-index convention above, the initial generalized quantum scattering diagram is
\[
  \Din=\bigl\{(e_i^{\omega\perp},\Phi_i(x^{e_i})):
                    i\in I_{\mathrm{act}}\bigr\}.
\]
We define
\[
  \mathfrak D_{R_{\mathbf z}}:=\Scat(\Din)
\]
and call it the generalized quantum scattering diagram associated to
$(L,\omega,r,\mathbf z);$ we abbreviate it to \(\mathfrak D\) when the coefficient ring is clear.  In the relabelled convention \eqref{eq:active-scattering-directions}, \(I_{\mathrm{act}}=\{1,\ldots,n\}.\)
\end{definition}

\begin{lemma}[Orderwise coefficient base change]
\label{lem:reciprocal-coefficient-base-change}
The diagram \(\mathfrak D_{R_{\mathbf z}}\) is defined orderwise over \(R_{\mathbf z}.\) If \(A\) is a \(\QQ_t\)-algebra and \(\sigma:R_{\mathbf z}\to A\) is a \(\QQ_t\)-algebra homomorphism, define \(\widehat\otimes_\sigma A\) as the inverse limit of tensor products of the finite-order quotients. Then
\[
 \operatorname{red}\!\left(
  \mathfrak D_{R_{\mathbf z}}\widehat\otimes_{\sigma}A
 \right)
 \sim \Scat(\sigma(\Din)).
\]
Truncation, path products, and consistent completion commute with this base change.
\end{lemma}

\begin{proof}
The incoming logarithms belong to $\widehat{\g}_{R_{\mathbf z}}.$ At each order the Kontsevich--Soibelman construction uses homogeneous components, finite sums, Lie brackets, and rational BCH coefficients.
It therefore stays in this normalized Lie form, whose adjoint action
preserves the completed $R_{\mathbf z}$-quantum torus. Applying \(\sigma\) to the finite construction gives the completion of \(\sigma(\Din).\) Uniqueness at every order and passage to the inverse limit prove the statement.
\end{proof}

\section{Central lifting and cylinder projection}\label{sec:central-lift}

We construct a positive cylindrical representative of the lifted
scattering class by adding positive factors at successive finite orders,
and then project it to the base space.

\subsection{The central extension}

Let
\[
  U:=\bigoplus_{i=1}^n\bigoplus_{j=1}^{r_i}\ZZ u_{ij}
\]
and set
\[
  L^\dagger:=L\oplus U.
\]
Extend the skew form by
\[
  \omega^\dagger(\ell+u,\ell'+u'):=\omega(\ell,\ell'),
  \qquad
  \ell,\ell'\in L,
  \quad u,u'\in U.
\]
Thus $U\subset\ker\omega^\dagger.$ Define
\[
  f_j^i:=e_i+u_{ij}\in L^\dagger.
\]
The positive monoid in the lifted lattice is generated by the $f_j^i$:
\[
  L_0^\dagger:=\bigoplus_{i,j}\ZZ f_j^i,
\]
and
\[
  (L_0^\dagger)^+
  :=\left\{\sum_{i=1}^n\sum_{j=1}^{r_i}a_{ij}f_j^i:a_{ij}\in\ZZ_{\ge0}\right\}\setminus\{0\}.
\]
Let
\[
  \Pr:L_\RR^\dagger\longrightarrow L_\RR
\]
be the projection along $U_\RR.$

\begin{definition}[Order]
For
\[
  v=\sum_{i,j}a_{ij}f_j^i\in (L_0^\dagger)^+,
\]
define
\[
  \ord^\dagger(v):=\sum_{i,j}a_{ij}.
\]
For
\[
  w=\sum_{i=1}^n b_i e_i\in L_0^+,
\]
define
\[
  \ord(w):=\sum_i b_i.
\]
This is the order used in the completion.  It is the coordinate sum with respect to the chosen positive basis, not the primitive index of a lattice vector.
\end{definition}

\begin{lemma}[Order preservation]
\label{lem:order-preservation}
For every $v\in (L_0^\dagger)^+,$
\[
  \ord(\Pr(v))=\ord^\dagger(v).
\]
\end{lemma}
Let $\g^\dagger_{\mathrm{norm}}$ be the normalized Lie algebra for $(L^\dagger,L_0^\dagger,\omega^\dagger),$ and put $\g^\dagger:=\g^\dagger_{\mathrm{norm}}
\otimes_{\QQ[t^{\pm1}]}\QQ(t).$
The incoming factors below lie in $\exp(\widehat{\g^\dagger_{\mathrm{norm}}}).$
We apply positivity to this normalized form and use its scalar extension for the Lie-algebra comparisons below.

Define the lifted incoming diagram
\[
  \Ddag_{\mathrm{in}}
  :=\left\{\left((f_j^i)^{\omega^\dagger\perp},\E(-x^{f_j^i})\right):1\le i\le n,\,1\le j\le r_i\right\}
\]
and let
\[
  \Ddag:=\Scat(\Ddag_{\mathrm{in}}).
\]
\begin{proposition}[Verification of the Davison--Mandel hypotheses]
\label{prop:DM-hypotheses-central-lift}
The central lifted incoming diagram \(\Ddag_{\mathrm{in}}\) satisfies every hypothesis of Theorem~\ref{thm:DM}, even though \(\omega^\dagger\) is degenerate. More explicitly:
\begin{enumerate}[label=(\roman*)]
\item \(L^\dagger\) is a finite-rank lattice and \(\omega^\dagger\) is an integral skew-symmetric form;
\item the cone
\[
 \sigma^\dagger:=\sum_{i,j}\RR_{\geq0}f_j^i
 \subset L^\dagger_\RR
\]
is strongly convex and rational polyhedral, and \((L_0^\dagger\cap\sigma^\dagger)\setminus\{0\}
=(L_0^\dagger)^+;\)
\item the incoming diagram is finite, every \(f_j^i\) is primitive, belongs to \((L_0^\dagger)^+,\) and has \(\omega^\dagger(f_j^i,-)=\omega(e_i,-)\ne0;\)
\item every support is the full incoming hyperplane \((f_j^i)^{\omega^\dagger\perp},\) and every incoming BPS polynomial is \(1\in\ZZ_{\geq0}[t^{\pm1}];\)
\item the order completion is the monoid-adic completion required by Theorem~\ref{thm:DM}, and the coefficient substitutions use the distinguished-\(-t\) convention \eqref{eq:explicit-plethystic-substitution}.
\item the graded Lie algebra satisfies the Davison--Mandel compatibility condition
\[
 \omega^\dagger(a,b)=0
 \quad\Longrightarrow\quad
 [\mathfrak g_a^\dagger,\mathfrak g_b^\dagger]=0;
\]
the sublattice \(U\) is therefore central rather than an additional interacting sector.
\end{enumerate}
\end{proposition}

\begin{proof}
Only (ii) and the possible effect of degeneracy need verification.
Projection to \(U\) sends \(f_j^i=e_i+u_{ij}\) to the basis vector
\(u_{ij}.\) Hence the \(f_j^i\) are primitive and linearly independent, so \(\sigma^\dagger\) is a strongly convex simplicial rational cone and its lattice points in \(L_0^\dagger\) form exactly the stated free positive monoid. Moreover,
\[
 \omega^\dagger(f_j^i,-)=\omega(e_i,-)\ne0
\]
by \eqref{eq:active-scattering-directions}; thus no incoming degree becomes radical.  Theorem~\ref{thm:DM} does not require the ambient skew form to be nondegenerate, and
\[
 [x^a,x^b]=(\omega^\dagger(a,b))_t x^{a+b}
\]
gives (vi). All remaining assertions follow immediately from the
definitions of the incoming walls, the order completion, and
\(\E.\)
\end{proof}

\begin{lemma}[Radical total degrees are not generated]
\label{lem:radical-degrees-not-generated}
Let \(\mathfrak h\subset\mathfrak g^\dagger\) be the completed Lie subalgebra generated by the active incoming homogeneous elements \(x^{\ell f_j^i},\) for every \(\ell\geq1.\) If
\[
 0\ne v\in(L_0^\dagger)^+
 \cap\ker(\omega^\dagger_1),
 \qquad
 \omega^\dagger_1(v)=\omega^\dagger(v,-),
\]
then \(\mathfrak h_v=0.\) Consequently, the Kontsevich--Soibelman
completion of the central lift never creates an effective wall whose total degree is radical.
\end{lemma}

\begin{proof}
No generator has radical degree by
\eqref{eq:active-scattering-directions}.  Consider an iterated bracket of at least two generators with total degree \(v,\) and write its
outermost bracket as \([A_a,B_b],\) where \(a+b=v.\) Its scalar bracket
coefficient is a multiple of
\[
 \omega^\dagger(a,b)
 =\omega^\dagger(a,v-a)
 =\omega^\dagger(a,v)=0.
\]
Thus every such Lie monomial vanishes.  Completed homogeneous
components are inverse limits of the same finite-order calculation, so the conclusion persists in the completion. The incoming logarithms have only degrees \(\ell f_j^i,\) so they belong to \(\mathfrak h.\) Every correction wall is produced by brackets of these homogeneous terms; hence no radical total degree occurs.
\end{proof}
\subsection{The cylinder representative}
For a multi-index
\(\alpha=(a_{ij})\in\ZZ_{\geq0}^{\sum_i r_i},\) set
\[
 |\alpha|:=\sum_{i,j}a_{ij},\qquad
 u(\alpha):=\sum_{i,j}a_{ij}u_{ij},\qquad
 w(\alpha):=\sum_i\left(\sum_j a_{ij}\right)e_i,
\]
and write \(\mathbf y^\alpha:=\prod_{i,j}y_{ij}^{a_{ij}}.\)
Define
\[
 \overline{\g}
 :=\bigoplus_{w\in L_0^+}\overline{\g}_w,
 \qquad
 \overline{\g}_w
 :=\bigoplus_{\alpha:\,w(\alpha)=w}
       \QQ(t)\,\mathbf y^\alpha x^w.
\]
We equip this space with the bracket inherited from the quantum
torus on \(L\), with the variables \(y_{ij}\) central:
\begin{equation}\label{eq:balanced-coefficient-bracket}
 [\mathbf y^\alpha x^{w(\alpha)},
  \mathbf y^\beta x^{w(\beta)}]
 =
 \bigl(\omega(w(\alpha),w(\beta))\bigr)_t
 \mathbf y^{\alpha+\beta}x^{w(\alpha+\beta)}.
\end{equation}
Since \(w(\alpha+\beta)=w(\alpha)+w(\beta),\) this makes
\(\overline{\g}\) an \(L_0^+\)-graded Lie algebra satisfying
\[
 [\overline{\g}_w,\overline{\g}_{w'}]
 \subseteq\overline{\g}_{w+w'},
 \qquad
 \omega(w,w')=0
 \Longrightarrow
 [\overline{\g}_w,\overline{\g}_{w'}]=0.
\]
Scattering diagrams in \(L_\RR\) use this \(w\)-grading;
the finer \(\alpha\)-grading records the central multi-degree.

For \(w=\sum_i b_i e_i\), the condition \(w(\alpha)=w\) is
equivalent to \(\sum_j a_{ij}=b_i\) for every \(i\), so each
\(\overline{\g}_w\) is finite-dimensional. Moreover,
\[
 |\alpha|
 =\ord^\dagger\bigl(w(\alpha)+u(\alpha)\bigr)
 =\ord\bigl(w(\alpha)\bigr).
\]
Thus completion by \(|\alpha|\) agrees with the order completion
introduced above.
This is the balanced subalgebra of the coefficient
quantum torus; arbitrary monomials \(\mathbf y^\alpha x^w\) with \(w\ne w(\alpha)\) are not introduced.

\begin{lemma}[Central graded Lie-algebra identification]
\label{lem:central-graded-Lie-identification}
The assignment
\begin{equation}\label{eq:central-Lie-isomorphism}
 \iota:
 \mathbf y^\alpha x^{w(\alpha)}
 \longmapsto
 x^{w(\alpha)+u(\alpha)}
 =x^{\sum_{i,j}a_{ij}f_j^i}
\end{equation}
is an order-preserving isomorphism from \(\overline{\g}\) to the lifted positive quantum-torus Lie algebra \(\g^\dagger.\) It extends to the completed Lie algebras and theirpro-unipotent groups, and it carries \(E\)-operators to the corresponding lifted \(E\)-operators.
\end{lemma}

\begin{proof}
The vectors \(f_j^i=e_i+u_{ij}\) are linearly independent: their
\(U\)-coordinates recover all coefficients \(a_{ij}.\) Hence \eqref{eq:central-Lie-isomorphism} is a bijection on homogeneous degrees. Since \(U\subset\ker\omega^\dagger,\)
\[
 \omega^\dagger(w(\alpha)+u(\alpha),
                  w(\beta)+u(\beta))
 =\omega(w(\alpha),w(\beta)).
\]
Comparison with \eqref{eq:balanced-coefficient-bracket} proves that \(\iota\) preserves brackets. Lemma~\ref{lem:order-preservation}
shows that it preserves the filtration. It therefore extends to the inverse-limit completions and exponentiates to an isomorphism of the completed scattering groups. Finally, \(\psi_\ell(\mathbf y^\alpha)
=\mathbf y^{\ell\alpha},\) so applying
\(\iota\) to \eqref{eq:plethystic-E-with-coefficients} gives the same series as first applying \(\iota\) to the monomial and then taking the \(E\)-operator. Explicitly,
\[
 \iota\!\left(\mathbf y^{\ell\alpha}
               x^{\ell w(\alpha)}\right)
 =x^{\ell(w(\alpha)+u(\alpha))},
\]
which is precisely compatibility with the \(\ell\)-fold lifted degree.
\end{proof}
Under $\iota$, the transported normalized form retains the full lifted
index
\[
 \left|\sum_{i,j}a_{ij}f_j^i\right|_0=\gcd(a_{ij}),
\]
which need not equal $|w(\alpha)|_0$ after forgetting the central degree.

Write
\begin{equation}\label{eq:balanced-projection-map}
 \pi_{\mathrm{bal}}:=\iota^{-1},
 \qquad
 \pi_{\mathrm{bal}}\!
 \left(x^{w(\alpha)+u(\alpha)}\right)
 =\mathbf y^\alpha x^{w(\alpha)}.
\end{equation}
This is not a quotient which forgets the \(U\)-degree: it records that
degree as a central coefficient monomial. It is an order-preserving
isomorphism of completed Lie algebras and groups. If a lifted diagram \(\mathfrak C\) has cylindrical walls
\[
 (\overline{\mathfrak d}\oplus U_\RR,f_{\mathfrak d}),
\]
define its cylinder projection, before reduction, by
\begin{equation}\label{eq:cylinder-projection-functor}
 \Pi_{\mathrm{cyl}}(\mathfrak C)
 :=
 \left\{
  \bigl(\overline{\mathfrak d},
        \pi_{\mathrm{bal}}(f_{\mathfrak d})\bigr)
 \right\}.
\end{equation}
Coincident projected supports are retained until their commuting local wall products have been combined.

Let \(\overline{\mathfrak D}_{\mathrm{in}}\) be the diagram on \(L_\RR,\) over \(\overline{\g},\) with incoming walls
\[
 \left(e_i^{\omega\perp},
       \E(-y_{ij}x^{e_i})\right),
 \qquad 1\leq i\leq n,\quad 1\leq j\leq r_i,
\]
and put
\[
 \Db:=\Scat(\overline{\mathfrak D}_{\mathrm{in}}).
\]
\[
 \Db:=\Scat(\overline{\mathfrak D}_{\mathrm{in}}).
\]
Choose a conical representative of \(\Db\) containing
\(\overline{\mathfrak D}_{\mathrm{in}}\) and having only outgoing
remaining walls. Define its cylinder pullback by
\begin{equation}\label{eq:cylinder-pullback-definition}
 (\overline{\mathfrak d},
   \overline f_{\overline{\mathfrak d}})
 \longmapsto
 \left(\overline{\mathfrak d}\oplus U_\RR,
       \iota(\overline f_{\overline{\mathfrak d}})\right).
\end{equation}
The right-hand side is understood as a coincident wall packet in the
sense of Definition~\ref{def:coincident-wall-packet}. More precisely, decompose
\[
 \log\iota(\overline f_{\overline{\mathfrak d}})
 =\sum_{\lambda}h_{\lambda}
\]
by actual lifted degree rays \(\lambda.\) Replace the displayed packet
by the ordinary walls
\[
 \bigl(\overline{\mathfrak d}\oplus U_\RR,
       \exp(h_{\lambda})\bigr),\qquad h_{\lambda}\ne0.
\]
The components commute, and only finitely many are nonzero at each order. Denote the resulting diagram by \(\Pr^*\Db\).

\begin{lemma}[Cylinder pullback]
\label{lem:cylinder-pullback}
The diagram \(\Pr^*\Db\) is a well-defined consistent scattering diagram, its incoming diagram is \(\Ddag_{\mathrm{in}},\) and
\[
 \Pr^*\Db\sim\Ddag.
\]
In particular, the equivalence class of \(\Ddag\) contains a reduced representative all of whose wall supports are cylinders \(\overline{\mathfrak d}\oplus U_\RR.\)
\end{lemma}

\begin{proof}
Let \(\mathbf y^\alpha x^{w(\alpha)}\) be a homogeneous monomial in a wall group of \(\Db.\) By centrality of \(U,\)
\[
 (w(\alpha)+u(\alpha))^{\omega^\dagger\perp}
 =w(\alpha)^{\omega\perp}\oplus U_\RR.
\]
All degrees in one base wall group are positive multiples of its
primitive degree. Their lifted degrees therefore have the same
oriented effective covector ray and pair to zero under
\(\omega^\dagger.\) The components \(h_{\lambda}\) commute, and each
lies in the Lie algebra of its actual lifted degree ray. Thus every
constituent of the packet in \eqref{eq:cylinder-pullback-definition} is a legitimate wall. Lemma~\ref{lem:order-preservation} preserves order. Since there are
only finitely many base walls and lifted degrees at each order, the
pullback is finite at each order as well.

Let \(\gamma^\dagger\) be an admissible path for \(\Pr^*\Db.\) The walls and their joints are full inverse images under \(\Pr,\) so \(\gamma=\Pr\circ\gamma^\dagger\) is admissible. A lifted packet is crossed exactly when its base wall is crossed. If a constituent has lifted degree \(v=w+u,\) then
\[
 \omega^\dagger(v,-(\gamma^\dagger)'(s))
 =\omega(w,-\gamma'(s)),
\]
so its crossing sign agrees with that of the base wall. At each
crossing event, the product of the lifted constituent factors is the
image under \(\iota\) of the base crossing factor. The order of the
events is unchanged. Consequently,


\begin{equation}\label{eq:path-product-cylinder-pullback}
 \PathProd_{\gamma^\dagger,\Pr^*\Db}
 =
 \iota\!\left(\PathProd_{\gamma,\Db}\right).
\end{equation}
For a closed path the right-hand side is the identity; hence \(\Pr^*\Db\) is consistent.

For an incoming wall,
\[
 (f_j^i)^{\omega^\dagger\perp}
 =e_i^{\omega\perp}\oplus U_\RR
\]
and, by Lemma~\ref{lem:central-graded-Lie-identification},
\[
 \iota\bigl(\E(-y_{ij}x^{e_i})\bigr)
 =\E(-x^{e_i+u_{ij}})
 =\E(-x^{f_j^i}).
\]
Thus the incoming diagram is exactly \(\Ddag_{\mathrm{in}}.\) Moreover, for a lifted degree \(v=w+u,\)
\[
 v\in\overline{\mathfrak d}\oplus U_\RR
 \quad\Longleftrightarrow\quad w\in\overline{\mathfrak d}.
\]
Thus the other walls remain outgoing.  The uniqueness theorem for consistent completion now gives \(\Pr^*\Db\sim\Ddag.\) Finally, reducing the
pullback only deletes trivial walls and combines walls with identical
cylindrical support, so the reduced representative remains
cylindrical.
\end{proof}

\begin{proposition}[Cylinder truncation, local finiteness, and path products]
\label{prop:cylinder-projection-formal}
Let \(\mathfrak C\) be an orderwise locally finite conical cylindrical diagram whose wall degrees lie in \((L_0^\dagger)^+,\) and write
\(\mathfrak C^{\leq N}\) for its quotient modulo order \(>N.\)
Then:
\begin{enumerate}[label=(\roman*)]
\item \(\Pi_{\mathrm{cyl}}(\mathfrak C)\) is orderwise locally finite, and truncation commutes with projection:
\[
 \Pi_{\mathrm{cyl}}(\mathfrak C^{\leq N})
 =
 \Pi_{\mathrm{cyl}}(\mathfrak C)^{\leq N};
\]
\item for every admissible base path \(\gamma\) and every generic constant-\(U\) lift \(\widetilde\gamma,\)
\begin{equation}\label{eq:formal-projected-path-product}
 \PathProd_{\gamma,\Pi_{\mathrm{cyl}}(\mathfrak C)}
 =
 \pi_{\mathrm{bal}}\!
 \left(\PathProd_{\widetilde\gamma,\mathfrak C}\right);
\end{equation}
\item with cylindrical subdivisions, reduction commutes with projection up to equivalence:
\[
 \operatorname{red}\Pi_{\mathrm{cyl}}(\mathfrak C)
 \sim
 \Pi_{\mathrm{cyl}}(\operatorname{red}\mathfrak C);
\]
\item if \(\mathfrak C_1\sim\mathfrak C_2\) are cylindrical, then \(\Pi_{\mathrm{cyl}}(\mathfrak C_1)
 \sim\Pi_{\mathrm{cyl}}(\mathfrak C_2).\)
\end{enumerate}
In particular, projection defines a map on equivalence classes of
cylindrical scattering diagrams.
\end{proposition}
\begin{proof}
By Lemma~\ref{lem:order-preservation}, projection preserves order and commutes with truncation.  It does not increase the number of nontrivial walls at any finite order, proving (i).

For any constant-\(U\) lift \(\widetilde\gamma\) of an admissible base path \(\gamma,\) the crossing times and their order agree. Moreover, for a lifted wall degree \(v=w+u,\)
\[
 \omega^\dagger(v,-\widetilde\gamma')
 =\omega(w,-\gamma'),
\]
so the crossing signs agree as well. Applying \(\pi_{\mathrm{bal}}\) to each truncated path product and passing to the inverse limit proves (ii), which immediately implies (iv). Finally, choose subdivisions pulled back from the base. Reduction preserves path products and these subdivisions preserve cylindricity, so (ii) also gives (iii).
 \end{proof}

Fix a hyperplane \(H\) and a point \(x\in H\) that is general for all diagrams \(\mathfrak C_a\) under consideration. Thus, at each finite order, a sufficiently small neighbourhood of \(x\) meets only walls whose supports agree locally with \(H.\) Choose an oriented rational ray \(\rho\subset H^\perp,\) set \(n_v:=\omega^\dagger(v,-),\) and choose a transverse vector \(\xi\) such that \(n(\xi)<0\) for every nonzero \(n\in\rho.\)

For each \(N,\) choose \(\varepsilon_N>0\) sufficiently small that the segment
\[
 \gamma_{x,\xi}(s)=x+s\xi,
 \qquad -\varepsilon_N\leq s\leq\varepsilon_N,
\]
is admissible to order \(N\) and crosses only the walls through \(x.\) Define
\[
 G^\xi_{a,N}
 :=\PathProd^{\leq N}_{\gamma_{x,\xi},\mathfrak C_a}.
\]
The local wall group is the direct product of the two commuting
subgroups whose logarithms have homogeneous degrees \(v\) with
\(n_v\in\rho\) and \(n_v\in-\rho,\) respectively. Let \(\pr_{\rho,N}\) and \(\pr_{-\rho,N}\) denote the corresponding projections, and define
\begin{equation}\label{eq:normalized-local-ray-products}
 F_{a,H,\rho,N}:=\pr_{\rho,N}(G^\xi_{a,N}),
 \qquad
 F_{a,H,-\rho,N}
 :=\bigl(\pr_{-\rho,N}(G^\xi_{a,N})\bigr)^{-1}.
\end{equation}
With our crossing-sign convention, \(F_{a,H,\pm\rho,N}\) is precisely the product, modulo order \(>N,\) of the wall elements through \(x\) whose effective covectors lie in \(\pm\rho.\) The point \(x\) is suppressed from the notation. Write \(r_{N',N}\) for the natural truncation map from order \(N'\) to order \(N.\)

\begin{lemma}[Normalized local wall products]
\label{lem:compatible-generic-germs}
Let \(\mathfrak C_1\sim\mathfrak C_2\) be orderwise finite diagrams with effective wall degrees. For the local data \((H,x,\rho)\) chosen above,
\begin{equation}\label{eq:local-pro-germ-inverse-system}
 F_{1,H,\pm\rho,N}=F_{2,H,\pm\rho,N},
 \qquad
 r_{N',N}(F_{a,H,\pm\rho,N'})
 =F_{a,H,\pm\rho,N}
\end{equation}
for \(a=1,2\) and \(N'\geq N.\) These products are independent of the sufficiently short transverse segment with the prescribed orientation, and of finite rational subdivisions for which \(x\) remains general.
\end{lemma}

\begin{proof}
All effective covectors of walls through \(x\) have kernel \(H.\)
If \(n_v=c\,n_{v'},\) then
\[
 \omega^\dagger(v,v')
 =n_v(v')=c\,n_{v'}(v')=0.
\]
Hence the local homogeneous Lie components commute.  Separating their degrees according to the effective covector rays
\(\rho\) and \(-\rho\) gives the direct-product decomposition of the local wall group.

Taking products over the walls of \(\mathfrak C_a\) through \(x\),
the Davison-Mandel crossing convention gives
\[
 G^\xi_{a,N}
 =
 \prod_{n_{v_{\mathfrak d}}\in\rho}
       f_{\mathfrak d}^{\leq N}
 \left(
 \prod_{n_{v_{\mathfrak d}}\in-\rho}
       f_{\mathfrak d}^{\leq N}
 \right)^{-1}.
\]
Thus normalization recovers the product of the wall elements on
each effective covector ray, without crossing signs.
Equivalence identifies the signed path products, so applying the
two projections proves the first equality.
Truncation commutes with these products and projections, proving
the second.

Changing the sufficiently short transverse segment, or subdividing
the walls while keeping \(x\) general, leaves these local products
unchanged.
\end{proof}

Define the completed local Lie algebra associated with $H$ by
\begin{equation}\label{eq:effective-ray-decomposition}
 \widehat{\mathfrak g}_H
 :=\prod_{\substack{v\in(L_0^\dagger)^+\\
                   n_v\ne0,\ \ker n_v=H}}
       \mathfrak g_v
 =\prod_{\rho\subset H^\perp}
       \widehat{\mathfrak g}_{H,\rho},
 \qquad
 \widehat{\mathfrak g}_{H,\rho}
 :=\prod_{\substack{v\in(L_0^\dagger)^+\\
                   \RR_{>0}n_v=\rho}}
       \mathfrak g_v.
\end{equation}
The outer product ranges over the two oriented rational rays
in \(H^\perp\); either factor may be zero.

\begin{lemma}[Raywise BPS uniqueness]
\label{lem:raywise-bps-uniqueness}
The Lie algebra \(\widehat{\mathfrak g}_H\) is abelian.
If \(g\in\exp(\widehat{\mathfrak g}_H)\) admits an orderwise finite
factorization
\begin{equation}\label{eq:local-bps-factorization}
 g=
 \prod_{\substack{v\in(L_0^\dagger)^+\\
                 n_v\ne0,\ \ker n_v=H}}
       \E(-p_v(t)x^v),
 \qquad p_v(t)\in\ZZ[t^{\pm1}],
\end{equation}
then the polynomials \(p_v(t)\) are uniquely determined by \(g\),
indexed by the actual lifted degree \(v\).
Their recovery commutes with truncation and applies separately
to each normalized local product \(F_{H,\rho}\).
\end{lemma}
\begin{proof}
For any two degrees occurring in
\(\widehat{\mathfrak g}_H\), their effective covectors are
proportional.  Writing \(n_v=c\,n_{v'}\), we obtain
\[
 \omega^\dagger(v,v')
 =n_v(v')=c\,n_{v'}(v')=0.
\]
The quantum-torus bracket therefore vanishes between the
corresponding homogeneous components.  This proves abelianness
and shows that \eqref{eq:effective-ray-decomposition} is a
decomposition of completed Lie algebras.

Since the factors commute, \(\log g\) is the sum of their logarithms,
and
\begin{equation}\label{eq:bps-log-divisibility-recursion}
 \log\E(-p_v(t)x^v)
 =
 \sum_{\ell\geq1}
 \frac{(-1)^{\ell-1}\psi_\ell(p_v(t))}
      {\ell(t^\ell-t^{-\ell})}x^{\ell v}.
\end{equation}
The coefficient of \(x^v\) in \(\log g\) is therefore
\(p_v(t)/(t-t^{-1})\), together with contributions from proper
divisors \(u\) satisfying \(v=\ell u\), \(\ell>1\).
Since
\[
 \ord^\dagger(u)=\frac{\ord^\dagger(v)}{\ell}
 <\ord^\dagger(v),
\]
these contributions are already determined by induction.
Subtracting them uniquely recovers \(p_v(t)\).

The monomials \(x^v\) are linearly independent in the lifted
quantum torus, so this argument distinguishes actual degrees
even when they have the same effective covector or differ by
a radical vector.  For \(\ord^\dagger(v)\leq N\), the recovery
uses only \(g\) modulo order \(>N\); hence it commutes with
truncation.  Restricting to degrees with \(n_v\in\rho\) proves
the assertion for \(F_{H,\rho}\).
\end{proof}

\begin{proposition}[Positive cylinders by finite-order construction]
\label{prop:global-positive-cylinder-gluing}
Let \(\mathfrak C\sim\mathfrak D_+\) be conical, orderwise finite
diagrams in the central lifted setting, with effective wall degrees.
Assume that \(\mathfrak C\) is cylindrical and \(\mathfrak D_+\) is positive.  Then their equivalence class contains an orderwise finite positive cylindrical diagram \(\mathfrak C_+.\)
Its normalized local BPS coefficients, indexed by actual lifted
degrees, depend only on the equivalence class.

If both diagrams contain a positive cylindrical incoming diagram
\(\mathfrak I,\) and every wall of
\(\mathfrak D_+\setminus\mathfrak I\) is outgoing, then
\(\mathfrak C_+\) may be chosen to contain \(\mathfrak I,\) with
every other wall outgoing.
\end{proposition}

\begin{proof}
We construct full positive \(E\)-factors in increasing BPS degree,
retaining every factor once inserted.
For the first assertion, put \(\mathfrak I=\varnothing;\)
for the relative assertion, retain the prescribed incoming diagram.
Throughout, resolve wall packets into ordinary walls on their
actual degree rays.

Set \(\mathfrak C_{<1}:=\mathfrak I.\)
Suppose that \(\mathfrak C_{<N}\) is cylindrical, consists of
\(\mathfrak I\) and finitely many full positive \(E\)-factors
of BPS degree less than \(N,\) and satisfies
\[
 \mathfrak C_{<N}\sim\mathfrak C
 \qquad\text{modulo order \(>N-1\)}.
\]

At order \(N,\) the diagrams \(\mathfrak C\) and
\(\mathfrak C_{<N}\) have finitely many cylindrical supports.
Moreover, since \(U\subset\ker\omega^\dagger,\) the supporting
hyperplane of every wall of \(\mathfrak D_+^{\leq N}\) is also
cylindrical. Let \(\mathcal H_N\) be the finite family of base
hyperplanes obtained from the supporting hyperplanes of these
three truncated diagrams.

For each \(\bar H\in\mathcal H_N,\) choose a finite rational
conical subdivision \(\Sigma_{\bar H,N}\) adapted to the
cylindrical base supports and to all intersections
\(\bar H\cap\bar H',\) with \(\bar H'\ne\bar H.\)
Such a subdivision is obtained from the finitely many homogeneous
linear inequalities defining those supports and intersections.
On each relatively open top-dimensional cell
\(\sigma\in\Sigma_{\bar H,N},\) the normalized local products
of both cylindrical diagrams are constant.
These subdivisions are auxiliary: previously inserted walls
are neither subdivided nor modified.

Put \(H:=\bar H\oplus U_\RR,\) and choose an oriented ray
\(\rho\subset H^\perp.\)
Write \(F_{\sigma,\rho}^{\leq N}\) and
\(B_{\sigma,\rho}^{\leq N}\) for the normalized local products
of \(\mathfrak C\) and \(\mathfrak C_{<N},\) respectively,
on \(\sigma\oplus U_\RR;\) an empty product is the identity.
The boundaries and joints of \(\mathfrak D_+^{\leq N}\) form a finite union of subsets of codimension at least two. Hence general points of \(\sigma\oplus U_\RR\) avoid them.

At each such point, equivalence and
Lemma~\ref{lem:compatible-generic-germs} identify
\(F_{\sigma,\rho}^{\leq N}\) with the normalized local product
of \(\mathfrak D_+.\)
Its BPS coefficients
\[
 p_{\sigma,v}(t)\in\ZZ_{\geq0}[t^{\pm1}],
 \qquad
 \ord^\dagger(v)\leq N,\quad n_v\in\rho,
\]
are therefore nonnegative.
By Lemma~\ref{lem:raywise-bps-uniqueness}, they are independent
of the chosen general point in \(\sigma\oplus U_\RR\) and
remain unchanged under further subdivision of \(\sigma.\)
No boundary of a noncylindrical wall is projected in this argument.

Consider the residual local product
\[
 R_{\sigma,\rho,N}
 :=F_{\sigma,\rho}^{\leq N}
   \bigl(B_{\sigma,\rho}^{\leq N}\bigr)^{-1}.
\]
The induction hypothesis makes it trivial modulo order \(>N-1.\)
Since the local ray groups are abelian, the BPS coefficients
of \(F_{\sigma,\rho}^{\leq N}\) and
\(B_{\sigma,\rho}^{\leq N}\)  agree in all degrees below \(N.\)
Their multiple-cover contributions consequently cancel, giving
\begin{equation}\label{eq:finite-order-cylinder-residual}
 \log R_{\sigma,\rho,N}
 =
 \sum_{\substack{\ord^\dagger(v)=N\\ n_v\in\rho}}
 \frac{a_{\sigma,v}(t)}{t-t^{-1}}x^v
 \qquad\pmod{\text{order }>N}.
\end{equation}

Let \(i_{\sigma,v}(t)\) be the local degree-\(v\) BPS coefficient
of \(\mathfrak I,\) with \(i_{\sigma,v}=0\) when
\(\mathfrak I=\varnothing.\)
All previously inserted factors outside \(\mathfrak I\)
have BPS degree less than \(N.\) Thus
\begin{equation}\label{eq:finite-order-cylinder-positive-coefficient}
 a_{\sigma,v}(t)
 =p_{\sigma,v}(t)-i_{\sigma,v}(t)
 \in\ZZ_{\geq0}[t^{\pm1}].
\end{equation}
Indeed, at a general point of \(\sigma\oplus U_\RR,\) this difference is precisely the degree-\(v\) BPS coefficient of the positive local product of \(\mathfrak D_+\setminus\mathfrak I.\)
The prescribed factors cancel exactly, while the lower-degree
multiple-cover contributions have already been removed.

Let \(\mathfrak W_N\) be the indexed family of full walls
\begin{equation}\label{eq:finite-order-cylinder-new-wall}
 \left(
   \overline{\sigma}\oplus U_\RR,\,
   \E(-a_{\sigma,v}(t)x^v)
 \right),
 \qquad
 \ord^\dagger(v)=N,\quad a_{\sigma,v}\ne0,
\end{equation}
as \(\bar H,\) \(\sigma,\) and \(\rho\) vary.
Each support is a closed rational polyhedral cone contained
in \(\ker n_v,\) so each pair defines a wall.
If \(v\) is nonprimitive, the wall element belongs to the group
of its primitive degree ray.
There are finitely many cells and finitely many lifted degrees
of order \(N;\) hence \(\mathfrak W_N\) is finite.
Set
\begin{equation}\label{eq:finite-order-cylinder-recursion}
 \mathfrak C_{<N+1}
 :=\mathfrak C_{<N}\sqcup\mathfrak W_N.
\end{equation}
The degree-\(N\) logarithms of the inserted factors are exactly
the summands in \eqref{eq:finite-order-cylinder-residual}.
Thus the normalized local products agree with those of
\(\mathfrak C\) to order \(N,\) separately on both oriented rays.
Consequently,
\[
 \mathfrak C_{<N+1}\sim\mathfrak C
 \qquad\text{modulo order \(>N\)}.
\]

For the relative assertion, suppose \(a_{\sigma,v}\ne0.\)
Every general point of \(\sigma\oplus U_\RR\) then lies on
a wall of \(\mathfrak D_+^{\leq N}\setminus\mathfrak I\) carrying a nonzero degree-\(v\) BPS factor.
The finitely many closed supports of these walls therefore
cover \(\overline{\sigma}\oplus U_\RR.\)
Each wall has actual degree ray \(\RR_{>0}v\) and is conical
and outgoing, so its support does not contain \(v.\)
Hence \(\overline{\sigma}\oplus U_\RR\) does not contain \(v\) either, proving that every inserted wall is outgoing.

Finally, define
\[
 \mathfrak C_+
 :=\mathfrak I\sqcup\bigsqcup_{N\geq1}\mathfrak W_N.
\]
Only stages \(1,\ldots,N\) contribute modulo order \(>N.\)
Thus \(\mathfrak C_+\) is an orderwise finite positive cylindrical
diagram, and the finite-order equivalences give
\(\mathfrak C_+\sim\mathfrak C.\)
Reduction using subdivisions pulled back from the base preserves
positivity and cylindricity.
Lemma~\ref{lem:raywise-bps-uniqueness} gives uniqueness and
truncation compatibility of the normalized local BPS coefficients;
the polyhedral representative itself need not be unique.

For \(v=w(\alpha)+u(\alpha),\) the balanced projection of
\eqref{eq:finite-order-cylinder-new-wall} is
\begin{equation}\label{eq:finite-order-base-new-wall}
 \left(
   \overline{\sigma},\,
   \E\bigl(-a_{\sigma,v}(t)\mathbf y^\alpha
                x^{w(\alpha)}\bigr)
 \right).
\end{equation}
Here \(\overline{\sigma}\subset w(\alpha)^{\omega\perp}\)
is rational polyhedral, and the BPS polynomial is nonnegative
in the free central variables.
By Lemma~\ref{lem:order-preservation}, the projected factor
has order \(N.\)
\end{proof}

\begin{lemma}[Positive cylinder representative]
\label{lem:cylinder-representative}
The class of \(\Ddag\) contains an orderwise finite positive cylinder
representative retaining \(\Ddag_{\mathrm{in}},\) with every other
wall outgoing.  Its wall elements are commuting products of
\(\E(-t^a x^v),\) where \(a\in\ZZ\) and
\(v\in(L_0^\dagger)^+.\)  Its orderwise reductions remain positive
and cylindrical when subdivisions are pulled back from the base.
\end{lemma}

\begin{proof}
Lemma~\ref{lem:cylinder-pullback} gives the cylindrical diagram
\(\Pr^*\Db.\) Independently, Theorem~\ref{thm:DM} gives an
equivalent positive diagram with the same prescribed incoming walls
and only outgoing additional walls.  
Using Proposition~\ref{prop:DM-hypotheses-central-lift} and Lemma~\ref{lem:radical-degrees-not-generated}, we apply Proposition~\ref{prop:global-positive-cylinder-gluing} to obtain the result.

\end{proof}

\begin{proposition}[Independence of the positive cylinder choice]
\label{prop:positive-cylinder-choice}
Let \(\mathfrak D_{\mathrm{cyl},+}^{(1)}\) and \(\mathfrak D_{\mathrm{cyl},+}^{(2)}\) be two cylindrical,
positive representatives of the equivalence class of \(\Ddag.\) Then
\[
 \Pi_{\mathrm{cyl}}
   (\mathfrak D_{\mathrm{cyl},+}^{(1)})
 \sim
 \Pi_{\mathrm{cyl}}
   (\mathfrak D_{\mathrm{cyl},+}^{(2)}).
\]
More strongly, at every finite order and general oriented germ, the
local BPS polynomial attached to each triple \((H,\rho,v),\) where
\(\RR_{>0}\omega^\dagger(v,-)=\rho,\) is the same for the two
representatives. Thus the projected equivalence class and its local lifted BPS data
are independent of the choice; the individual wall supports need not be.
\end{proposition}

\begin{proof}
The first assertion is
Proposition~\ref{prop:cylinder-projection-formal}(iv). For the stronger
assertion, fix \(H\) and \(\rho.\) Equivalence first identifies the
signed products \(G^\xi_{1,N}=G^\xi_{2,N}.\) Lemma
\ref{lem:compatible-generic-germs} therefore identifies the normalized
systems on \(\rho\) and \(-\rho.\) Applying
Lemma~\ref{lem:raywise-bps-uniqueness} separately to these systems
identifies every actual lifted degree, and the divisor recursion
commutes with truncation. Hence the inverse-limit BPS data agree. Applying \(\pi_{\mathrm{bal}}\) records, rather than forgets, each lifted degree as a central coefficient monomial. This proves the last statement with the precise degree convention explained in Remark~\ref{rem:projected-bps-degree} below.
\end{proof}

\subsection{Projection}

For the base normalized Lie algebra $\g,$ put
$\g^*:=\g\otimes_{\QQ[t^{\pm1}]}\QQ(t)[y_{ij}].$ Its balanced completed subalgebra is the completion of
$\overline{\g}.$ Fix any positive
cylinder representative
\(\mathfrak D_{\mathrm{cyl},+}\) supplied by
Lemma~\ref{lem:cylinder-representative}.

\begin{definition}[Projected diagram]
\label{def:Dstar}
Define
\[
 \Dstar
 :=\operatorname{red}
   \Pi_{\mathrm{cyl}}(\mathfrak D_{\mathrm{cyl},+}).
\]
Equivalently, a positive factor
\[
 \left(\overline{\mathfrak d}\oplus U_\RR,
       \E(-t^a x^v)\right),
 \qquad v=w(\alpha)+u(\alpha),
\]
projects to
\[
 \left(\overline{\mathfrak d},
       \E(-t^a\mathbf y^\alpha x^{w(\alpha)})\right).
\]
\end{definition}

Proposition~\ref{prop:positive-cylinder-choice} shows that the
equivalence class of \(\Dstar\), together with its normalized local
BPS data with central multi-degree retained, is independent of
\(\mathfrak D_{\mathrm{cyl},+}\).  Thus Definition~\ref{def:Dstar}
specifies a canonical equivalence class, not an unrecorded choice of
representative.

\begin{remark}[Lifted versus projected BPS degree]
\label{rem:projected-bps-degree}
There are three distinct, and compatible, meanings of the degree of
the BPS data.  On a lifted local pro-germ the canonical polynomials are indexed by the actual degrees
\[
 v=w(\alpha)+u(\alpha)\in(L_0^\dagger)^+.
\]
The balanced projection is injective on these degrees:
\[
 x^{w(\alpha)+u(\alpha)}
 \longmapsto \mathbf y^\alpha x^{w(\alpha)}.
\]
Thus it does not merge two lifted BPS invariants.  If one chooses to
group all coefficient monomials having the same geometric projected degree \(w\), the intrinsic coefficient-refined datum is
\begin{equation}\label{eq:projected-bps-polynomial}
 P_{H,w}(t,\mathbf y)
 :=\sum_{\substack{\alpha\in\ZZ_{\geq0}^{\sum_i r_i}\\
                    w(\alpha)=w}}
 p_{H,w+u(\alpha)}(t)\mathbf y^\alpha .
\end{equation}
The sum is finite at every order, and plethystic additivity gives
\[
 \prod_{\alpha:\,w(\alpha)=w}
 \E(-p_{H,w+u(\alpha)}(t)\mathbf y^\alpha x^w)
 =\E(-P_{H,w}(t,\mathbf y)x^w).
\]
Consequently independence is degree-by-degree in \(v\), equivalently
coefficient-by-coefficient in \(P_{H,w}(t,\mathbf y)\).

After a noninjective coefficient specialization
\(y_{ij}\mapsto q_{ij}\), or after forgetting the \(U\)-multi-degree,
different \(\alpha\)'s can merge.  The universal construction gives the evaluated aggregate
\[
 P_{H,w}(t,\mathbf q)
 =\sum_{w(\alpha)=w}
   p_{H,w+u(\alpha)}(t)\mathbf q^\alpha
\]
for each projected degree \(w\), independently of the universal
positive representative.  Arbitrary coefficient evaluation need not
commute with the plethystic substitutions, so no recovery of these
aggregates from the specialized diagram alone is asserted.  Individual
lifted summands are retained by the universal central lift.
\end{remark}

\begin{remark}[Root positivity versus exchange-coefficient positivity]
\label{rem:root-vs-elementary-positivity}
The positivity obtained from the central lift is positivity in the
independent lifted/root variables, equivalently coefficientwise
positivity in the central multi-degree before specialization.  It does
not in general imply positivity after rewriting a symmetric polynomial
in the elementary reciprocal coefficients $\mathbf z_{i,s}$.  For
example, if $q_1q_2=1$ and $q_1+q_2=\mathbf z$, then
$q_1^2+q_2^2=\mathbf z^2-2$.  Accordingly, all positivity statements
below are formulated as root-positive or lifted-multidegree-positive
statements unless an additional coefficientwise theorem is explicitly
invoked.
\end{remark}

\begin{lemma}[Projected walls are legitimate]
\label{lem:legitimate-projection}
The diagram $\Dstar$ is a well-defined scattering diagram on $L_\RR$ for the Lie algebra $\g^*$.
\end{lemma}

\begin{proof}
By Lemma~\ref{lem:cylinder-representative}, each support is
\(\overline{\mathfrak d}\oplus U_\RR\); hence
\(\overline{\mathfrak d}\) is a rational polyhedral cone of codimension
one in \(L_\RR\).  If \(v=w(\alpha)+u(\alpha)\) is a homogeneous lifted
wall degree, then for every \(\ell+u\in L^\dagger\),
\[
 \omega^\dagger(v,\ell+u)=\omega(w(\alpha),\ell).
\]
Since
\(\overline{\mathfrak d}\oplus U_\RR
 \subseteq v^{\omega^\dagger\perp}\),
taking \(u=0\) gives the exact support condition
\[
 \overline{\mathfrak d}
 \subseteq w(\alpha)^{\omega\perp}.
\]
Moreover,
\(\mathbf y^\alpha z^{w(\alpha)}\) belongs to the balanced subalgebra
of \(\g^*\).  Subdivision, products on a fixed actual degree ray,
and deletion of identity walls preserve legitimacy.
Orderwise local finiteness, including the case in which several
cylinders project to one support or one BPS product is displayed as
several factors, is
Proposition~\ref{prop:cylinder-projection-formal}(i).
\end{proof}

\begin{lemma}[Consistency of the projected diagram]
\label{lem:Dstar-consistent}
The projected diagram $\Dstar$ is consistent and
\[
  \Dstar\sim\Db.
\]
\end{lemma}

\begin{proof}
The positive cylinder representative is consistent because it is
equivalent to \(\Ddag\).  For an admissible closed base loop
\(\gamma\), choose any constant-\(U\) lift; admissibility of the base
loop implies admissibility of the lift because all walls, boundaries,
and joints are full inverse images.  Equation
\eqref{eq:formal-projected-path-product} gives
\[
 \PathProd_{\gamma,\Pi_{\mathrm{cyl}}
   (\mathfrak D_{\mathrm{cyl},+})}
 =
 \pi_{\mathrm{bal}}\!
 \left(\PathProd_{\widetilde\gamma,
   \mathfrak D_{\mathrm{cyl},+}}\right)=1.
\]
Reduction does not change path products, so \(\Dstar\) is consistent.

The cylinder pullback satisfies the exact identity
\[
 \Pi_{\mathrm{cyl}}(\Pr^*\Db)=\Db,
\]
because \(\pi_{\mathrm{bal}}\circ\iota=1\) on every wall element and
projection is inverse to taking the full cylinder on every support.
Moreover
\(\mathfrak D_{\mathrm{cyl},+}\sim\Pr^*\Db\).
Proposition~\ref{prop:cylinder-projection-formal}(iii)--(iv) therefore
gives
\[
 \Dstar
 \sim\Pi_{\mathrm{cyl}}(\mathfrak D_{\mathrm{cyl},+})
 \sim\Pi_{\mathrm{cyl}}(\Pr^*\Db)
 =\Db.
\]
\end{proof}

Combining the previous two lemmas gives the following theorem.

\begin{theorem}[Positive cylinder projection theorem]
\label{thm:cylinder-projection}
Under \eqref{eq:active-scattering-directions}, \(\Ddag\) admits a
positive, orderwise finite cylinder representative retaining
the prescribed incoming walls, with all other walls outgoing.
Its projection \(\Dstar\) is consistent and equivalent to \(\Db\).
Projection preserves truncations and path products.  The projected
equivalence class and normalized local BPS data, indexed by actual
lifted degree, are independent of the chosen representative.
\end{theorem}
\begin{proof}
The existence of the positive cylinder representative follows from
Proposition~\ref{prop:global-positive-cylinder-gluing} and
Lemma~\ref{lem:cylinder-representative}.  The projection properties
are given by Lemma~\ref{lem:legitimate-projection} and
Proposition~\ref{prop:cylinder-projection-formal}; consistency and
the equivalence \(\Dstar\sim\Db\) follow from
Lemma~\ref{lem:Dstar-consistent}.  Independence of the normalized
local BPS data is Proposition~\ref{prop:positive-cylinder-choice},
with the specialization qualification in
Remark~\ref{rem:projected-bps-degree}.
\end{proof}

\begin{corollary}[Plethystic coefficient positivity]
\label{cor:plethystic-coefficient-positivity}
Let
\[
A=\ZZ[t^{\pm1},y_1,\ldots,y_s],\qquad
A_{\geq0}=\ZZ_{\geq0}[t^{\pm1},y_1,\ldots,y_s],
\]
equipped with the substitutions
\(\psi_\ell(t)=(-1)^{\ell-1}t^\ell\) and
\(\psi_\ell(y_a)=y_a^\ell\) from
\eqref{eq:explicit-plethystic-substitution}.  Suppose a finite
incoming diagram has
effective, linearly independent wall directions
\(v_1,\ldots,v_r\), freely
generating the positive monoid used for the completion, and consists of
the full incoming walls
\[
 \left(v_i^{\omega\perp},\E_A(-p_i(t,\mathbf y)x^{v_i})\right),
 \qquad p_i(t,\mathbf y)\in A_{\geq0}.
\]
Then its consistent completion admits a reduced representative whose
wall elements are commuting products of factors
\[
 \E_A(-t^a\mathbf y^\alpha x^v),
 \qquad
 a\in\ZZ,\quad
 \alpha\in\ZZ_{\geq0}^s,\quad v\in L_0^+.
\]
Thus positivity holds coefficient by coefficient in the plethystic
variables.
\end{corollary}

\begin{proof}
Write
\[
 p_i(t,\mathbf y)
 =\sum_{\alpha\in F_i}c_{i,\alpha}(t)\mathbf y^\alpha,
 \qquad
 c_{i,\alpha}(t)\in\ZZ_{\geq0}[t^{\pm1}],
\]
where index set \(F_i\) is finite.  Introduce a free central lattice with one
basis vector \(u_{i,\alpha}\) for every pair
\((i,\alpha)\), and use the lifted incoming walls
\[
 \left(v_i^{\omega\perp}\oplus U_\RR,\,
       \E(-c_{i,\alpha}(t)x^{v_i+u_{i,\alpha}})\right),
 \qquad \alpha\in F_i.
\]
Theorem~\ref{thm:DM} applies to this ordinary lifted diagram.
The proofs of the cylinder-pullback and local-product argument of
Lemmas~\ref{lem:cylinder-pullback} and
\ref{lem:cylinder-representative} apply with the pairs
\((i,\alpha)\) as the free generator set: only the independent central
coordinates, centrality of \(U\), and the free positive filtration are
used.  They give it a positive cylindrical representative.

Project by the coefficient map
\[
 x^{u_{i,\alpha}}\longmapsto\mathbf y^\alpha.
\]
For every \(\ell\geq1\), this map sends
\[
 x^{\ell u_{i,\alpha}}
 \longmapsto\mathbf y^{\ell\alpha}
 =\psi_\ell(\mathbf y^\alpha).
\]
Hence, directly from
\eqref{eq:plethystic-E-with-coefficients}, the projected incoming
product is
\[
 \prod_{\alpha\in F_i}
 \E_A(-c_{i,\alpha}(t)\mathbf y^\alpha x^{v_i})
 =\E_A(-p_i(t,\mathbf y)x^{v_i}),
\]
and every projected positive factor has the form stated above.
Cylinder projection preserves consistency by
Lemma~\ref{lem:Dstar-consistent}.  This proves the claim without any
one-parameter specialization argument.
\end{proof}

\section{Positive factorization of wall functions}\label{sec:positivity}

The projection theorem allows us to descend the positivity theorem of Davison--Mandel from the lifted ordinary quantum scattering diagram to the generalized quantum scattering diagram.

\begin{theorem}[Positive factorization over the splitting field]
\label{thm:positive-factorization}
Under the active-index convention
\eqref{eq:active-scattering-directions}, let $\mathfrak D$ be the
generalized quantum scattering diagram of
Definition~\ref{def:GQSD}.  After base change to the splitting field
$K_{\mathbf z}$, the diagram $\mathfrak D$ admits an equivalent
representative whose wall functions are products of factors
\[
  \E\left(-t^a\prod_{i=1}^n\prod_{j=1}^{r_i}q_{ij}^{a_{ij}}x^v\right),
  \qquad
  a\in\ZZ,
  \quad a_{ij}\in\ZZ_{\ge0},
  \quad v\in L_0^+.
\]
Each displayed factor denotes the orderwise evaluation
\(y_{ij}\mapsto q_{ij}\) of its free-root counterpart; no
plethystic substitution is imposed on \(K_{\mathbf z}\).
\end{theorem}

\begin{proof}
Use the positive cylinder representative of
Lemma~\ref{lem:cylinder-representative}.  Its factors have the form
$\E(-t^a x^v)$ with $v\in (L_0^\dagger)^+$.  Write
\[
  v=\sum_{i,j}a_{ij}f_j^i
  =\sum_{i,j}a_{ij}(e_i+u_{ij}).
\]
Then
\[
  \Pr(v)=\sum_i\left(\sum_j a_{ij}\right)e_i
\]
and
\[
  v-\Pr(v)=\sum_{i,j}a_{ij}u_{ij}.
\]
Under projection, the factor becomes
\[
  \E\left(-t^a \prod_{i,j} y_{ij}^{a_{ij}}x^{\Pr(v)}\right).
\]
Apply the orderwise coefficient map \(y_{ij}\mapsto q_{ij}\) into
\(K_{\mathbf z}\).  For the incoming walls this gives
\[
  \prod_{j=1}^{r_i}\E(-q_{ij}x^{e_i})=\Phi_i(x^{e_i}),
\]
so the specialized projected diagram is equivalent to the generalized
quantum diagram \(\mathfrak D\) after scalar extension to
\(K_{\mathbf z}\).  The desired factorization follows.
Here the map means orderwise evaluation of the universal
coefficient series: \(y_{ij}^{\ell}\) is sent to
\(q_{ij}^{\ell}\) in its \(\ell\)-th substitution term.  Thus the
evaluated series is exactly the factor displayed above.  This is
ordinary coefficientwise evaluation of an already formed series; no
additional operation on the algebraic splitting field is defined.
\end{proof}

\begin{remark}[Specialization is taken on automorphisms]
\label{rem:no-specialization-wall-element}
The formal wall element itself contains denominators
$t^m-t^{-m}$ and therefore need not admit a value at $t=1$.
We take classical limits of wall-crossing maps through their adjoint
actions. The finite-product formula is regular at \(t=1;\) it therefore defines an orderwise specialization of wall-crossing automorphisms whenever the exchange coefficients are regular at
\(t=1.\) In particular, no argument below uses a formal expression
such as \(\Phi_i|_{t=1}\).
\end{remark}

\begin{remark}[Normalized classical limit]
\label{rem:normalized-classical-limit}
On the normalized Lie algebra of
\eqref{eq:normalized-quantum-torus-Lie}, define
\[
 \operatorname{cl}(H):=\lim_{t\to1}(t-t^{-1})H,
 \qquad
 \text{for } H\in \g
\]
Then $\operatorname{cl}(\hat{x}^v)=\frac{x^v}{|v|_0}.$
Every generator and iterated bracket has at most a simple pole at
$t=1$, since the commutator supplies a factor $t-t^{-1}$.
This defines a filtered Lie homomorphism to the classical algebra
with bracket $\{x^a,x^b\}=\omega(a,b)x^{a+b}$ and extends orderwise
to the completions, as in \cite[Section~2.2.4]{Davison_Mandel21}. Equation~\eqref{eq:normalized-Lie-adjoint} shows that it intertwines the adjoint actions.  It is a specialization of the normalized form, not an evaluation on its full $\QQ(t)$-scalar extension.
\end{remark}

\section{Seed scattering diagrams}\label{sec:seed-diagrams}

We now attach the preceding construction to a compatible generalized
quantum seed with formal reciprocal exchange coefficients. All seed
diagrams and mutation statements in this section and the next are formed over \(R_{\mathbf z},\) with scalar extensions used to prove positivity.  Coefficient evaluation is deferred to Section~\ref{sec:coefficient-evaluation}. We distinguish the compatibility diagonal from its inverse when comparing the scattering and Fock--Goncharov normalizations.

\subsection{Compatible seed data and normalized incoming walls}

Fix the mutation data \((\mathbf r,\mathbf z)\) of Section~\ref{sec:generalized-diagram}.

Fix an integer \(m\ge n.\) Fix a lattice \(N\) of rank \(m\) and let
\(M=N^*=\operatorname{Hom}_{\mathbb Z}(N,\mathbb Z).\) Write \(I=\{1,\ldots,m\},\) \(F=\{n+1,\ldots,m\},\) and \(I_{\uf}=I\setminus F=\{1,\ldots,n\}.\)

A \emph{compatible marked seed} is a tuple
\[
 S=(E,\widetilde B,\Lambda,\mathbf r,\mathbf z),
\]
where:
\begin{enumerate}[label=(\roman*), noitemsep]
\item \(E_S=\{e_{i;S}\}_{i\in I}\) is a \(\mathbb Z\)-basis of \(N;\)
\item \(\widetilde B_S=(b_{ji;S})_{1\le j\le m,\,1\le i\le n}
      \in\operatorname{Mat}_{m\times n}(\mathbb Z)\) is the extended exchange matrix in the seed coordinates;
\item \(\Lambda_S=(\lambda_{ij;S})_{1\le i,j\le m}\) is the matrix, in the dual basis \(E_S^*=\{e_{i;S}^*\}_{i\in I},\) of an integral
      skew-symmetric form $\Lambda$ on \(M,\) satisfying
      \[
       \widetilde B_S^T\Lambda_S=[I_n\mid0].
      \]
\end{enumerate}
When we need to emphasize the basis, we write \(e_{i;S}:=e_i\) and
\(e_{i;S}^*:=e_i^*;\) these are elements of the \emph{same} \(N\) and
\(M\) for every seed.

Define
\[
  p^*:N_{\uf}\longrightarrow M,
  \qquad
  p^*(e_i):=\sum_{j=1}^m b_{ji;S}e_j^*,
\] and write
$v_i=p^*(e_i).$ Compatibility gives, for every $m\in M$,
\begin{equation}\label{eq:seed-compatibility-normalized}
  \Lambda(v_i,m)=\langle e_i,m\rangle.
\end{equation}
More generally, for $\alpha\in\ZZ^n$ and $\beta\in\ZZ^m$,
\begin{equation}\label{eq:compatibility-dictionary-basic}
  \Lambda(\widetilde B\alpha,\beta)
  =\alpha^TI_n \overline\beta
  =(\alpha,\overline\beta),
\end{equation}
where $\overline\beta$ denotes the first $n$ coordinates of $\beta$. By the compatibility of $(\widetilde B,\Lambda),$ we have the following lemma.

\begin{lemma}[Compatibility implies injectivity]
\label{lem:compatible-pstar-injective}
The map
\[
  p^*\vert_{N_{\uf}}:N_{\uf}\longrightarrow M
\]
is injective. Equivalently, the columns of $\widetilde B$ are linearly
independent over $\ZZ$ (and over $\QQ$).
\end{lemma}

\begin{lemma}[Compatible seed indices are automatically active]
\label{lem:compatible-seed-indices-active}
Let \(S=(E,\widetilde B,\Lambda,\mathbf{r},\mathbf{z})\) be a compatible seed. For every unfrozen index \(i\in I_{\uf},\)
\[
 \Lambda(v_i,-)= e_i\ne0,
 \qquad v_i=p^*(e_i).
\]
In particular, no unfrozen index of a compatible generalized seed is
removed by the active-index convention. An ineffective index can occur
only in the preceding abstract scattering datum, where its monomial is
central and its wall action is the identity.
\end{lemma}

\begin{proof}
Equation~\eqref{eq:seed-compatibility-normalized} is an equality of
covectors on \(M\):
\[
 \Lambda(v_i,-)=\langle e_i,-\rangle.
\]
Evaluating at \(e_i^*\in M\) gives
\(\Lambda(v_i,e_i^*)=1>0.\) Thus \(v_i\) is not in the radical
of \(\Lambda,\) and its annihilator is a genuine codimension-one
hyperplane. The last assertion is exactly the centrality observation
preceding Definition~\ref{def:GQSD}.
\end{proof}

Set
\[
  \Psi_{t}(X)
  :=\exp\left(
       \sum_{a\ge1}\frac{(-1)^{a-1}X^a}
       {a(t^{ a}-t^{- a})}
     \right).
\]
If
\[
  P_i(y)=\sum_{s=0}^{r_i} z_{i,s}y^s,
\]
define the incoming wall element by
\begin{align}
  \Phi_i(x^{v_i})
  &:=
  \exp\left(
    \sum_{a\geq1}
    \frac{(-1)^{a-1}s_{i,a}x^{av_i}}
         {a(t^{ a}-t^{-a})}
  \right) \notag\\
  &=\prod_{j=1}^{r_i}
    \Psi_{t}(q_{ij}x^{v_i})
  \qquad\text{after base change to }K_{\mathbf z},
  \label{eq:Delta-normalized-wall-element}
\end{align}
where \[ s_{i,a}:=\sum_{j=1}^{r_i}q_{ij}^{a} \] is the $a$-th power sum of the formal roots of $P_i$. 
Although the exponential expression in
\eqref{eq:Delta-normalized-wall-element} is singular at $t=1$, its
adjoint action is the regular finite product
\begin{equation}\label{eq:Delta-normalized-adjoint}
\begin{split}
  \Ad_{(\Phi_i)^\varepsilon}(x^m)
  =x^m\prod_{a=1}^{|\langle e_i,m\rangle|}
  P_i\!\left(t^{\varepsilon(2a-1)}x^{v_i}\right),
  \qquad
  \varepsilon=\sgn\langle e_i,m\rangle=\sgn\Lambda( v_i,m) .
\end{split}
\end{equation}
Indeed, the numerator arising from the commutator is divisible by
$t^{a}-t^{- a}$ because of
\eqref{eq:seed-compatibility-normalized}. Formula
\eqref{eq:Delta-normalized-adjoint} is the seed-scattering version of
the $\{a\}$-power in \cite[Section~3.1]{FPY24}.

Put
\[
 P_S:=\sum_{i=1}^n\ZZ_{\geq0}v_i,
 \qquad
 M_{S,0}^+:=P_S\setminus\{0\}.
\]
By Lemma~\ref{lem:compatible-pstar-injective}, every
\(v\in P_S\) has a unique expression \(v=\sum_i a_iv_i.\)  We use
the seed filtration
\begin{equation}\label{eq:seed-order}
 \ord_S(v):=\sum_i a_i,
 \qquad v=\sum_i a_iv_i,
\end{equation}
and retain its seed subscript throughout. Every
seed scattering group, quantum torus module, and completed base change
is the inverse limit of its finite \(\ord_S\)-quotients.
Let $M_{S,0}=\sum_i\ZZ v_i.$ Put
\begin{equation}\label{eq:seed-pairing-index}
 \gamma_S(v):=\gcd_{1\le i\le n}(a_i),
 \qquad v=\sum_i a_iv_i\in M_{S,0}^+,
 \qquad
 \widehat x^v:=\frac{x^v}{(\gamma_S(v))_t}.
\end{equation}
Define
\[
 \g_S:=\operatorname{Lie}_{\QQ_t}
       \langle\widehat x^v:v\in M_{S,0}^+\rangle
       \subseteq\bigoplus_{v\in M_{S,0}^+}\QQ(t)x^v,
 \qquad
 \g_{S,R_{\mathbf z}}
 :=\g_S\otimes_{\QQ_t}R_{\mathbf z}.
\]
Compatibility gives $\Lambda(v,M)=\gamma_S(v)\ZZ;$ in particular,
this index is preserved by the toric isometries used below.
Consequently each generator acts with Laurent coefficients:
\[
 [\widehat x^v,x^m]
 =\frac{(\Lambda(v,m))_t}{(\gamma_S(v))_t}x^{v+m}
 \in \QQ_t^\Lambda[M].
\]
The same holds for all iterated brackets. Moreover,
$\gamma_S(av_i)=a$, so
\[
 \log\Phi_i
 =\sum_{a\ge1}\frac{(-1)^{a-1}s_{i,a}}{a}\widehat x^{av_i}
 \in\widehat{\g}_{S,R_{\mathbf z}}.
\]

Let $\hat{G}_S$ denote the exponential group $\widehat{\g}_{S,R_{\mathbf z}}.$
The normalized classical limit of
Remark~\ref{rem:normalized-classical-limit} likewise applies, with
$\operatorname{cl}(\widehat x^v)=x^v/\gamma_S(v)$.

\begin{lemma}[Faithfulness of the seed scattering action]
\label{lem:seed-adjoint-faithful}
The adjoint representation of $\hat{G}_S$ on
all translated positive-order modules is faithful.
\end{lemma}

\begin{proof}

Let
\[
  0\ne v=\sum_{i=1}^n a_i v_i
\]
be a nonzero homogeneous degree of the seed scattering
Lie algebra.

Since \(v\ne0,\) there exists an index \(j\) such that
\(a_j\ne0.\) By the normalized compatibility condition
\eqref{eq:seed-compatibility-normalized},
\[
  \Lambda(v,e_j^*)
  =\sum_{i=1}^n a_i\Lambda(v_i,e_j^*)
  =a_j\ne0.
\]
Consequently,
\[
  [x^v,x^{e_j^*}]
  =
  \left(
    t^{a_j}
    -t^{-a_j}
  \right)x^{v+e_j^*}
  \ne0.
\]
Thus no nonzero homogeneous degree of the seed
scattering Lie algebra is central.

Now suppose that \(g=\exp(h)\) acts trivially.
Assume \(h\ne0,\) and let \(d\) be the lowest
nonzero order of \(h.\) Write
\[
  h_d=\sum_{\ord_S(v)=d}
        c_v x^v
\]
for its homogeneous component of order \(d\).

Since
\[
  \exp(\operatorname{ad}h)=\operatorname{id},
\]
comparison of the terms of order \(d\) gives
\[
  \operatorname{ad}(h_d)=0.
\]
In particular, for every \(m\in M,\)
\[
  0=[h_d,x^m]
   =\sum_{\ord_S(v)=d}
     c_v
     \left(
       t^{\Lambda(v,m)}
       -t^{-\Lambda(v,m)}
     \right)x^{v+m}.
\]

Distinct degrees \(v\) yield distinct monomials
\(x^{v+m}.\) Hence, for every \(v\) with \(c_v\ne0,\)
\[
  t^{\Lambda(v,m)}
  -t^{-\Lambda(v,m)}=0
  \qquad\text{for all }m\in M.
\]
Because \(t\) is formal, this implies
\[
  \Lambda(v,m)=0
  \qquad\text{for all }m\in M,
\]
contradicting the noncentrality established above.

Therefore \(h=0,\) and hence \(g=1.\)
The adjoint representation is faithful.

\end{proof}

\begin{definition}[Seed scattering diagram]
\label{def:seed-diagram}
The generalized quantum scattering diagram attached to the compatible
seed $S=(E_S,\widetilde B_S,\Lambda_S,\mathbf{r},\mathbf{z})$ is
\[
  \mathfrak D_{S,R_{\mathbf z}}
  :=\Scat(\mathfrak D_{S,\mathrm{in}}),
\]
where
\[
  \mathfrak D_{S,\mathrm{in}}
  =\left\{(e_i^\perp,\Phi_i(x^{v_i})):1\le i\le n\right\}.
\]
It is a diagram for $\g_{S,R_{\mathbf z}}$ in $M_\RR,$ acting on the completed quantum torus associated with $(M,\Lambda_S).$ By
Lemma~\ref{lem:compatible-seed-indices-active}, every displayed unfrozen direction is active, so no mutable seed index is deleted or reclassified. When the coefficient ring is clear, we write \(\mathfrak D_S\) for its
scalar extension.
\end{definition}

We record the mutation of the compatible seed, both to fix conventions and to keep the later reconstruction argument independent of any prior definition of cluster variables.
Put
\[
  R=\operatorname{diag}(r_1,\ldots,r_n).
\]
For $k\in I_{\uf}$ and $\varepsilon\in\{+1,-1\},$ let $E_{k,\varepsilon}^{\widetilde B R}\in\operatorname{Mat}_{m\times m}(\ZZ)$ be the identity away from its $k$th column, with
\begin{equation}\label{eq:E-mutation-matrix}
 \bigl(E_{k,\varepsilon}^{\widetilde B R}\bigr)_{ik}
 =\begin{cases}
   -1,&i=k,\\
   [-\varepsilon r_kb_{ik}]_+,&i\ne k,
  \end{cases}
\end{equation}
and let $F_{k,\varepsilon}^{R\widetilde B}\in\operatorname{Mat}_{n\times n}(\ZZ)$ be the identity away from its $k$th row, with
\begin{equation}\label{eq:F-mutation-matrix}
 \bigl(F_{k,\varepsilon}^{R\widetilde B}\bigr)_{kj}
 =\begin{cases}
   -1,&j=k,\\
   [\varepsilon r_kb_{kj}]_+,&j\ne k.
  \end{cases}
\end{equation}
\begin{remark}
The \(m\times m\) matrix \(E_{k,\varepsilon}^{\widetilde B R}\) acts on torus exponent columns. The \(n\times n\) matrix \(F_{k,\varepsilon}^{R\widetilde B}\) acts on the exchange-column indices; it is not a second coordinate change on \(M\). Only the positive sign specifies the forward basis convention.
\end{remark}

\begin{definition}[Mutation of compatible seed]
\label{def:compatible-data-mutation}
The mutation of \(S\) in direction \(k\in I_{\uf}\) is
\[
 S'=\mu_k(S)=(E',\widetilde B',\Lambda',\mathbf r,\mathbf z),
\]
where \(E'=\{e_i'\}_{i\in I}\) is the new \(\mathbb Z\)-basis of \(N\)
defined by
\begin{equation}\label{eq:seed-mutation-basis}
 e_i'=
 \begin{cases}
  e_i+r_k[-b_{ik}]_+e_k,&i\ne k,\\[2pt]
  -e_k,&i=k,
 \end{cases}
\end{equation}
with \(b_{ik}\) the \((i,k)\)-entry of \(\widetilde B\) and
\([x]_+=\max(x,0),\) and where 
\begin{equation}\label{eq:compatible-data-mutation}
 \widetilde B'=E_{k,+}^{\widetilde B R}\widetilde B F_{k,+}^{R\widetilde B},
 \qquad
 \Lambda'=(E_{k,+}^{\widetilde B R})^T\Lambda E_{k,+}^{\widetilde B R}.
\end{equation}
\end{definition}

\begin{remark}
The formulas for the mutated matrices in \eqref{eq:compatible-data-mutation} are independent of the choice of sign. That is, replacing the subscript \(+\) with \(-\) yields the same matrices:
\[
 E_{k,-}^{\widetilde B R}\widetilde B F_{k,-}^{R\widetilde B} = \widetilde B',
 \qquad
 (E_{k,-}^{\widetilde B R})^T\Lambda E_{k,-}^{\widetilde B R} = \Lambda'.
\]
Thus the definition of \(S'\) is well-defined with respect to the signed mutation matrices. (This is an instance of the sign-independence property established in Lemma~\ref{lem:T-preserves-Lambda}.)
\end{remark}

The following identity is the \(R\)-weighted generalized counterpart of the ordinary \(g\)-vector/\(c\)-vector tropical duality of \cite{nakanishi2023}.

\begin{lemma}[One-mutation symmetrizer identity]
\label{lem:one-mutation-symmetrizer-identity}
Let \(B=(b_{ij})_{1\leq i,j\leq n}\) be the principal part of \(\widetilde B,\) and let \(E_{k,\varepsilon}^{BR}\) be the \(n\times n\) principal block of \(E_{k,\varepsilon}^{\widetilde B R}.\) Then $B$ is skew-symmetric, and for each of the two mutation signs,
\begin{equation}\label{eq:EF-Langlands-intertwining}
 \boxed{
 \bigl(E_{k,\varepsilon}^{BR}\bigr)^T
 = F_{k,\varepsilon}^{RB}.}
\end{equation}
\end{lemma}

\begin{proof}
Multiplying the compatibility equation on the right by
\(\widetilde B\) gives
\[
 \widetilde B^T\Lambda\widetilde B= B.
\]
This equation implies that $B$ is skew-symmetric.

We verify \eqref{eq:EF-Langlands-intertwining} entry by entry.  Away
from row \(k\), both sides have the same rows.  At the
diagonal entry \((k,k)\), both sides equal \(-1.\) If \(j\ne k\),
then
\begin{align*}
 \bigl((E_{k,\varepsilon}^{BR})^T\bigr)_{kj}
 &=[-\varepsilon r_kb_{jk}]_+\\
 &=[ \varepsilon r_kb_{kj}]_+\\
 &=\bigl( F_{k,\varepsilon}^{RB}\bigr)_{kj}.
\end{align*}

\end{proof}

The usual entrywise calculation gives
\begin{equation}\label{eq:generalized-matrix-mutation-entrywise}
 b'_{ij}=\begin{cases}
 -b_{ij},&i=k\text{ or }j=k,\\
 b_{ij}+r_k\bigl([ - b_{ik}]_+b_{kj}
                 +b_{ik}[b_{kj}]_+\bigr),&i,j\ne k.
 \end{cases}
\end{equation}
In particular, the $k$th scattering direction satisfies
\begin{equation}\label{eq:vk-mutation}
  v'_k=p^*(e_{k}')=-v_k.
\end{equation}

\begin{theorem}[Wall-function factorization for seed diagrams]
\label{thm:seed-positive}
The equivalence class of \(\mathfrak D_{S,R_{\mathbf z}}\) has a
reduced orderwise \(R_{\mathbf z}\)-form. After base change to
\(K_{\mathbf z}\), it contains a reduced representative
\(\mathfrak D_S^+\) whose wall functions are products of factors
\[
  \E\left(
    -t^a\prod_{i,j}q_{ij}^{a_{ij}}x^m
  \right),
  \qquad a\in\ZZ,
  \quad a_{ij}\in\ZZ_{\ge0}.
\]
\end{theorem}

\begin{proof}
Start with the reduced Kontsevich--Soibelman representative over
\(R_{\mathbf z},\) whose existence follows order by order from Lemma~\ref{lem:reciprocal-coefficient-base-change}. This is the \(R_{\mathbf z}\)-form claimed in the statement.

Compatibility gives
\[
 \Lambda(v_i,m)=\langle e_i,m\rangle
 \qquad(i\in I_{\uf},\ m\in M).
\]
Thus the incoming walls of \(\mathfrak D_{S,R_{\mathbf z}}\) are
\[
 \bigl(e_i^\perp,\Phi_i(x^{v_i})\bigr),
 \qquad
 \Phi_i(x^{v_i})
 =\prod_{j=1}^{r_i}\Psi_t(q_{ij}x^{v_i}),
\]
and the diagram is defined over the uniform quantum parameter \(t\);
no componentwise rescaling is needed.  The relative-integrality
condition of Remark~\ref{rem:DM-relative-integrality} is automatic
because \(\Lambda\) is integral on all of \(M\times M\), and the
linear independence of the \(v_i\)'s is
Lemma~\ref{lem:compatible-pstar-injective}.

Applying Theorem~\ref{thm:positive-factorization} to
\(\mathfrak D_{S,R_{\mathbf z}}\), after base change to
\(K_{\mathbf z}\), gives an equivalent reduced representative whose
wall functions are products of factors
\[
 \E_t\left(
   -t^a\prod_{i=1}^{n}\prod_{j=1}^{r_i}q_{ij}^{a_{ij}}x^m
 \right),
 \qquad a\in\ZZ,\quad a_{ij}\ge0,\quad m\in L_0^+.
\]
The incoming walls of this representative are precisely the displayed
\(\Phi_i(x^{v_i})\), and all other walls are outgoing; this is the
representative \(\mathfrak D_S^+\).
\end{proof}

\subsection{Theta functions and positivity}
\begin{definition}[Pointed order sectors]
\label{def:pointed-seed-completion}
For \(p\in M\), put
\[
 \widehat{\mathcal T}_{S,p}
 :=\prod_{a\in P_S}R_{\mathbf z}x^{p+a},\qquad
 F_{S,p}^{\ge d}
 :=\prod_{\substack{a\in P_S\\\ord_S(a)\ge d}}
       R_{\mathbf z}x^{p+a}.
\]
This is a complete filtered module. Let
\(\mathcal T_S^{\mathrm{bd}}\) be the ring of formal sums whose
support is contained in a finite union of translates of \(P_S\).
Multiplication uses the based quantum-torus rule.
No ideal-adic completion of the full Laurent torus is intended.
\end{definition}

Every coefficient in a product in \(\mathcal T_S^{\mathrm{bd}}\)
is a finite sum: the order is positive on \(P_S^+\), and there
are only finitely many decompositions of a fixed positive degree.
Thus this ring contains the finite quantum torus and every pointed
sector, without imposing a common topology on all of \(M\).

\begin{lemma}[Separated pointed order completion]
\label{lem:pointed-order-completion}
Each \(\widehat{\mathcal T}_{S,p}\) is separated and complete.
If a linear operator \(R\) satisfies
\[
 R(F_{S,p}^{\ge d})\subseteq F_{S,p}^{\ge d+1},
\]
then \(1+R\) is a continuous module automorphism, with inverse
\(\sum_{j\ge0}(-R)^j\).
\end{lemma}

\begin{proof}
The order-\(N\) quotient has the finite basis
\(\{x^{p+a}:\ord_S(a)\le N\}\). The sector is the inverse
limit of these quotients, whose kernels have zero intersection.
On this quotient \(R^{N+1}=0\), so the geometric inverse is a
finite sum. The inverses commute with the quotient maps.
\end{proof}

\begin{lemma}[Orderwise finiteness of broken lines]
\label{lem:broken-line-order-finiteness}
Fix \(p\in M\), a general endpoint \(Q\), and \(N\ge0\).
Only finitely many broken lines with initial exponent \(p\) have
final exponent \(p+a\), with \(a\in P_S\) and
\(\ord_S(a)\le N\).
\end{lemma}

\begin{proof}
Use the finite order-\(N\) wall diagram and a finite rational
subdivision of its supports. Each genuine bend adds a nonzero
degree in \(P_S\), so at most \(N\) genuine bends occur.
There are finitely many lists of wall cells and surviving
homogeneous terms of these lengths. For any such list, reconstruct
the line backwards from \(Q\). The direction of each segment is
fixed by its monomial; its intersection with the next wall cell
is either empty or a single point. Generality excludes segments
contained in a supporting hyperplane and encounters with joints.
Consequently each list gives at most one broken line.
\end{proof}

\begin{proposition}[Pointed theta functions and positivity]
\label{prop:theta-positivity}
For \(p\in M\) and general \(Q\),
\[
 \vartheta_{p,Q}
 =x^p+\sum_{a\in P_S^+}c_{p,a}x^{p+a}
 \in\widehat{\mathcal T}_{S,p}.
\]
On each fixed sector \(\widehat{\mathcal T}_{S,p_0}\), the
linear transition \(x^{p_0+a}\mapsto\vartheta_{p_0+a,Q}\)
is a continuous filtration-preserving module automorphism.
After root extension, every \(c_{p,a}\) admits an expression in
\(\ZZ_{\ge0}[t^{\pm1},q_{ij}]\).
\end{proposition}

\begin{proof}
Lemma~\ref{lem:broken-line-order-finiteness} defines each coefficient
by a finite sum at its order. Normalized adjoint actions retain
the \(R_{\mathbf z}\)-coefficient form. The straight broken line
contributes \(x^p\) with coefficient \(1\); every genuine bend
adds a nonzero positive degree. Thus the transition on a fixed
sector is \(1+R\), with \(R\) strictly raising \(\ord_S\).
It is well defined on arbitrary sector series because below any
order bound only finitely many input and output degrees occur.
Lemma~\ref{lem:pointed-order-completion} gives its inverse.
The positive adjoint expansion used in
\cite[Lemma~3.8]{Davison_Mandel21}, applied to the free-root
representative of Theorem~\ref{thm:seed-positive}, gives the last
assertion by summing the finite broken-line contributions.
\end{proof}

This is a statement about a specified pointed module. It does not
define a full canonical theta algebra by completing its proposed
basis.

\begin{lemma}[Chamber constancy]
\label{lem:theta-chamber-constancy}
For two general endpoints \(Q,Q'\) in the same open chamber of a
consistent diagram,
\(\vartheta_{p,Q}=\vartheta_{p,Q'}\) in
\(\widehat{\mathcal T}_{S,p}\).
\end{lemma}

\begin{proof}
Join the endpoints inside the chamber. At every finite order
the path has identity wall product. The theta path-transport
identity of Section~\ref{subsec:broken-lines} therefore identifies
the two expansions. Pass to the inverse limit in the fixed sector.
\end{proof}

\begin{corollary}[Positive structure constants]
\label{cor:positive-structure}
For $p_1,p_2\in M$,
\[
  \vartheta_{p_1,Q}\vartheta_{p_2,Q}
  =t^{\Lambda(p_1,p_2)}\vartheta_{p_1+p_2,Q}
   +\sum_{v\in M_{S,0}^+}\alpha(p_1,p_2;p_1+p_2+v)\vartheta_{p_1+p_2+v,Q},
\]
and over \(K_{\mathbf z}\) the coefficients
\(\alpha(p_1,p_2;p_1+p_2+v)\) are Laurent polynomials in \(t\) and
polynomials in the lifted root variables \(q_{ij}\) with nonnegative
integer coefficients.
\end{corollary}

\begin{proof}
The sector transition in Proposition~\ref{prop:theta-positivity} first gives a unique convergent theta expansion in the sector with initial exponent \(p_1+p_2\). Theta path transport and uniqueness make its coefficients independent of the endpoint. For the coefficient indexed by \(p\), choose an auxiliary general endpoint sufficiently close to \(p\) and use the broken-line multiplication formula \cite[Proposition~3.3]{Davison_Mandel21}: the coefficient of \(\vartheta_{p,Q}\) is obtained by summing the ordered products of the final monomials of pairs of broken lines whose final exponents add to \(p.\) At a fixed target exponent this is a finite sum. By Theorem~\ref{thm:seed-positive} and the positive adjoint expansion in \cite[Lemma~3.8]{Davison_Mandel21}, every bend coefficient is a Laurent polynomial with nonnegative integral coefficients; multiplying two final based monomials only contributes a Laurent monomial in \(t.\) Hence every summand, and therefore every structure constant, is positive. This is the coefficient-wise version of \cite[Theorem~3.10]{Davison_Mandel21}; it uses the broken-line formula directly and does not invert the unitriangular monomial--theta transition matrix.

The unique pair of straight broken lines has final product
\[
 x^{p_1}x^{p_2}=t^{\Lambda(p_1,p_2)}x^{p_1+p_2},
\]
which gives the displayed leading coefficient. Any other pair has
acquired a nonzero element of \(M_{S,0}^+,\) giving precisely the
remaining indices in the statement.
\end{proof}
\section{Mutation of generalized quantum scattering diagrams}
\label{sec:mutation}

All seeds share the fixed lattices \(N\) and \(M=N^*,\) and the
fixed skew-symmetric form \(\Lambda\) on \(M.\) Seed subscripts
record the chosen bases. Fix \(k\in I_{\uf},\) put \(S'=\mu_k(S),\)
and abbreviate \(e_i=e_{i;S}\) and \(v_i=p_S^*(e_i)\).
The basis convention of Definition~\ref{def:compatible-data-mutation}
gives, as identities in the fixed lattice \(M\),
\begin{equation}\label{eq:fixed-lattice-mutated-degrees}
 v_k'=-v_k,\qquad
 v_i'=v_i+[r_kb_{ki}]_+v_k\quad(i\ne k).
\end{equation}
Thus
\[
 L_0:=\sum_{i=1}^n\mathbb Z v_i=\sum_{i=1}^n\mathbb Z v_i',\qquad
 P_S:=\bigoplus_{i=1}^n\mathbb Z_{\ge0}v_i,\qquad
 P_{S'}:=\bigoplus_{i=1}^n\mathbb Z_{\ge0}v_i'.
\]
Write \(P_S^+=P_S\setminus\{0\}\), and similarly for \(S'\).
The two diagrams lie in the same space \(M_{\mathbb R}\), but use
the respective positive monoids and their order completions.

\subsection{Linear branches and the exceptional wall}

Set
\[
 H_{k,\pm}:=\{m\in M_{\mathbb R}:
                  \pm\langle e_k,m\rangle\ge0\},\qquad
 A_k(m):=m+r_k\langle e_k,m\rangle v_k.
\]
Define
\begin{equation}\label{eq:fixed-lattice-tropical-mutation}
 T_{k,S}(m):=
 \begin{cases}
 A_k(m),&m\in H_{k,+},\\
 m,&m\in H_{k,-}.
 \end{cases}
\end{equation}
We write \(T_k=T_{k,S}\), \(T_{k,+}=A_k\), and
\(T_{k,-}=\mathrm{id}\). These linear branches agree on \(e_k^\perp\).

\begin{lemma}[Linear branches and duality]
\label{lem:explicit-contragredient-branches}
\label{lem:T-versus-E}
The map \(A_k\) is an integral lattice automorphism. In the old dual
basis,
\begin{equation}\label{eq:T-versus-E}
 [A_k]_{E_S^*}=E_{k,+}^{\widetilde B_SR}
                    E_{k,-}^{\widetilde B_SR}.
\end{equation}
Moreover,
\[
 A_k^{-1}=I-r_kv_k\otimes e_k,\qquad
 A_k^{-T}=I-r_ke_k\otimes v_k,
 \qquad
 \langle A_k^{-T}n,A_km\rangle=\langle n,m\rangle.
\]
\end{lemma}

\begin{proof}
Since \(\langle e_k,v_k\rangle=b_{kk}=0,\) the operator \(v_k\otimes e_k\) has square zero. This proves the inverse formulas and the pairing identity. Both mutation matrices fix \(\mathbf e_i\) for \(i\ne k,\) while their product sends \(\mathbf e_k\) to \(\mathbf e_k+r_k\mathbf b_k.\) This proves \eqref{eq:T-versus-E}.
\end{proof}

\begin{lemma}[Preservation of the quantum form]
\label{lem:T-preserves-Lambda}
Both branches \(T_{k,\pm}\) preserve \(\Lambda.\)
\end{lemma}

\begin{proof}
For \(a_i=\langle e_k,m_i\rangle,\) compatibility gives
\[
 \Lambda(A_km_1,A_km_2)
 =\Lambda(m_1,m_2)-r_k a_1a_2
                       +r_k a_1a_2
 =\Lambda(m_1,m_2).
\]
The other branch is the identity.
\end{proof}

The branches therefore act on the based quantum torus by \(T_{k,\sigma}(x^m)=x^{T_{k,\sigma}m}.\) They also reindex formal coefficient families indexed by \(M.\) This reindexing does not, by itself, assert membership in the \(P_{S'}\)-completion.

\begin{lemma}[Inverse tropical transport]
\label{lem:T-seed-inverse}
The map \(T_{k,S}\) is a bijection, with inverse
\[
 T_{k,S}^{-1}(m)=
 \begin{cases}
 m-r_k\langle e_k,m\rangle v_k,&m\in H_{k,+},\\
 m,&m\in H_{k,-}.
 \end{cases}
\]
If the same forward basis mutation is applied at \(S',\) then
\begin{equation}\label{eq:fixed-lattice-double-transport}
 T_{k,S'}T_{k,S}=A_k,\qquad
 T_{k,S}^{-1}=A_k^{-1}T_{k,S'}.
\end{equation}
\end{lemma}

\begin{proof}
Both branches preserve \(\langle e_k,m\rangle,\) so the displayed inverse cancels the shear on \(H_{k,+}.\) At \(S',\) one has \(e_k'=-e_k\) and \(v_k'=-v_k.\) Thus \(T_{k,S'}\) is \(A_k\) on the old \(H_{k,-}\) and the identity on the old \(H_{k,+},\) proving \eqref{eq:fixed-lattice-double-transport}.
\end{proof}

\begin{remark}\label{rem:fixed-lattice-double-basis}
Two forward basis mutations need not restore the basis in the fixed lattice. Indeed,
\[
 e_k''=e_k,\qquad e_i''=e_i-r_kb_{ik}e_k\quad(i\ne k),
 \qquad (e_i'')^*=A_k(e_i^*).
\]
Consequently the reverse transport uses \(T_{k,S}^{-1}.\) If it is expressed using a second forward mutation, the linear identification \(A_k^{-1}\) must be included.
\end{remark}

\paragraph{Exceptional-wall localization.}
Let \(\mathcal T_{R_{\mathbf z}}=R_{\mathbf z}^{\Lambda}[M]\). Let \(S_k\) be the multiplicative set generated by
\[
 P_k(t^{(2a-1)}x^{v_k}),\qquad a\in\mathbb Z.
\]
For every \(m\in M,\)
\begin{equation}\label{eq:shifted-exchange-normality}
 x^mP_k(t^c x^{v_k})
 =P_k(t^{c+2\Lambda(m,v_k)}x^{v_k})x^m.
\end{equation}
The shift is a multiple of \(2,\) so these conjugations preserve \(S_k.\) Since its elements commute with one another, \(S_k\) is a left and right Ore set. Put
\[
 \mathcal T_{S,k}^{\mathrm{loc}}
 :=S_k^{-1}\mathcal T_{R_{\mathbf z}}.
\]
Reciprocity makes the corresponding factors in \(x^{-v_k}\) invertible in the same localization. The following identities concern finite rational adjoint actions, not evaluation of a singular dilogarithm.

\begin{lemma}[Exceptional-wall transition]
\label{lem:exceptional-conjugation}
On \(\mathcal T_{S,k}^{\mathrm{loc}},\)
\begin{equation}\label{eq:exceptional-fixed-lattice-transition}
 T_{k,+}^{-1}\circ
 \operatorname{Ad}_{\Phi_k(x^{-v_k})}\circ T_{k,-}
 =
 \operatorname{Ad}_{\Phi_k(x^{v_k})^{-1}}.
\end{equation}
\end{lemma}

\begin{proof}
Write \(v=v_k,\) \(r=r_k,\) and \(c=\langle e_k,m\rangle.\) For \(c>0,\) the finite adjoint formula gives
\[
 \operatorname{Ad}_{\Phi_k(x^{-v})}(x^m)
 =x^m\prod_{a=1}^{c}P_k(t^{-(2a-1)}x^{-v})^{-1}.
\]
For \(D_a=(2a-1),\) reciprocity gives
\[
 P_k(t^{-D_a}x^{-v})^{-1}
 =t^{rD_a}x^{rv}P_k(t^{D_a}x^v)^{-1}.
\]
The scalar \(t^{rc^2}\) cancels the scalar \(t^{-rc^2}\) from \(x^mx^{rcv}.\) Hence the preceding
adjoint expression is
\[
 x^{m+rcv}\prod_{a=1}^{c}P_k(t^{(2a-1)}x^v)^{-1}.
\]
Applying \(A_k^{-1},\) which fixes \(v,\) gives the required
inverse adjoint action.

If \(c<0,\) put \(d=-c.\) Now
\[
 \operatorname{Ad}_{\Phi_k(x^{-v})}(x^m)
 =x^m\prod_{a=1}^{d}P_k(t^{(2a-1)}x^{-v}).
\]
Use
\(P_k(t^{D_a}x^{-v})=t^{rD_a}x^{-rv}P_k(t^{-D_a}x^v).\)
The same cancellation gives leading exponent \(m-rdv=m+rcv\);
applying \(A_k^{-1}\) leaves
\[
 x^m\prod_{a=1}^{d}P_k(t^{-(2a-1)}x^v),
\]
as required. The case \(c=0\) is immediate.
\end{proof}

For a crossing from \(H_{k,-}\) to \(H_{k,+},\) the Davison--Mandel crossing exponents of the old direction \(v_k\) and the new direction \(-v_k\) are respectively \(-1\) and \(+1.\) Equation~\eqref{eq:exceptional-fixed-lattice-transition} is precisely the corresponding transition; its inverse handles the reverse crossing.

\subsection{Classical detection and the mutation transform}

We choose the positive representative \(\mathfrak D_S^+\) of Theorem~\ref{thm:seed-positive}, retaining all full incoming walls and only outgoing additional walls. Reduction will be used for orderwise local comparisons.

There are no additional positive factors of pure degree \(\ell v_k\) in this representative. Indeed, quotient the graded Lie algebra by degrees involving some \(v_i,\) \(i\ne k.\) Only the commuting \(k\)-th incoming data remain, whose consistent completion is the single incoming wall. At every general local wall point, equivalence and Lemma~\ref{lem:raywise-bps-uniqueness} leave no additional positive pure-\(k\) factor after the prescribed incoming factor has been removed.

\begin{definition}[Piecewise-linear wall transform]
\label{def:mutated-diagram}
For each wall \((\mathfrak d,f_{\mathfrak d})\) other than the \(k\)-th incoming wall, insert the nonempty codimension-one pieces
\[
 \left(T_{k,\sigma}(\mathfrak d\cap H_{k,\sigma}),
       T_{k,\sigma}(f_{\mathfrak d})\right),
 \qquad \sigma\in\{+,-\}.
\]
Replace the \(k\)-th incoming wall by
\[
 (e_k^\perp,\Phi_k(x^{-v_k})).
\]
Denote this wall collection by \(T_k(\mathfrak D_S^+).\) Its membership in the target completion is established below.
\end{definition}

For a diagram \(\mathfrak C,\) let \(\EssSupp(\mathfrak C)\) record its general oriented wall points with nonidentity normalized local product, order by order.

\begin{lemma}[Oriented support is an equivalence invariant]
\label{lem:oriented-support-invariant}
Equivalent diagrams without ineffective factors have the same normalized oriented local products, and hence the same essential oriented support.
\end{lemma}

\begin{proof}
Fix a general point \(x\) of a supporting hyperplane \(H,\) an oriented ray \(\rho\subset H^\perp,\) and an order bound \(N.\) Apply Lemma~\ref{lem:compatible-generic-germs} to a common finite subdivision near \(x.\) It identifies the normalized \(\rho\)-products. The identities commute with truncation and therefore pass to the inverse limit. Reduction only removes identity local products and changes their polyhedral presentation.
\end{proof}

\begin{lemma}[Nonvanishing of positive reciprocal products]
\label{lem:positive-reciprocal-product-nonvanishing}
Let \(A=\mathbb Z[X_1^{\pm1},\ldots,X_s^{\pm1}],\) and let \(P\) be a finitely generated sharp monoid with a positive additive order. An orderwise finite nonempty product
\[
 \prod_\nu(1+X^{a_\nu}Z^{u_\nu})^{c_\nu}
 \prod_\mu(1-X^{b_\mu}Z^{w_\mu})^{-d_\mu},
 \qquad
 u_\nu,w_\mu\in P^+,\quad c_\nu,d_\mu\in\mathbb Z_{>0},
\]
is not \(1\) in \(\widehat{A[P]},\) nor after extending coefficients
to \(\operatorname{Frac}(A).\) Its first nonzero order does not cancel.
\end{lemma}

\begin{proof}
Let \(r\) be the least order of a degree appearing in a factor.
Modulo terms of order greater than \(r,\) the product is
\[
 1+\sum_{\operatorname{ord}(u_\nu)=r}
      c_\nu X^{a_\nu}Z^{u_\nu}
  +\sum_{\operatorname{ord}(w_\mu)=r}
      d_\mu X^{b_\mu}Z^{w_\mu}.
\]
Collecting equal degrees and Laurent monomials only adds positive
integers. The nonconstant part is nonzero, and the coefficient
extension is injective.
\end{proof}

Let \(\mathfrak D_{\mathrm{cl},S}\) be a classical generalized scattering representative over \(\mathbb F_{\mathbf z}^{\mathrm{cl}},\) with the full incoming walls and outgoing additional walls. Its incoming actions are
\begin{equation}\label{eq:classical-generalized-incoming-action}
 x^m\longmapsto
 x^mP_i(x^{v_i})^{\langle e_i,m\rangle}.
\end{equation}

\begin{lemma}[Classical detection of positive walls]
\label{lem:classical-support-detection}
Regular adjoint specialization at \(t=1\) identifies \(\mathfrak D_S^+\) with a representative of \(\mathfrak D_{\mathrm{cl},S}\). No nontrivial normalized oriented local product disappears, including its first nonzero order. In particular,
\[
 \EssSupp(\mathfrak D_S^+)
 =
 \EssSupp(\mathfrak D_{\mathrm{cl},S}).
\]
\end{lemma}

\begin{proof}
Put \(\ell_i=\lfloor r_i/2\rfloor,\) choose algebraically independent variables \(X_{i,a},\) and set
\[
 K_{\mathrm{rt}}=\mathbb Q(X_{i,a}:i,\ 1\le a\le\ell_i).
\]
Use the reciprocal root parameters
\[
 (q_{i1},\ldots,q_{ir_i})=
 \begin{cases}
 (X_{i,1},X_{i,1}^{-1},\ldots,X_{i,\ell_i},X_{i,\ell_i}^{-1}),
       &r_i=2\ell_i,\\
 (X_{i,1},X_{i,1}^{-1},\ldots,X_{i,\ell_i},X_{i,\ell_i}^{-1},1),
       &r_i=2\ell_i+1.
 \end{cases}
\]
For \(Y_{i,a}=X_{i,a}+X_{i,a}^{-1},\) the polynomial \(P_i(Y)\) is \(\prod_a(1+Y_{i,a}Y+Y^2),\) with an additional factor \(1+Y\) in odd degree. If \(E_{i,s}\) denotes the \(s\)-th elementary symmetric
polynomial in the \(Y_{i,a}\), coefficient comparison gives
\[
 z_{i,s}=E_{i,s}+Q_{i,s}(E_{i,1},\ldots,E_{i,s-1}),
 \qquad 1\le s\le\ell_i.
\]
The variables \(Y_{i,a}\) are algebraically independent: the extension from their field to \(K_{\mathrm{rt}}\) is algebraic, since \(X_{i,a}^2-Y_{i,a}X_{i,a}+1=0,\) and both fields have the same transcendence degree. The triangular identities therefore embed \(R_{\mathbf z}^{\mathrm{cl}}\) in \(K_{\mathrm{rt}}.\) Moreover, each \(Y_{i,a}\) is algebraic over the field of its elementary symmetric polynomials. Hence \(K_{\mathrm{rt}}/\mathbb F_{\mathbf z}^{\mathrm{cl}}\) is finite and faithfully flat, as is \(K_{\mathrm{rt}}(t)/\mathbb F_{\mathbf z}.\)

 A positive elementary factor has regular adjoint specialization
\begin{equation}\label{eq:positive-factor-classical-action}
 \operatorname{sp}_{t=1}
 \operatorname{Ad}_{\E_{t}
       (-t^a\mathbf q^\alpha x^v)}(x^m)
 =
 x^m(1+(-1)^a\mathbf q^\alpha x^v)^{
                         (-1)^a\Lambda(v,m)}.
\end{equation}

Set \(m_{0,S}=\sum_{i=1}^m e_{i;S}^*.\) For a nonzero degree \(v=\sum_{1\le i\le n}a_iv_i,\) one has
\[
 \Lambda(v,m_{0,S})=\sum_{1\le i\le n}a_i>0.
\]
At a general wall point, all factors on a fixed effective covector ray commute and fix the monomials occurring in the other factors. Consequently their normalized action on \(x^{m_{0,S}},\) divided by \(x^{m_{0,S}},\) is a product of \((1+\mathbf q^\alpha x^v)^c\) and \((1-\mathbf q^\alpha x^v)^{-c},\) with \(c>0.\) After reciprocal substitution this is covered by Lemma~\ref{lem:positive-reciprocal-product-nonvanishing}. A nonempty positive local product therefore remains nontrivial, and its least BPS order is still visible. The normalized leading-term
recursion, or the faithful seed action, identifies this order with the first nonzero order of the local wall action.

For incoming factors, \eqref{eq:positive-factor-classical-action} gives exactly \eqref{eq:classical-generalized-incoming-action}. Regular adjoint specialization commutes with finite-order products and with the inverse limit, by Remarks~\ref{rem:no-specialization-wall-element} and~\ref{rem:normalized-classical-limit}. The specialized diagram is consistent, retains its incoming data, and has only outgoing additional walls. The classical action is faithful by the same injectivity calculation as Lemma~\ref{lem:seed-adjoint-faithful}. Classical scattering uniqueness therefore identifies its equivalence class with \(\mathfrak D_{\mathrm{cl},S},\) over \(K_{\mathrm{rt}}.\) The preceding nonvanishing and Lemma~\ref{lem:oriented-support-invariant} identify the essential oriented supports.

Finally, the canonical finite-order local actions come from the universal normalized coefficient form. The injection \(R_{\mathbf z}^{\mathrm{cl}}\hookrightarrow K_{\mathrm{rt}}\) allows equality and nontriviality of these finite expressions to descend. No evaluation on all of \(\mathbb Q(t),\) and no positivity in the elementary symmetric variables, is asserted.
\end{proof}

\begin{lemma}[Order control under tropical transport]
\label{lem:tropical-order-control}
Let \(B\) be either branch of \(T_{k,S}\), and fix \(h\in L_0\).
For every \(N\ge0\), there is \(N'(N,h)\) such that
\[
 a\in P_S,\quad Ba+h\in P_{S'},\quad
 \ord_{S'}(Ba+h)\le N
 \quad\Longrightarrow\quad
 \ord_S(a)\le N'(N,h).
\]
The bound can be chosen for both branches simultaneously and is
independent of the coefficient values.
\end{lemma}

\begin{proof}
The set \(Q_N=\{b\in P_{S'}:\ord_{S'}(b)\le N\}\) is finite.
Each integral invertible branch has at most one preimage of each
\(b-h\). Take the maximum source order over the preimages lying
in \(P_S\), for both branches, and include \(0\) in this maximum.
\end{proof}

Once target positivity of the degrees has been proved, the case
\(h=0\) controls wall products and their truncations. For a theta
series with initial exponent \(p\), use
\(h=Bp-T_{k,S}(p)\): a final exponent \(p+a\) then has target
increment \(Ba+h\). In this use \(h\in L_0\), because
\(Bp-p\) and \(T_{k,S}(p)-p\) are multiples of \(v_k\).
This is order control for the expressions under consideration,
not a continuous identification of the two entire completions.

\begin{proposition}[Well-defined mutation transform]
\label{prop:mutated-diagram-well-defined}
The collection \(T_k(\mathfrak D_S^+)\) is an orderwise finite
scattering diagram in the \(P_{S'}\)-completion.
\end{proposition}

\begin{proof}
By \eqref{eq:fixed-lattice-mutated-degrees} and the unimodularity of \(F_{k,+}^{RB_S},\) the \(v_i'\) form a basis of \(L_0.\) Thus \(P_{S'}\) is saturated in \(L_0\). Both branches of \(T_k\) preserve \(L_0\), since
\[
 A_k(v_i)=v_i+r_kb_{ki}v_k
\]
and the inverse shear has the same integrality property.

We first work with finite local expressions, before claiming membership in the target completion. The classical generalized mutation theorem, expressed in our fixed-lattice basis convention, gives
\begin{equation}\label{eq:classical-fixed-lattice-mutation}
 T_k(\mathfrak D_{\mathrm{cl},S})
 \sim\mathfrak D_{\mathrm{cl},S'}.
\end{equation}
See \cite[Theorem~3.25]{CKM} or \cite[Theorem~6.27]{Mou24}. Under \eqref{eq:forward-seed-coordinate-change}, a target coordinate column \(a\) represents \(\sum_j(E_{k,+}a)_j e_{j;S}^*\) in \(M\). Thus the marked target diagram in that coordinate formulation is precisely \(\mathfrak D_{\mathrm{cl},S'}\) in the fixed space used here.

Consider a general nonexceptional transformed wall point, in the interior of one of \(H_{k,\pm},\) and an actual degree \(d\) on one of its effective covector rays. Reindexing by that branch preserves the local adjoint expressions and commutes with their classical
specialization. Lemma~\ref{lem:classical-support-detection} makes the classical normalized ray product nontrivial. By \eqref{eq:classical-fixed-lattice-mutation}, the target classical local product has a nonzero logarithmic term of some degree \(c\in P_{S'}^+\) on the same effective ray. Hence
\[
 \mathbb R_{>0}\Lambda(d,-)=\mathbb R_{>0}\Lambda(c,-).
\]
Compatibility makes \(\Lambda:L_0\otimes\mathbb Q\to N\otimes\mathbb Q\) injective. Thus \(d=\lambda c\) with \(\lambda\in\mathbb Q_{>0};\) saturation gives \(d\in P_{S'}^+.\) The exceptional wall has degrees in \(\mathbb Z_{>0}(-v_k)\subset P_{S'}^+\). The pure-\(k\) observation above excludes any other wall with supporting hyperplane \(e_k^\perp\).

Now apply Lemma~\ref{lem:tropical-order-control} with \(h=0\). Every target order bound uses a finite common source-order bound. Source-order finiteness and the finite half-space subdivision give target-order finiteness and convergence of every transformed wall element. Target positivity has already been proved, so this establishes the target completion before any global target path product is used.
\end{proof}

\subsection{Consistency and incoming data}

\begin{lemma}[Linear transport on a half-space]
\label{lem:linear-halfspace-transport}
Let \(\gamma\) be an admissible path contained in \(H_{k,\sigma}^{\circ},\) where $\sigma\in \{1,-1\}.$ Then
\[
 \PathProd_{T_{k,\sigma}\gamma,T_k(\mathfrak D_S^+)}
 =
 T_{k,\sigma}\bigl(
   \PathProd_{\gamma,\mathfrak D_S^+}\bigr).
\]
Equivalently,
\[
 \WC_{T_{k,\sigma}\gamma,T_k(\mathfrak D_S^+)}
 \circ T_{k,\sigma}
 =
 T_{k,\sigma}\circ\WC_{\gamma,\mathfrak D_S^+}.
\]
\end{lemma}

\begin{proof}
For each target order, the preceding proposition supplies a finite source-order bound. At those bounds there are finitely many crossing factors. Their order is unchanged, and
\[
 \Lambda(T_{k,\sigma}v,-T_{k,\sigma}\dot\gamma)
 =\Lambda(v,-\dot\gamma)
\]
preserves the crossing signs. Multiply the elementary identities and
pass to the inverse limit.
\end{proof}

\begin{lemma}[Incoming walls under mutation]
\label{lem:incoming-data-mutation}
After a common half-space subdivision, \(T_k(\mathfrak D_S^+)\) and \(\mathfrak D_{S'}\) have the same incoming walls and wall elements. Every remaining wall is outgoing.
\end{lemma}

\begin{proof}
For \(i\ne k\) with \(b_{ki}\ne0,\) set \(\sigma_i=\operatorname{sgn}(b_{ki})\). Equation~\eqref{eq:fixed-lattice-mutated-degrees} gives \(T_{k,\sigma_i}(v_i)=v_i'.\) Compatibility implies \(b_{ik}=-b_{ki}.\) Using \(e_i'=e_i+r_k[-b_{ik}]_+e_k,\) one checks
\[
 T_{k,\sigma_i}(e_i^\perp\cap H_{k,\sigma_i})
 =(e_i')^\perp\cap H_{k,\sigma_i}.
\]
For example, if \(b_{ki}>0,\) then \(b_{ik}<0,\) and for \(m\in e_i^\perp\cap H_{k,+},\)
\[
 \langle e_i',A_km\rangle
 =\langle e_i,m\rangle+
   r_k\bigl(b_{ik}+[-b_{ik}]_+\bigr)\langle e_k,m\rangle=0.
\]
The negative case uses the identity branch and \(b_{ik}>0.\) The unchanged coefficient data and the toric isometry give
\[
 T_{k,\sigma_i}(\Phi_i(x^{v_i}))
 =\Phi_i(x^{v_i'}).
\]
The complementary half \(H_{k,-\sigma_i}\) does not contain \(v_i\), and its image is outgoing. If \(b_{ki}=0,\) both halves have the same direction and wall element and recombine. For \(i=k,\) one has \((e_k')^\perp=e_k^\perp\) and \(v_k'=-v_k,\) exactly the replacement wall.

Finally, for any other conical wall and any branch,
\[
 T_{k,\sigma}(v)\in
 T_{k,\sigma}(\mathfrak d\cap H_{k,\sigma})
 \quad\Longleftrightarrow\quad
 v\in\mathfrak d\cap H_{k,\sigma}.
\]
Thus outgoing walls remain outgoing. Subdivide the target initial walls in the same way; their complementary halves are also outgoing. This identifies the incoming data used for uniqueness.
\end{proof}

\begin{theorem}[Mutation compatibility]
\label{thm:mutation}
Let \(S'=\mu_k(S).\) As diagrams in the fixed space \(M_{\mathbb R}\),
\[
 T_{k,S}(\mathfrak D_S^+)\sim\mathfrak D_{S'}.
\]
\end{theorem}

\begin{proof}
Proposition~\ref{prop:mutated-diagram-well-defined} places the transformed collection in the target sharp completion. Fix a target order and the corresponding finite source-order bound. Choose a closed path \(\gamma,\) based in \(H_{k,\sigma_0}^{\circ},\) generic for the relevant finite diagrams and transverse to \(e_k^\perp.\) Pull its half-space arcs back by the appropriate branches of \(T_k^{-1}.\) They concatenate to the closed path \(\widetilde\gamma=T_k^{-1}\gamma;\) a small perturbation of its corners does not change the crossing data.

Lemma~\ref{lem:linear-halfspace-transport} compares the open arcs. At a crossing from \(H_{k,-}\) to \(H_{k,+}\), the old and new crossing exponents are \(-1\) and \(+1\). Lemma~\ref{lem:exceptional-conjugation} gives the branch transition; its inverse treats the opposite crossing. The intermediate linear branch maps cancel, giving
\begin{equation}\label{eq:mutation-loop-intertwining}
 T_{k,\sigma_0}^{-1}\circ
 \WC_{\gamma,T_k(\mathfrak D_S^+)}
 \circ T_{k,\sigma_0}
 =
 \WC_{\widetilde\gamma,\mathfrak D_S^+}
 =\mathrm{id}.
\end{equation}
Faithfulness for the target seed, Lemma~\ref{lem:seed-adjoint-faithful}, makes the target path-group element trivial. Passing through the finite orders proves consistency.

Lemma~\ref{lem:incoming-data-mutation} identifies the common incoming data after subdivision, with every other wall outgoing. Theorem~\ref{thm:KS-uniqueness} gives the required equivalence. All comparisons of canonical truncated products descend through the faithful reciprocal root extension.
\end{proof}

\subsection{Broken lines and theta functions}

Use \(T_{k,S}(\mathfrak D_S^+)\) as the specified representative of the target equivalence class, and denote its theta sums by \(\vartheta^{S',T}_{p,Q}.\)

\begin{lemma}[Exceptional bending correspondence]
\label{lem:exceptional-bending-correspondence}
At a transverse crossing of \(e_k^\perp\), the two branches of
\(T_{k,S}\) identify the collected monomial terms in the old
and transformed incoming-wall bending rules, with equal
coefficients. The correspondence includes the constant terms.
\end{lemma}

\begin{proof}
Put \(v=v_k\), \(r=r_k\), \(Y=x^v\), and
\(D_a=2a-1\). Write the monomial before the crossing as
\(b x^m\), with \(c=\langle e_k,m\rangle\ne0\), and put
\(d=|c|\). Since the velocity is \(-m\), a crossing from
\(H_{k,-}\) to \(H_{k,+}\) has \(c=-d<0\).
The old crossing exponent is \(-1\), and its finite expansion is
\[
 b x^m\prod_{a=1}^dP_k(t^{-D_a}Y).
\]
Choose \(s_a\in\{0,\ldots,r\}\) in the \(a\)-th factor and
write \(s=\sum_a s_a\). Quantum-torus multiplication makes the
chosen based term
\[
 b\left(\prod_a z_{k,s_a}\right)
 t^{\sum_a(d-D_a)s_a}x^{m+sv}.
\]
The outgoing branch sends its exponent to
\(A_k(m+sv)=m+(s-rd)v\).
In the transformed diagram the incoming monomial is still
\(b x^m\), the crossing exponent for \(-v\) is \(+1\), and
the expansion is
\[
 b x^m\prod_{a=1}^dP_k(t^{D_a}Y^{-1}).
\]
Choose \(s'_a=r-s_a\). Its exponent is
\(m-\sum_as'_av=m+(s-rd)v\), and its coefficient is
\[
 b\left(\prod_a z_{k,r-s_a}\right)
 t^{\sum_a(D_a-d)(r-s_a)}.
\]
This is the preceding coefficient: reciprocity gives
\(z_{k,r-s_a}=z_{k,s_a}\), while
\(\sum_aD_a=d^2\) cancels the constant part of the exponent.

For the reverse crossing, \(c=d>0\). The old expansion is
\(b x^m\prod_aP_k(t^{D_a}Y)\), whose tuple has coefficient
\(b\prod_a z_{k,s_a}\,t^{\sum_a(D_a-d)s_a}\).
The transformed incoming exponent is \(A_km=m+rdv\), and its
expansion is
\[
 b x^{m+rdv}\prod_{a=1}^dP_k(t^{-D_a}Y^{-1}).
\]
The same complementary tuple has outgoing exponent \(m+sv\)
and Weyl exponent
\(\sum_a(d-D_a)(r-s_a)=\sum_a(D_a-d)s_a\).
Thus it has exactly the old coefficient.

The tuple complement is an involution. Tuples with the same sum
are collected into one monomial coefficient; their complements
also have the same sum \(rd-s\), so the equality survives this
collection. A vanished evaluated coefficient vanishes on both
sides. A straight crossing on one side may become a genuine bend
on the other; allowing the constant terms accounts for this case.
\end{proof}

\begin{proposition}[Mutation of broken lines]
\label{prop:broken-lines-mutation}
The map \(T_{k,S}\) gives a coefficient-preserving bijection from broken lines for \(\mathfrak D_S^+\) to those for \(T_{k,S}(\mathfrak D_S^+),\) sending their initial exponent and endpoint from \((p,Q)\) to \((T_{k,S}(p),T_{k,S}(Q)).\) If \(Q\in H_{k,\sigma}^{\circ},\) a final monomial \(cx^m\) is sent to \(cx^{T_{k,\sigma}m}.\)
\end{proposition}

\begin{proof}
Subdivide a broken line at its crossings of \(e_k^\perp.\) On a segment in \(H_{k,\sigma},\) replace its position by \(T_{k,\sigma}\gamma\) and its monomial \(cx^m\) by \(cx^{T_{k,\sigma}m}.\) The branches agree on the separating hyperplane and preserve its defining coordinate; the resulting path
is continuous and has velocity \(-T_{k,\sigma}m.\) Its unbounded initial segment lies in the half-space determined by \(\langle e_k,p\rangle,\) so its initial exponent is \(T_k(p).\)

Away from \(e_k^\perp\), the branch is an integral isometry.
It reindexes each homogeneous bending term without altering its
coefficient or crossing sign. At \(e_k^\perp\), apply
Lemma~\ref{lem:exceptional-bending-correspondence}, including its
constant terms. Its complementary choices verify the bending
rule for the transformed path, not merely equality of the total
wall automorphisms. Consecutive segments with the same monomial
are merged, so this also treats a bend that becomes straight.

Apply the inverse branches of \(T_{k,S}\) to reverse the construction.
Lemma~\ref{lem:T-seed-inverse} supplies the actual inverse, including the required linear correction if it is expressed as a second forward
mutation. The two constructions are inverse on individual broken lines. For an order-\(N\) target theta sum, the final exponents use only a finite source-order truncation by the affine case of Lemma~\ref{lem:tropical-order-control}. Thus the correspondence also identifies the orderwise sums.
\end{proof}

\begin{corollary}[Theta functions under mutation]
\label{cor:theta-mutation}
For general \(Q\in H_{k,\sigma}^{\circ},\)
\[
 T_{k,\sigma}(\vartheta^S_{p,Q})
 =
 \vartheta^{S',T}_{T_{k,S}(p),\,T_{k,S}(Q)}.
\]
The right-hand side is the theta function of the target equivalence
class at the indicated exponent and endpoint.
\end{corollary}

\begin{proof}
Sum the coefficient-preserving correspondence of Proposition~\ref{prop:broken-lines-mutation} at each finite order. Equivalence of consistent diagrams preserves theta sums.
\end{proof}

\begin{corollary}[Inverse theta transport]
\label{def:theta-transports-mutation} \label{cor:theta-transport-involutive}The correspondence of Corollary~\ref{cor:theta-mutation} and the correspondence defined by \(T_{k,S}^{-1}\) are mutually inverse on individual theta functions and their finite linear combinations.
\end{corollary}

\begin{proof}
The inverse broken-line construction restores both the initial
exponent and the endpoint, and preserves every coefficient. Linearity gives the assertion for finite sums.
\end{proof}

\begin{remark}\label{rem:theta-transport-completions}
These statements do not identify the full \(P_S\)- and \(P_{S'}\)-completions by continuity. Their cones, and hence their topologies, can differ even though the ambient lattice is fixed. Lemma~\ref{lem:tropical-order-control} supplies finite-order bounds for the transported expressions, not continuity of an unrestricted map between all series. The algebraic chart transitions below use finite rational wall actions; a stronger statement about canonical theta algebras requires a separate comparison of the relevant spaces of series.
\end{remark}

\section{Quantum generalized cluster algebras from scattering diagrams}
\label{sec:scattering-algebra}

We define the quantum generalized cluster algebras from the rational
chart transitions of the generalized quantum scattering diagram.
We then identify its classical seed pattern and reachable cluster
complex, express cluster monomials as theta functions, and prove the Laurent phenomenon and root positivity. Throughout, we retain the coefficient form
\[
 R_{\mathbf z}
 =\QQ[t^{\pm1}]\otimes_{\QQ}R_{\mathbf z}^{\mathrm{cl}}.
\]

\subsection{Quantum torus charts and elementary scattering pullbacks}
\label{sec:quantum-torus-charts}

All seeds use the fixed lattices \(N,\) \(M=N^*,\) and the fixed
skew-symmetric form \(\Lambda\) of Section~\ref{sec:mutation}. For a seed \(S_u\), define
\begin{equation}\label{eq:seed-coordinate-identification}
 \beta_u:\ZZ^m\longrightarrow M,\qquad
 \beta_u(\mathbf a)=\sum_{i=1}^m a_i e_{i;u}^*,
 \qquad
 \Lambda_u(\mathbf a,\mathbf b)
 :=\Lambda(\beta_u(\mathbf a),\beta_u(\mathbf b)).
\end{equation}
Write \(\boldsymbol\epsilon_1,\ldots,\boldsymbol\epsilon_m\) for the standard coordinate basis, and let \(\mathbf b_{k;u}\) be column \(k\) of \(\widetilde B_u\). Thus \(v_{k;u}=\beta_u(\mathbf b_{k;u})\).
For coordinate vectors, the notation \([\mathbf a]_+\) is understood componentwise.

The seed-\(u\) chart and its fraction skew field are
\[
 \mathcal T_u=\mathcal T_{u,R_{\mathbf z}}
 :=\bigoplus_{\mathbf a\in\ZZ^m}R_{\mathbf z}X_u(\mathbf a),
 \qquad
 X_u(\mathbf a)X_u(\mathbf b)
 =t^{\Lambda_u(\mathbf a,\mathbf b)}X_u(\mathbf a+\mathbf b),
\]
\[
 \mathcal F_u
 :=\operatorname{Frac}\!
       \left(\mathcal T_u\otimes_{R_{\mathbf z}}\mathbb F_{\mathbf z}\right).
\]
These are labelled copies of the quantum torus and its skew field for \((M,\Lambda).\) In a specified chart, \(x^p\) means \(X_u(\beta_u^{-1}p).\) This lattice identification is distinct from the chart embeddings into the initial skew field defined below. Rational maps are considered in \(\mathcal F_u,\) while coefficient evaluation is performed on coefficient forms and regular expressions.

For adjacent matrix seeds \(u\xrightarrow{k}u',\) abbreviate
\[
 E_{u,\varepsilon}:=E_{k,\varepsilon}^{\widetilde B_uR},
 \qquad \varepsilon\in\{+1,-1\}.
\]
Since \(\Lambda_{u'}=E_{u,\varepsilon}^T\Lambda_uE_{u,\varepsilon},\) there is a toric coordinate isomorphism
\begin{equation}\label{eq:scattering-toric-map}
 \phi_{k;u,\varepsilon}:\mathcal F_{u'}\longrightarrow\mathcal F_u,
 \qquad
 X_{u'}(\mathbf a)\longmapsto X_u(E_{u,\varepsilon}\mathbf a).
\end{equation}
When \(S_{u'}\) is obtained by the forward fixed-lattice basis mutation of Section~\ref{sec:mutation},
\begin{equation}\label{eq:forward-seed-coordinate-change}
 \beta_{u'}=\beta_u E_{u,+}.
\end{equation}
Thus \(\phi_{k;u,+}\) fixes the intrinsic monomials \(x^p.\) The negative sign gives an alternative toric factorization of the same rational mutation; it does not specify a second mutated basis.

We define the rational wall automorphisms explicitly. For \(\mathbf a\in\ZZ^m,\) put \(c=a_k.\) If \(c\ne0\), write \(\eta=\sgn(c)\) and set
\begin{equation}\label{eq:finite-rational-wall-action}
 \mathsf W_{k;u,\varepsilon}(X_u(\mathbf a))
 =
 X_u(\mathbf a)
 \prod_{j=1}^{|c|}
 P_k\!\left(
     t^{\varepsilon\eta(2j-1)}
     X_u(\varepsilon\mathbf b_{k;u})
 \right)^{-\eta}.
\end{equation}
For \(c=0\), set \(\mathsf W_{k;u,\varepsilon}(X_u(\mathbf a))=X_u(\mathbf a)\). The following lemma specifies its rational domain and its comparison with formal wall crossing.

\begin{lemma}[Rational wall actions and formal expansions]
\label{lem:rational-formal-wall-comparison}
Let \(v=v_{k;u}\), and let \(S_{k,u}\) be generated by
\[
 P_k(t^c x^v),\qquad c\in\ZZ.
\]
It is a left and right Ore set in \(\mathcal T_u\). The finite
actions \(\mathsf W_{k;u,\varepsilon}\) are automorphisms of
\(\mathcal D_{k,u}:=\mathcal T_u[S_{k,u}^{-1}]\) and of
\(\mathcal F_u\). There is an injective expansion homomorphism
\[
 \mathrm{Exp}_{S_u}:\mathcal D_{k,u}
       \longrightarrow\mathcal T_{S_u}^{\mathrm{bd}}.
\]
For the positive sign, expansion intertwines the rational action
with \(\Ad_{\Phi_k(x^v)^{-1}}\) in the pointed sectors.
For the negative sign, expansion is taken after rewriting the
factors in \(x^{-v}\) by reciprocity.
\end{lemma}

\begin{proof}
The factors are nonzero and commute. Conjugation by \(x^m\)
shifts \(c\) by \(2\Lambda(m,v)\), as in
\eqref{eq:shifted-exchange-normality}, and therefore permutes the
generators of \(S_{k,u}\). The finite common-denominator argument
given after that identity applies equally to all integer shifts.
It proves both Ore conditions. All denominators are regular since
the quantum torus is a domain.

Since \(v\in P_{S_u}^+\), each shifted polynomial is \(1\) plus
strictly positive order terms. Its inverse is the unique geometric
series in \(\widehat{\mathcal T}_{S_u,0}\). These inverses define
\(\mathrm{Exp}_{S_u}\); a numerator is finite, so its expansion
has support in finitely many translates of \(P_{S_u}\).
If \(f d^{-1}\) has zero expansion, multiply on the right by the
unit \(\mathrm{Exp}_{S_u}(d)\). Then the finite torus element
\(f\) has zero expansion. Linear independence of its based
monomials gives \(f=0\), proving injectivity.

The finite adjoint identity for \(\Phi_k\) gives
\eqref{eq:finite-rational-wall-action} on every monomial.
For the positive sign it is the expansion of the formal action
\(\Ad_{\Phi_k(x^v)^{-1}}\); that action fixes \(x^v\) and hence
every shifted denominator. Its inverse has the same properties,
so it restricts to an automorphism of \(\mathcal D_{k,u}\).
For the negative sign use
\[
 P_k(t^c x^{-v})=t^{cr}x^{-rv}P_k(t^{-c}x^v),
 \qquad r=r_k.
\]
Thus the negative-direction factors are units in the same
localization. The exceptional identity gives
\(\mathsf W_{k;u,-}=A_k\circ\mathsf W_{k;u,+}\);
the shear fixes \(v\) and preserves the denominator set. It too
is a localization automorphism. All the asserted comparisons
follow on generators and their inverses.
\end{proof}

We write \(\mathsf W^{\mathrm{rat}}\) when its rational domain
needs emphasis, and \(\mathsf W^{\mathrm{for}}_{k;u,+}\) for its
positive-direction formal expansion. The lemma compares finite
rational expressions with pointed series; it does not extend a
negative-direction dilogarithm to the whole \(P_{S_u}\)-completion.

In the formal notation for these rational actions,
\begin{equation}\label{eq:elementary-wall-map}
 \mathsf W_{k;u,\varepsilon}
 =
 \Ad_{\bigl(\Phi_{k,u}
        (X_u(\varepsilon\mathbf b_{k;u}))\bigr)^{-\varepsilon}}.
\end{equation}
Equation~\eqref{eq:elementary-wall-map} denotes \eqref{eq:finite-rational-wall-action}; it does not place the formal dilogarithm itself in \(\mathcal F_u.\) For \(\varepsilon=+1,\) this is the old incoming-wall crossing from
\(H_{k,-}\) to \(H_{k,+},\) whose Davison--Mandel crossing sign is \(-1.\)

\begin{definition}[Rational scattering mutation pullback]
\label{def:scattering-mutation-pullback}
The signed rational scattering pullback is
\begin{equation}\label{eq:scattering-mutation-pullback}
 \mu_{k;u,\varepsilon}^{\mathrm{sc}}
 :=\mathsf W_{k;u,\varepsilon}^{\mathrm{rat}}\circ\phi_{k;u,\varepsilon}:
 \mathcal F_{u'}\longrightarrow\mathcal F_u.
\end{equation}
\end{definition}

\begin{proposition}[The exchange relation from one wall]
\label{prop:scattering-exchange-relation}
For \(i\in I\),
\begin{equation}\label{eq:scattering-exchange-relation}
 \mu_{k;u,\varepsilon}^{\mathrm{sc}}
       (X_{u'}(\boldsymbol\epsilon_i))
 =
 \begin{cases}
 X_u(\boldsymbol\epsilon_i),&i\ne k,\\[3pt]
 \displaystyle\sum_{s=0}^{r_k}z_{k,s}
 X_u\!\left(
 -\boldsymbol\epsilon_k
 +s[\varepsilon\mathbf b_{k;u}]_+
 +(r_k-s)[-\varepsilon\mathbf b_{k;u}]_+
 \right),&i=k.
 \end{cases}
\end{equation}
Consequently the pullback is independent of \(\varepsilon\).
\end{proposition}

\begin{proof}
If \(i\ne k,\) then \(E_{u,\varepsilon}\boldsymbol\epsilon_i
=\boldsymbol\epsilon_i,\) and \(c=0\) in \eqref{eq:finite-rational-wall-action}. 
For \(i=k\), put
\begin{equation}\label{eq:E-on-ek}
 \mathbf a_k:=E_{u,\varepsilon}\boldsymbol\epsilon_k
 =-\boldsymbol\epsilon_k+r_k[-\varepsilon\mathbf b_{k;u}]_+.
\end{equation}
Compatibility gives
\[
 (\mathbf a_k)_k=-1,\qquad
 \Lambda_u(\mathbf a_k,\varepsilon\mathbf b_{k;u})
 =\varepsilon.
\]
Hence
\begin{align*}
 \mathsf W_{k;u,\varepsilon}(X_u(\mathbf a_k))
 &=X_u(\mathbf a_k)
   \sum_{s=0}^{r_k}z_{k,s}t^{-s\varepsilon}
            X_u(s\varepsilon\mathbf b_{k;u})\\
 &=\sum_{s=0}^{r_k}z_{k,s}
            X_u(\mathbf a_k+s\varepsilon\mathbf b_{k;u}).
\end{align*}
The Weyl factor cancels \(t^{-s\varepsilon}.\) Using \([-\varepsilon\mathbf b]_++\varepsilon\mathbf b
=[\varepsilon\mathbf b]_+\) proves the exchange formula. Changing \(\varepsilon\) and replacing \(s\) by \(r_k-s\) leaves it unchanged because \(z_{k,s}=z_{k,r_k-s}.\) Equality on the torus generators proves equality of the skew-field maps.
\end{proof}

We write \(\mu_{k;u}^{\mathrm{sc}}\) for this common map.

\subsection{The scattering atlas and its algebra}

Let \(\mathbb T_n\) be the \(n\)-regular labelled tree with root \(t_0.\) Assign the initial seed \(S\) to \(t_0,\) and define the fixed-lattice bases by forward mutation along the unique paths starting at \(t_0.\) Thus \eqref{eq:forward-seed-coordinate-change} holds on edges oriented away from the root.
The matrix seed mutation is involutive; on the reverse orientation
we use the rational coordinate pullback for the already assigned
endpoint charts.

For
\(t_0=u_0\xrightarrow{i_1}u_1\to\cdots
 \xrightarrow{i_s}u_s=u,\) define
\begin{equation}\label{eq:scattering-chart-embedding}
 \iota_u
 :=\mu_{i_1;u_0}^{\mathrm{sc}}\circ\cdots\circ
   \mu_{i_s;u_{s-1}}^{\mathrm{sc}}:
   \mathcal F_u\longrightarrow\mathcal F_{t_0},
 \qquad \iota_{t_0}=\id.
\end{equation}
Set
\begin{equation}\label{eq:scattering-cluster-variable}
 X_{i;u}^{\mathrm{sc}}
 :=\iota_u(X_u(\boldsymbol\epsilon_i)).
\end{equation}
For \(j\in F,\) these variables are unchanged along edges; write
\(X_j=X_{j;u}^{\mathrm{sc}}.\) They are invertible quantum variables and need not be central.

\begin{definition}[Quantum generalized cluster algebra from scattering]
\label{def:scattering-quantum-cluster-algebra}
The quantum generalized cluster algebra and its upper algebra are
\begin{align}
 \mathcal A_q^{\mathrm{sc}}(S)
 &:=
 R_{\mathbf z}\left\langle
 X_j^{\pm1}\ (j\in F),\
 X_{i;u}^{\mathrm{sc}}\ (i\in I_{\uf},\ u\in\mathbb T_n)
 \right\rangle
 \subset\mathcal F_{t_0},
 \label{eq:scattering-quantum-cluster-algebra}\\
 \mathcal U_q^{\mathrm{sc}}(S)
 &:=
 \bigcap_{u\in\mathbb T_n}\iota_u(\mathcal T_u)
 \subset\mathcal F_{t_0}.
 \label{eq:scattering-upper-algebra}
\end{align}
\end{definition}

We also write \(\mathcal A_q^{\mathrm{ord}}(S):=\mathcal A_q^{\mathrm{sc}}(S).\) The Laurent theorem below proves its inclusion in \(\mathcal U_q^{\mathrm{sc}}(S)\).

\begin{theorem}[Adjacent-chart involutivity]
\label{thm:scattering-involutivity}
For an edge \(u\xrightarrow{k}u',\) the two pullbacks are inverse:
\begin{equation}\label{eq:scattering-involution}
 \mu_{k;u'}^{\mathrm{sc}}\circ\mu_{k;u}^{\mathrm{sc}}
 =\id_{\mathcal F_{u'}},
 \qquad
 \mu_{k;u}^{\mathrm{sc}}\circ\mu_{k;u'}^{\mathrm{sc}}
 =\id_{\mathcal F_u}.
\end{equation}
Re-rooting therefore gives an isomorphic scattering atlas and algebra.
\end{theorem}

\begin{proof}
Choose sign \(\varepsilon\) on the forward edge and \(-\varepsilon\) on the reverse edge. Matrix mutation gives
\[
 E_{u',-\varepsilon}=E_{u,\varepsilon}=:E,\qquad
 E^2=I,\qquad
 \mathbf b_{k;u'}=-\mathbf b_{k;u},\qquad
 E\mathbf b_{k;u}=\mathbf b_{k;u}.
\]
Thus \(\phi_{k;u',-\varepsilon}=\phi_{k;u,\varepsilon}^{-1}.\) Put \(\phi=\phi_{k;u,\varepsilon}.\) Under \(\phi,\) the monomial \(X_{u'}(-\varepsilon\mathbf b_{k;u'})\) becomes \(X_u(\varepsilon\mathbf b_{k;u}),\) while the \(k\)-th coordinate of \(E\mathbf a\) is \(-a_k.\) Substitution into \eqref{eq:finite-rational-wall-action} therefore gives the rational conjugation identity
\[
 \phi\circ\mathsf W_{k;u',-\varepsilon}\circ\phi^{-1}
 =\mathsf W_{k;u,\varepsilon}^{-1}.
\]
This is the coordinate form of the exceptional-wall identity in
Lemma~\ref{lem:exceptional-conjugation}.
It follows that
\[
 \mathsf W_{k;u',-\varepsilon}\circ\phi^{-1}
 =
 (\mathsf W_{k;u,\varepsilon}\circ\phi)^{-1}.
\]
Sign independence proves \eqref{eq:scattering-involution}. Re-rooting at \(v\) replaces every embedding by \(\iota_v^{-1}\circ\iota_u,\) which proves the last assertion.
\end{proof}

\begin{remark}\label{rem:chart-inverse-versus-forward-basis}
The reverse edge uses the already assigned basis at its endpoint. It is not a second application of the forward basis rule. If \(S''=\mu_k(\mu_k(S))\) is instead formed by two forward basis mutations in the fixed lattice, then
\[
 e_i''=e_i-r_kb_{ik}e_k\ (i\ne k),\qquad e_k''=e_k,\qquad
 (e_i'')^*=A_{k,S}(e_i^*).
\]
The toric identification \(A_{k,S}^{-1}\) restores the original basis. Similarly, reverse tropical transport is \(T_{k,S}^{-1},\) as in Lemma~\ref{lem:T-seed-inverse}. These statements do not identify the full seed-dependent order completions.
\end{remark}

\subsection{The principal classical seed pattern}
\label{subsec:principal-classical-pattern}

For comparison with \cite[Section~2.1]{FPY24}, take \(m=2n\) and the compatible principal extension
\begin{equation}\label{eq:principal-compatible-extension}
 \widetilde B_{t_0}=\begin{pmatrix}B\\ I_n\end{pmatrix},
 \qquad
 \Lambda_{t_0}
 =\begin{pmatrix}0&-I_n\\ I_n&-B\end{pmatrix},
 \qquad B^T=- B.
\end{equation}
Then \(\widetilde B_{t_0}^T\Lambda_{t_0}=[I_n\mid0].\) Write \(\widetilde B_u=\binom{B_u}{(c_{ai;u})},\) and denote the classical chart coordinates by \(x_{1;u},\ldots,x_{n;u},Y_1,\ldots,Y_n.\)
Set
\[
 y_{i;u}:=\prod_{a=1}^nY_a^{c_{ai;u}},
 \qquad
 \widehat y_{i;u}:=y_{i;u}\prod_{j=1}^nx_{j;u}^{b_{ji;u}}.
\] Let \(\mathbf z_{\mathrm{ind}}\) be the independent reciprocal coefficients \(z_{i,s}\), \(1\le s\le\lfloor r_i/2\rfloor.\) In the tropical semifield \(\mathbb P=\operatorname{Trop}(\mathbf Y,\mathbf z_{\mathrm{ind}}),\)
put
\begin{equation}\label{eq:principal-tropical-denominator}
 D_{k;u}:=\bigoplus_{s=0}^{r_k}z_{k,s}y_{k;u}^s
 =\prod_{a=1}^nY_a^{r_k\min(c_{ak;u},0)}.
\end{equation}

\begin{proposition}[Principal classical mutation]
\label{prop:principal-classical-seed-pattern}
The classical chart transitions for \eqref{eq:principal-compatible-extension} form the reciprocal generalized seed pattern with initial coefficients \(y_{i;t_0}=Y_i\) in \(\mathbb P.\)
\end{proposition}

\begin{proof}
Specializing the finite exchange formula at \(t=1\) gives
\begin{equation}\label{eq:principal-classical-exchange}
 x_{k;u'}
 =
 x_{k;u}^{-1}
 \left(\prod_{j=1}^nx_{j;u}^{[-b_{jk;u}]_+}\right)^{r_k}
 \frac{\sum_{s=0}^{r_k}z_{k,s}\widehat y_{k;u}^s}{D_{k;u}},
 \qquad x_{i;u'}=x_{i;u}\quad(i\ne k).
\end{equation}
Indeed, the frozen factor in a summand is
\[
 \prod_aY_a^{r_k[-c_{ak;u}]_++s c_{ak;u}}
 =D_{k;u}^{-1}y_{k;u}^s.
\]
The endpoint coefficients are \(1\), so the tropical minimum of each independent \(\mathbf z\)-exponent is zero. Taking the minimum of \(s c_{ak;u}\) over \(0\le s\le r_k\) proves \eqref{eq:principal-tropical-denominator}.

The lower rows of matrix mutation satisfy
\[
 c_{ai;u'}=
 \begin{cases}
 -c_{ak;u},&i=k,\\
 c_{ai;u}+r_k\bigl([-c_{ak;u}]_+b_{ki;u}
                         +c_{ak;u}[b_{ki;u}]_+\bigr),&i\ne k.
 \end{cases}
\]
Therefore
\begin{equation}\label{eq:principal-classical-y-mutation}
 y_{i;u'}=
 \begin{cases}
 y_{k;u}^{-1},&i=k,\\
 y_{i;u}y_{k;u}^{r_k[b_{ki;u}]_+}D_{k;u}^{-b_{ki;u}},&i\ne k.
 \end{cases}
\end{equation}
Together with matrix mutation, these are \cite[Definition~2.3]{FPY24}. The classical expressions are subtraction-free in the independent initial variables and coefficients. Thus every cluster coordinate inverted along a finite path has nonzero classical image, and the chart expressions specialize recursively along that path.
\end{proof}

This is specialization of regular chart expressions. Passing to the coefficient group ring \(\QQ\mathbb P\) additionally inverts the independent exchange parameters.

\subsection{The universal cluster complex and cluster theta functions}
\label{subsec:universal-cluster-geometry}

All cones lie in the fixed space \(M_{\RR}\). For a seed \(S_u\), set
\[
 C_u^+:=\{p\in M_{\RR}:
          \langle e_{i;u},p\rangle\ge0\text{ for }i\in I_{\uf}\},
 \qquad
 L_F:=\sum_{j\in F}\RR e_{j;t_0}^*.
\]
Along the outward path \(t_0=u_0\xrightarrow{i_1}\cdots\xrightarrow{i_s}u_s=u,\) define
\begin{equation}\label{eq:fixed-lattice-path-transport}
 \mathscr T_u
 :=T_{i_s,S_{u_{s-1}}}\circ\cdots\circ T_{i_1,S_{u_0}},
 \qquad
 \mathscr T_{t_0}=\id,\qquad
 C_u:=\mathscr T_u^{-1}(C_u^+).
\end{equation}
On a reverse edge, transport is the inverse of its forward map. In particular, it is \(T_{k,S}^{-1},\) with the linear correction of Lemma~\ref{lem:T-seed-inverse} if written as a second forward mutation.

\begin{proposition}[Universal reachable cluster fan]
\label{prop:universal-cluster-fan}
The cones \(C_u\) and their faces form a rational fan with common lineality space \(L_F\). Its quotient in \(M_\RR/L_F\) is simplicial and unimodular. Each \(C_u\) is a reachable chamber closure of the universal quantum diagram, and adjacent seed labels give chambers sharing the corresponding facet.
\end{proposition}

\begin{proof}
Lemma~\ref{lem:compatible-pstar-injective} supplies the classical injectivity assumption. By Lemma~\ref{lem:classical-support-detection}, the quantum and classical diagrams have the same essential support. The classical chamber results \cite[Lemma~3.30 and Theorem~3.31]{CKM}, together with the fixed-lattice mutation of Section~\ref{sec:mutation}, give convexity, the common-face property, and adjacency for the cones \(C_u\). Thus the classical chamber closures are also quantum chamber closures.

Each intermediate image of a reachable chamber is a chamber.
Its interior avoids the incoming hyperplanes, so every successive inverse tropical map restricts to a single linear branch there.
Consequently \(\mathscr T_u^{-1}\circ\beta_u\) on the interior of
\[
 \mathcal C^+:=\RR_{\ge0}^n\times\RR^{m-n}
\] extends uniquely to a linear isomorphism
\[
 G_u:\RR^m\longrightarrow M_{\RR}.
\]
Each branch is integral and unimodular. Moreover, each preserves the fixed form \(\Lambda\), by Lemma~\ref{lem:T-preserves-Lambda}.
Combining this with \eqref{eq:seed-coordinate-identification} gives
\begin{equation}\label{eq:universal-cluster-cone}
 \begin{gathered}
 g_{i;u}:=G_u(\boldsymbol\epsilon_i),\qquad
 C_u=G_u(\mathcal C^+)
     =\sum_{i=1}^n\RR_{\ge0}g_{i;u}+L_F,\\
 \Lambda(G_u\mathbf a,G_u\mathbf b)
 =\Lambda_u(\mathbf a,\mathbf b).
 \end{gathered}
\end{equation}
The frozen dual vectors are unchanged by seed mutation and fixed by each tropical branch. Hence \(G_u(\boldsymbol\epsilon_j)=e_{j;t_0}^*\) for \(j\in F.\) Unimodularity now gives the stated lineality and quotient structure.
\end{proof}

\paragraph{Extended \(g\)-vectors.}
In initial coordinates put
\begin{equation}\label{eq:extended-g-vector-coordinates}
 \widetilde{\mathbf g}_{i;u}
 :=\beta_{t_0}^{-1}(g_{i;u}),\qquad
 \widetilde G_u:=\beta_{t_0}^{-1}G_u.
\end{equation}
Let \(\mathsf C_u\) be the principal generalized \(C\)-matrix;
the notation distinguishes it from the chamber \(C_u\).

\begin{lemma}[Chamber markings and extended \(g\)-vectors]
\label{lem:intrinsic-extended-g-vectors}
The columns of \(\widetilde G_u\) are the extended
\(g\)-vectors of \cite[Section~2.3]{FPY24}. They start at
\(\widetilde G_{t_0}=I_m\). For \(u\xrightarrow{k}u'\),
\begin{equation}\label{eq:extended-g-mutation-recurrence}
 \begin{aligned}
 \widetilde{\mathbf g}_{i;u'}&=\widetilde{\mathbf g}_{i;u}
                  &&(i\ne k),\\
 \widetilde{\mathbf g}_{k;u'}
 &=-\widetilde{\mathbf g}_{k;u}
   +r_k\left(\sum_{j=1}^m[-b_{jk;u}]_+
                     \widetilde{\mathbf g}_{j;u}
       -\sum_{j=1}^n[-c_{jk;u}]_+\mathbf b_{j;t_0}\right).
 \end{aligned}
\end{equation}
If \(\varepsilon=\operatorname{sgn}
(\mathsf C_u\mathbf f_k)\), then
\begin{equation}\label{eq:extended-g-signed-recurrence}
 \widetilde G_{u'}=\widetilde G_uE_{k,\varepsilon}^{\widetilde B_uR},
 \qquad
 \widetilde G_u\widetilde B_u=\widetilde B_{t_0}\mathsf C_u.
\end{equation}
\end{lemma}

\begin{proof}
The classical chamber marking theorem
\cite[Section~4.1]{CKM} identifies the linear marking
\(\mathscr T_u^{-1}\beta_u\) with the leading monomial degrees
of the classical seed pattern. This applies with frozen exponents
retained. To record the recurrence explicitly, use auxiliary
principal coefficient variables \(y_j\) and assign
\[
 \deg x_j=\boldsymbol\epsilon_j\quad(1\le j\le m),\qquad
 \deg y_j=-\mathbf b_{j;t_0}\quad(1\le j\le n),\qquad
 \deg z_{j,s}=0.
\]
The classical exchange relation is homogeneous for this grading.
At the \(k\)-th exchange, its monomial prefactor contributes
\(-\widetilde{\mathbf g}_{k;u}
 +r_k\sum_j[-b_{jk;u}]_+\widetilde{\mathbf g}_{j;u}\).
Its principal coefficient normalization contributes
\(-r_k\widetilde B_{t_0}[-\mathsf C_u\mathbf f_k]_+\).
Their sum is \eqref{eq:extended-g-mutation-recurrence}; the other
columns are unchanged. These are precisely the initial condition
and recurrence (2.10) of \cite{FPY24}. Their uniqueness identifies
the marked columns, without using the quantum theta theorem below.

Proposition~2.17 of \cite{FPY24}, or its entrywise induction using
this recurrence and matrix mutation, now gives the second identity
in \eqref{eq:extended-g-signed-recurrence}. If
\(\varepsilon=+1\), the coefficient correction in
\eqref{eq:extended-g-mutation-recurrence} is zero. If
\(\varepsilon=-1\), it equals
\(r_k\widetilde B_{t_0}\mathsf C_u\mathbf f_k
 =r_k\widetilde G_u\mathbf b_{k;u}\).
In the two cases the remaining column is respectively the
\(k\)-th column of \(\widetilde G_uE_{k,+}\) and
\(\widetilde G_uE_{k,-}\). This proves the signed recurrence.
\end{proof}

Consequently
\[
 \widetilde G_u=
 \begin{pmatrix}G_u^{\mathrm{cl}}&0\\ *&I_{m-n}\end{pmatrix},
 \qquad
 \beta_{t_0}^{-1}(C_u)
 =\sum_{i=1}^n\RR_{\ge0}\widetilde{\mathbf g}_{i;u}
       +\sum_{j\in F}\RR\boldsymbol\epsilon_j.
\]
Thus the quotient fan is the generalized \(g\)-vector fan.
For the principal extension \eqref{eq:principal-compatible-extension},
\cite[Lemma~2.18]{FPY24} gives
\(\widetilde{\mathbf g}_{i;u}
 =\binom{\mathbf g_{i;u}^{\mathrm{cl}}}{0}\) for mutable \(i\).
These coordinates depend on the initial extended exchange matrix,
\(\mathbf r\), and the mutation path.

We call this fan the \emph{universal reachable cluster complex}. Repeated seed labels may describe the same cone; only the geometric cones are counted. No assertion is made that the reachable chambers exhaust the complement of the scattering support. For any choice of root \(u\), let \(\Delta_u^+\) denote the corresponding reachable fan in \(M_{\RR},\) and put
\[
 \mathcal R_u:=M\cap\bigcup_{\sigma\in\Delta_u^+}\sigma.
\]
The edge transports identify these sets of reachable lattice indices.

For \(\mathbf a=(a_1,\ldots,a_m)^T\in\ZZ_{\ge0}^n\times\ZZ^{m-n},\) define the normalized cluster monomial
\begin{equation}\label{eq:universal-normalized-cluster-monomial}
 X_u^{\mathrm{sc}}(\mathbf a)
 :=\iota_u(X_u(\mathbf a))
 =t^{-\sum_{i<j}a_i a_j(\Lambda_u)_{ij}}
     \prod_{i=1}^m(X_{i;u}^{\mathrm{sc}})^{a_i},
\end{equation}
where the product has increasing indices and $t^{-\sum_{i<j}a_i a_j(\Lambda_u)_{ij}}$ is called Weyl factor.

\begin{lemma}[Classical specialization of theta coefficients]
\label{lem:theta-classical-detection}
For any seed \(S_u\) and general endpoint \(Q,\) quantum theta functions specialize orderwise to their classical counterparts. In independent reciprocal-root coordinates \(\boldsymbol\xi\), every monomial coefficient satisfies
\begin{equation}\label{eq:theta-coefficient-classical-detection}
 c(t,\boldsymbol\xi)=0
 \quad\Longleftrightarrow\quad
 c(1,\boldsymbol\xi)=0.
\end{equation}
\end{lemma}

\begin{proof}
Use the reciprocal parametrization of Lemma~\ref{lem:classical-support-detection}, writing its independent parameters as \(\xi_{i,a}.\) Each root \(q_{ij}\) is \(1,\) some \(\xi_{i,a},\) or its inverse. Proposition~\ref{prop:theta-positivity}, including the Weyl normalization of final monomials, gives
\[
 c(t,\boldsymbol\xi)
 \in\ZZ_{\ge0}[t^{\pm1},\xi_{i,a}^{\pm1}].
\]
At \(t=1,\) coefficients with the same reciprocal-root monomial are
added. A nonzero polynomial of this form cannot vanish, proving
\eqref{eq:theta-coefficient-classical-detection}.

At each order the broken-line sums are finite. Regular adjoint specialization commutes with multiplication and coefficient extraction by Remarks~\ref{rem:no-specialization-wall-element}
and~\ref{rem:normalized-classical-limit}. Lemma~\ref{lem:classical-support-detection} identifies the specialized diagram with the classical scattering class, and equivalence preserves theta sums. This proves the specialization assertion.
\end{proof}

We next establish finiteness for reachable indices independently of the quantum chart identification. This makes the rational maps in the subsequent comparison applicable to actual Laurent polynomials.

\begin{lemma}[Finiteness for reachable theta indices]
\label{lem:reachable-theta-finiteness}
Let \(p\in\mathcal R_u,\) and let \(Q_u\in\operatorname{Int}(C_u^+)\) be general. Then
\[
 \vartheta^u_{p,Q_u}\in\mathcal T_u.
\]
After reciprocal-root extension, its monomial support equals that
of \(\vartheta^{\mathrm{cl},u}_{p,Q_u}.\) If \(p\in C_u^+\cap M,\) then
\[
 \vartheta^u_{p,Q_u}=x^p.
\]
\end{lemma}

\begin{proof}
The classical reachable-index theorem \cite[Theorem~4.8]{CKM} says that \(\vartheta^{\mathrm{cl},u}_{p,Q_u}\) is a finite Laurent polynomial. By Lemma~\ref{lem:theta-classical-detection}, precisely the same monomial coefficients vanish in the quantum series. Thus the quantum support is finite. Its coefficients already lie in \(R_{\mathbf z}\) by the orderwise coefficient construction of
Lemma~\ref{lem:reciprocal-coefficient-base-change} and Proposition~\ref{prop:theta-positivity}; the injective root extension
is used here to detect vanishing, not to infer coefficient descent.
Hence the finite series lies in \(\mathcal T_u.\)

If \(p\in C_u^+,\) the classical no-bending identity \cite[Proposition~4.5]{CKM} gives \(\vartheta^{\mathrm{cl},u}_{p,Q_u}=x^p.\) All other quantum coefficients vanish by detection. The straight contribution has coefficient \(1,\) and every genuine bend adds a nonzero element of the sharp monoid \(P_{S_u}.\) Thus the coefficient of \(x^p\) is exactly \(1\). Frozen exponents may be arbitrary integers: frozen directions pair trivially with all scattering degrees, and multiplication by their monomials changes only the usual Weyl factor.
\end{proof}

\begin{lemma}[Theta transport and chart pullbacks]
\label{lem:theta-chart-compatibility}
Let \(u\xrightarrow{k}u'\) be an outward edge, and let \(p\in\mathcal R_u.\) For general \(Q_u\in\operatorname{Int}(C_u^+)\) and \(Q_{u'}\in\operatorname{Int}(C_{u'}^+),\)
\begin{equation}\label{eq:theta-chart-compatibility}
 \mu_{k;u}^{\mathrm{sc}}
 \left(\vartheta^{u'}_{T_{k,S_u}(p),Q_{u'}}\right)
 =
 \vartheta^u_{p,Q_u}
 \qquad\text{in }\mathcal F_u.
\end{equation}
The reverse identity uses the inverse pullback and \(T_{k,S_u}^{-1}.\)
\end{lemma}

\begin{proof}
Write \(T=T_{k,S_u}\) and \(v=v_{k;u}.\) Both \(p\) and \(T(p)\) are reachable indices. Lemma~\ref{lem:reachable-theta-finiteness} therefore makes the theta functions at the two positive chambers Laurent polynomials before any rational pullback is applied.

Choose a general point in the relative interior of the \(k\)-th facet of \(C_u^+.\) Choose nearby general points \(Q_+\) and \(Q_-\) on its positive and negative sides, respectively. They can be chosen with
\[
 Q_+\in\operatorname{Int}(C_u^+),\qquad
 Q_-\in H_{k,-}^{\circ}\cap\operatorname{Int}(C_{u'}^+).
\]
Indeed, \(e_{k;u'}=-e_{k;u},\) and the other defining inequalities remain strict near the interior of this facet. Since \(T\) is the identity on \(H_{k,-},\) Corollary~\ref{cor:theta-mutation} and equivalence of the target representative give
\[
 \phi_{k;u,+}
 \left(\vartheta^{u'}_{T(p),Q_-}\right)
 =
 \vartheta^u_{p,Q_-}.
\]
Here \eqref{eq:forward-seed-coordinate-change} identifies the based monomials on both sides without an extra scalar.

Cross the common facet from \(Q_-\) to \(Q_+.\) At every finite order its normalized local wall action is that of the \(k\)-th incoming wall. Its crossing sign is \(-1.\) The theta path-transport identity of Section~\ref{subsec:broken-lines} consequently gives
\[
 \vartheta^u_{p,Q_+}
 =
 \mathsf W_{k;u,+}\!
 \left(\phi_{k;u,+}
       (\vartheta^{u'}_{T(p),Q_-})\right)
\]
as an identity of formal expansions.
This local comparison may be performed at each order; theta sums
are constant as the endpoints move inside the indicated chambers.

The input \(\phi_{k;u,+}(\vartheta^{u'}_{T(p),Q_-})\)
is a finite Laurent polynomial by
Lemma~\ref{lem:reachable-theta-finiteness}. Its rational wall
image therefore lies in the single Ore localization
\(\mathcal D_{k,u}\) of
Lemma~\ref{lem:rational-formal-wall-comparison}.
The displayed formal identity says that its injective expansion
is the expansion of \(\vartheta^u_{p,Q_+}\), itself a finite
Laurent polynomial. Injectivity gives equality in
\(\mathcal D_{k,u}\subset\mathcal F_u\).
By definition this rational wall image is the scattering
pullback on the left of \eqref{eq:theta-chart-compatibility}.

Lemma~\ref{lem:theta-chamber-constancy} replaces \(Q_+\) and
\(Q_-\) by arbitrary general endpoints in their respective
positive chambers, proving \eqref{eq:theta-chart-compatibility}.
Theorem~\ref{thm:scattering-involutivity} and
Lemma~\ref{lem:T-seed-inverse} give the reverse identity.
\end{proof}

\begin{theorem}[Cluster monomials are theta functions]
\label{thm:universal-cluster-theta}
Let \(\mathbf a\in\ZZ_{\ge0}^n\times\ZZ^{m-n},\) and let \(Q_0\in\operatorname{Int}(C_{t_0}^+)\) be general. Then
\[
 X_u^{\mathrm{sc}}(\mathbf a)
 =\vartheta^{t_0}_{G_u\mathbf a,Q_0}.
\]
In particular, \(X_{i;u}^{\mathrm{sc}}=\vartheta^{t_0}_{g_{i;u},Q_0}\) for \(i\in I_{\uf}.\)
\end{theorem}

\begin{proof}
Use the path \(t_0=u_0\xrightarrow{i_1}\cdots\xrightarrow{i_s}u_s=u.\) Set \(p_s=\beta_u(\mathbf a),\) and define recursively
\[
 p_{j-1}=T_{i_j,S_{u_{j-1}}}^{-1}(p_j).
\]
Every \(p_j\) is a reachable index in the seed-\(u_j\) fan, and
\[
 p_0=\mathscr T_u^{-1}(\beta_u(\mathbf a))=G_u\mathbf a.
\]
For general \(Q_j\in\operatorname{Int}(C_{u_j}^+),\) Lemma~\ref{lem:reachable-theta-finiteness} gives \(\vartheta^{u}_{p_s,Q_s}=X_u(\mathbf a).\) Apply Lemma~\ref{lem:theta-chart-compatibility} successively:
\[
 \mu_{i_j;u_{j-1}}^{\mathrm{sc}}
      (\vartheta^{u_j}_{p_j,Q_j})
 =\vartheta^{u_{j-1}}_{p_{j-1},Q_{j-1}}.
\]
All these theta functions are finite Laurent polynomials in their own charts, so the compositions are the rational maps defining \(\iota_u.\) Their composite gives the asserted identity, including the normalization of \eqref{eq:universal-normalized-cluster-monomial}.
\end{proof}

\begin{theorem}[Laurent phenomenon and root positivity]
\label{thm:universal-Laurent-positivity}
For reachable seeds \(u,v\) and \(\mathbf a\in\ZZ_{\ge0}^n\times\ZZ^{m-n},\)
\[
 \iota_v^{-1}\bigl(X_u^{\mathrm{sc}}(\mathbf a)\bigr)
 \in\mathcal T_{v,R_{\mathbf z}}.
\]
After extension to the universal splitting field, each Laurent coefficient admits an expression in \(\ZZ_{\ge0}[t^{\pm1},q_{ij}].\)
\end{theorem}

\begin{proof}
Re-root the existing atlas at \(v\), using the inverse pullbacks on reversed edges. The inverse tropical maps are used on the same edges, so the preceding theta identification applies to the path from \(v\) to \(u.\) It expresses \(\iota_v^{-1}(X_u^{\mathrm{sc}}(\mathbf a))\) as a theta function with a reachable index and endpoint in \(C_v^+.\) Lemma~\ref{lem:reachable-theta-finiteness} places this function in \(\mathcal T_{v,R_{\mathbf z}}.\) Its finiteness follows from classical finite support and quantum coefficient detection, without using the quantum Laurent phenomenon as an input. The coefficient assertion is Proposition~\ref{prop:theta-positivity} applied in the seed-\(v\) chart.
\end{proof}

In particular,
\begin{equation}\label{eq:universal-cluster-Laurent-inclusion}
 \mathcal A_q^{\mathrm{ord}}(S)
 =\mathcal A_q^{\mathrm{sc}}(S)
 \subseteq\mathcal U_q^{\mathrm{sc}}(S)
 \subseteq\mathcal T_{t_0,R_{\mathbf z}}.
\end{equation}
The polynomial coefficient form permits coefficient evaluation.
Root positivity refers to expressions in the roots \(q_{ij}\);
it does not assert positivity in the elementary symmetric variables
\(z_{i,s}\).
The nonvanishing test above uses independent universal reciprocal-root
parameters and is not asserted after arbitrary coefficient evaluation.

\begin{remark}[Formal and prescribed coefficients]
\label{rem:formal-and-prescribed-coefficients}
The algebra \(\mathcal A_q^{\mathrm{sc}}(S)\) is constructed here over the formal coefficient ring \(R_{\mathbf z}.\) Proposition~\ref{prop:scattering-exchange-relation} gives its recursive exchange rule, and Theorem~\ref{thm:scattering-involutivity} verifies the inverse chart transitions.

This extends the exchange rule of \cite{BCDX18,BCDX23} to formal reciprocal coefficients. Its Laurent property and root positivity are established in Theorem~\ref{thm:universal-Laurent-positivity}.
The comparison with the prescribed-coefficient algebras of \cite{BCDX18,BCDX23} is made after coefficient evaluation in
Section~\ref{sec:coefficient-evaluation}.
\end{remark}

\section{Coefficient evaluation}\label{sec:coefficient-evaluation}

We evaluate the universal coefficient forms and identify the resulting algebras with the prescribed-coefficient algebras of \cite{BCDX18,BCDX23}. Retain the fixed lattices \(N,M,\) the form \(\Lambda,\) and the
coordinate maps \(\beta_u\) of Section~\ref{sec:scattering-algebra}. Fix
\[
 h_{i,s}(t)\in\ZZ[t^{\pm1}],\qquad
 h_{i,0}(t)=h_{i,r_i}(t)=1,\qquad
 h_{i,s}(t)=h_{i,r_i-s}(t),
\]
and put
\[
 P_{i,h}(Y):=\sum_{s=0}^{r_i}h_{i,s}(t)Y^s.
\]
Reciprocity defines a \(\QQ[t^{\pm1}]\)-algebra homomorphism
\[
 \sigma_h:R_{\mathbf z}\longrightarrow k_t:=\QQ[t^{\pm1}],
 \qquad t\longmapsto t,\qquad
 z_{i,s}\longmapsto h_{i,s}(t).
\]
All coefficient evaluations below are taken on the retained \(R_{\mathbf z}\)-forms.

\subsection{Evaluated charts and scattering diagrams}

For each seed \(S_u,\) define
\[
 \mathcal T_{u,h}
 :=\bigoplus_{\mathbf a\in\ZZ^m}k_tX_{u,h}(\mathbf a),
 \qquad
 X_{u,h}(\mathbf a)X_{u,h}(\mathbf b)
 =t^{\Lambda_u(\mathbf a,\mathbf b)}X_{u,h}(\mathbf a+\mathbf b),
\]
\[
 \mathcal F_{u,h}
 :=\operatorname{Frac}
       \bigl(\mathcal T_{u,h}\otimes_{k_t}\QQ(t)\bigr).
\]
Coefficient evaluation is the quantum-torus homomorphism
\begin{equation}\label{eq:seed-chart-coefficient-evaluation}
 \operatorname{ev}_{h,u}:\mathcal T_{u,R_{\mathbf z}}
 \longrightarrow\mathcal T_{u,h},\qquad
 \sum_{\mathbf a}c_{\mathbf a}X_u(\mathbf a)
 \longmapsto
 \sum_{\mathbf a}\sigma_h(c_{\mathbf a})X_{u,h}(\mathbf a).
\end{equation}
It extends continuously to each translated \(P_{S_u}\)-completion; we use the same notation for this extension. In the seed-\(u\) chart, an intrinsic exponent \(p\in M\) denotes \(x_h^p:=X_{u,h}(\beta_u^{-1}p).\)
For an edge \(u\xrightarrow{k}u',\) define
\[
 \phi_{k;u,\varepsilon}^{(h)}(X_{u',h}(\mathbf a))
 :=X_{u,h}(E_{u,\varepsilon}\mathbf a),\qquad
 \mu_{k;u,\varepsilon}^{(h)}
 :=\mathsf W_{k;u,\varepsilon}^{(h)}
       \circ\phi_{k;u,\varepsilon}^{(h)}.
\]
Here \(\mathsf W_{k;u,\varepsilon}^{(h)}\) is the finite rational formula \eqref{eq:finite-rational-wall-action}, with \(P_k\) replaced by \(P_{k,h}\) and \(X_u\) by \(X_{u,h}.\) Its denominators are products of shifted exchange polynomials. They are nonzero: \(\mathbf b_{k;u}\ne0,\) and the constant and leading coefficients of \(P_{k,h}\) are \(1.\) The finite proofs of Proposition~\ref{prop:scattering-exchange-relation} and Theorem~\ref{thm:scattering-involutivity} therefore apply over \(k_t.\) Thus these pullbacks are independent of \(\varepsilon\) and are inverse on opposite orientations of an edge. Write \(\mu_{k;u}^{(h)}\) for the resulting isomorphism \(\mathcal F_{u',h}\to\mathcal F_{u,h}.\)

Compose these pullbacks along the rooted tree as in \eqref{eq:scattering-chart-embedding}, obtaining
\[
 \iota_u^{(h)}:\mathcal F_{u,h}\longrightarrow\mathcal F_{t_0,h},
 \qquad
 X_{i;u}^{(h)}:=\iota_u^{(h)}(X_{u,h}(\boldsymbol\epsilon_i)),
 \qquad
 X_u^{(h)}(\mathbf a):=\iota_u^{(h)}(X_{u,h}(\mathbf a)).
\]
Let \(\mathcal A_{q,h}^{\mathrm{ord}}(S)\) be the \(k_t\)-subalgebra of \(\mathcal F_{t_0,h}\) generated by the \(X_{i;u}^{(h)}\) with \(i\in I_{\uf}\) and by the invertible frozen variables. We compare this evaluated recursive pattern with the generalized quantum cluster pattern in Subsection~\ref{sec:BCDX-after-evaluation}.

\begin{lemma}[Finite-order coefficient evaluation]
\label{lem:finite-order-base-change}
For every \(N\), the map \(\sigma_h\) commutes with the
order-\(N\) normalized Lie construction, finite BCH and path
products, and adjoint actions on each truncated pointed sector.
After common finite subdivision, evaluation commutes with
reduction up to equivalence and a final deletion of identity walls.
\end{lemma}

\begin{proof}
The degree lattice, its order, and the normalized Laurent Lie form
are independent of \(\mathbf z\). On its nilpotent order-\(N\)
quotient, BCH and the adjoint exponential are finite expressions
using brackets, products, and rational numbers. All commute with
the coefficient homomorphism \(\sigma_h\), by
Lemma~\ref{lem:reciprocal-coefficient-base-change}. The same holds
for a finite path product and for the action on a fixed pointed
order quotient. The maps commute with all truncation maps.

On a common finite subdivision, reduction multiplies coincident
factors on the same degree ray and deletes identities. Evaluation
commutes with multiplication and sends identities to identities;
it can create further identities, which are then deleted. Thus the
two reduced presentations have the same local products. No
equality of their wall lists, preservation of all essential support,
or homomorphism on \(\mathbb F_{\mathbf z}\) is required.
\end{proof}

Define the evaluated scattering diagram by
\begin{equation}\label{eq:A-diagram-orderwise-base-change}
 \mathfrak D_{S,h}
 :=\operatorname{red}\!\left(
   \mathfrak D_{S,R_{\mathbf z}}
   \widehat\otimes_{\sigma_h}k_t\right).
\end{equation}
The completed base change is formed on the finite-order quotients of the normalized Lie form \(\mathfrak g_{S,R_{\mathbf z}}\) before taking the inverse limit. Reduction is understood orderwise. The incoming wall element is denoted by \(\Phi_{i,h}(x_h^{v_i});\) its adjoint action is \eqref{eq:Delta-normalized-adjoint} with \(P_i\) replaced by \(P_{i,h}.\) Write \(\vartheta^u_{p,Q}\) for the universal theta function and \(\vartheta^{u,h}_{p,Q}\) for that of \(\mathfrak D_{S_u,h}.\)

By Theorem~\ref{thm:universal-Laurent-positivity}, the restriction
\[
 \operatorname{ev}_h
 :=\operatorname{ev}_{h,t_0}|_{\mathcal A_q^{\mathrm{sc}}(S)}
\]
is defined on the entire universal algebra. This uses its inclusion in \(\mathcal T_{t_0,R_{\mathbf z}};\) an extension of \(\sigma_h\) to \(\mathbb F_{\mathbf z}\) is unnecessary.

\subsection{Evaluation of theta functions and cluster variables}

\begin{theorem}[Evaluation of scattering and cluster functions]
\label{thm:coefficient-evaluation}
Coefficient evaluation commutes orderwise with path products, their adjoint actions, and theta functions. Moreover,
\begin{equation}\label{eq:evaluation-cluster-algebra-image}
 \operatorname{ev}_h(X_{i;u}^{\mathrm{sc}})=X_{i;u}^{(h)},
 \qquad
 \operatorname{ev}_h\bigl(\mathcal A_q^{\mathrm{sc}}(S)\bigr)
 =\mathcal A_{q,h}^{\mathrm{ord}}(S).
\end{equation}
Every normalized reachable cluster monomial is the corresponding theta function of \(\mathfrak D_{S,h}.\)
\end{theorem}

\begin{proof}
Fix an order bound \(N\) and apply Lemma~\ref{lem:finite-order-base-change} to the normalized Lie quotient and the relevant pointed sectors. Consequently, for every common admissible path \(\gamma,\)
\[
 \operatorname{ev}_{h,u}\circ
 \WC_{\gamma,\mathfrak D_{S_u,R_{\mathbf z}}}
 =\WC_{\gamma,\mathfrak D_{S_u,h}}
     \circ\operatorname{ev}_{h,u}.
\]
The finite-order group identities also specialize, so the evaluated diagram is consistent. Trivial evaluated factors may be removed. The identities are compatible with truncation and hence hold in the inverse limit.

For a common general endpoint \(Q,\) Proposition~\ref{prop:theta-positivity} gives finitely many broken-line contributions at each order. Their coefficients evaluate
to the broken-line sum for the evaluated diagram. Zero contributions disappear and coincident contributions are added; reduction preserves the sum. Thus
\begin{equation}\label{eq:evaluation-theta}
 \operatorname{ev}_{h,u}(\vartheta^u_{p,Q})
 =\vartheta^{u,h}_{p,Q}
\end{equation}
at every order and in the completion.

To identify cluster variables, use induction along the tree. For an outward edge \(u\xrightarrow{k}u',\) put
\[
 \mathbf d_s:=-\boldsymbol\epsilon_k
       +s[\mathbf b_{k;u}]_+
       +(r_k-s)[-\mathbf b_{k;u}]_+,
 \qquad
 \mathbf a_s:=\mathbf d_s+\boldsymbol\epsilon_k\in\ZZ_{\ge0}^m.
\]
Right multiplication of the exchange identity by the old coordinate gives the identity in \(\mathcal A_q^{\mathrm{sc}}(S)\)
\begin{equation}\label{eq:polynomial-exchange-for-evaluation}
 X_{k;u'}^{\mathrm{sc}}X_{k;u}^{\mathrm{sc}}
 =\sum_{s=0}^{r_k}z_{k,s}
      t^{\Lambda_u(\mathbf d_s,\boldsymbol\epsilon_k)}
      \iota_u(X_u(\mathbf a_s)).
\end{equation}
Every term on the right is a normalized polynomial in the current cluster variables. Under the induction hypothesis its evaluation is the corresponding expression in the prescribed cluster at \(u.\) The specialized old coordinate is nonzero, since the specialized
chart maps are isomorphisms. Right cancellation in \(\mathcal F_{t_0,h}\) therefore identifies the new coordinate with \(X_{k;u'}^{(h)}.\) The initial variables and frozen inverses evaluate as prescribed, proving the first equality in \eqref{eq:evaluation-cluster-algebra-image}.The generating sets give the second equality.

Evaluation preserves the Weyl normalization of cluster monomials. Combining \eqref{eq:evaluation-theta} with Theorem~\ref{thm:universal-cluster-theta} yields
\begin{equation}\label{eq:specialized-cluster-theta}
 X_u^{(h)}(\mathbf a)
 =\vartheta^{t_0,h}_{G_u\mathbf a,Q_0},
 \qquad
 \mathbf a\in\ZZ_{\ge0}^n\times\ZZ^{m-n},
\end{equation}
for general \(Q_0\in\operatorname{Int}(C_{t_0}^+).\)
\end{proof}

\begin{corollary}[Laurent phenomenon after evaluation]
\label{cor:specialized-Laurent}
Every normalized reachable cluster monomial of \(\mathcal A_{q,h}^{\mathrm{ord}}(S)\) is a finite Laurent polynomial over \(k_t\) in every reachable quantum cluster.
\end{corollary}

\begin{proof}
For \(f\in\mathcal A_q^{\mathrm{sc}}(S),\) the universal Laurent theorem gives \(\iota_u^{-1}(f)\in\mathcal T_{u,R_{\mathbf z}}.\) Re-root the two atlases at \(u,\) using the inverse chart maps on reversed edges. The same polynomial-exchange induction as above identifies the evaluations of every cluster variable and frozen inverse. Hence, on the whole universal algebra,
\begin{equation}\label{eq:evaluation-chart-compatibility}
 \operatorname{ev}_{h,u}\bigl(\iota_u^{-1}(f)\bigr)
 =(\iota_u^{(h)})^{-1}\bigl(\operatorname{ev}_h(f)\bigr).
\end{equation}
The left-hand side is the evaluation of a finite Laurent polynomial. Apply this identity to the universal cluster monomial lifting the specified specialized monomial.
\end{proof}

\begin{corollary}[Extended \(g\)-vectors after evaluation]
\label{cor:evaluation-extended-g-vectors}
Coefficient evaluation preserves the extended \(g\)-vectors \(\widetilde{\mathbf g}_{i;u}\) of \eqref{eq:extended-g-vector-coordinates}.
\end{corollary}

\begin{proof}
The pointed theta expansion and the universal Laurent theorem give a finite expression
\[
 X_{i;u}^{\mathrm{sc}}
 =X_{t_0}(\widetilde{\mathbf g}_{i;u})
  +\sum_{\mathbf a\in\ZZ_{\ge0}^n\setminus\{0\}}
     c_{\mathbf a}(t,\mathbf z)
     X_{t_0}\!\left(\widetilde{\mathbf g}_{i;u}
                     +\widetilde B_{t_0}\mathbf a\right),
 \qquad c_{\mathbf a}(t,\mathbf z)\in R_{\mathbf z}.
\]
Injectivity of \(p_{S_{t_0}}^*\) makes every displayed increment nonzero. Evaluation fixes the leading coefficient \(1\) and changes only the other coefficients. By
\eqref{eq:evaluation-cluster-algebra-image}, the evaluated expression is \(X_{i;u}^{(h)}\), with the same pointed exponent.
\end{proof}

\subsection{Comparison with the generalized quantum cluster algebras}
\label{sec:BCDX-after-evaluation}

Recall \(R=\operatorname{diag}(r_1,\ldots,r_n).\) The conventions of \cite{BCDX18,BCDX23} correspond to ours by
\begin{equation}\label{eq:BCDX-coefficient-dictionary}
 \begin{aligned}
 t&=q^{1/2},&
 d_i^{\mathrm{BCDX}}&=r_i,\\
 \widetilde B_u^{\mathrm{BCDX}}&=\widetilde B_uR,&
 \Lambda_u^{\mathrm{BCDX}}&=\Lambda_u.
 \end{aligned}
\end{equation}
In particular, the normalized \(k\)-th exchange column in their formula is
\[
 \frac{1}{d_k^{\mathrm{BCDX}}}
   (\widetilde B_u^{\mathrm{BCDX}})_{\bullet k}
 =\mathbf b_{k;u}.
\]
Our compatibility identity gives
\begin{equation}\label{eq:BCDX-compatible-diagonal}
 (\widetilde B_u^{\mathrm{BCDX}})^T\Lambda_u
 =[R\mid0],
 \qquad
 -\Lambda_u\widetilde B_u^{\mathrm{BCDX}}
 =\begin{pmatrix}R\\0\end{pmatrix}.
\end{equation}
Thus the entries denoted by \(\widetilde d_i\) in \cite{BCDX23} are \(r_i,\) while their exchange degrees \(d_i\) are \(r_i.\)

Let \(\mathcal A_{q,h}^{\mathrm{BCDX}}(S)\) denote the generalized quantum cluster algebra of \cite{BCDX18,BCDX23} for these data and the prescribed coefficients \(h_{i,s}(t),\) with scalars extended from \(\ZZ[q^{\pm1/2}]\) to \(k_t.\) Using the identity on the initial based quantum torus, regard it as a subalgebra of \(\mathcal F_{t_0,h};\) the frozen variables are inverted in both constructions.

\begin{corollary}[Recovery of the prescribed-coefficient algebra]
\label{cor:BCDX-after-evaluation}
Coefficient evaluation has image
\begin{equation}\label{eq:BCDX-evaluation-image}
 \operatorname{ev}_h\bigl(\mathcal A_q^{\mathrm{sc}}(S)\bigr)
 =\mathcal A_{q,h}^{\mathrm{ord}}(S)
 =\mathcal A_{q,h}^{\mathrm{BCDX}}(S),
\end{equation}
with the same mutation pullbacks.
\end{corollary}

\begin{proof}
For an edge \(u\xrightarrow{k}u',\) let \(B_u\) be the principal part of \(\widetilde B_u.\) The matrix \(R^{-1}F_{k,\varepsilon}^{RB_u}R\) has \(k\)-th row entries \([\varepsilon r_jb_{kj;u}]_+\) for \(j\ne k\) and entry \(-1\) at \(k.\) It is therefore the ordinary right mutation matrix for \(\widetilde B_uR.\) Hence \eqref{eq:compatible-data-mutation} gives
\[
 \widetilde B_{u'}R
 =E_{u,\varepsilon}(\widetilde B_uR)
       (R^{-1}F_{k,\varepsilon}^{RB_u}R)
 =\mu_k^{\mathrm{std}}(\widetilde B_uR),
\]
where \(\mu_k^{\mathrm{std}}\) denotes ordinary matrix mutation. The skew matrix in both constructions mutates to \(E_{u,\varepsilon}^T\Lambda_uE_{u,\varepsilon}.\)

By Proposition~\ref{prop:scattering-exchange-relation}, the evaluated pullback fixes the coordinates with \(i\ne k\) and sends the \(k\)-th coordinate to
\[
 \sum_{s=0}^{r_k}h_{k,s}(t)
 X_{u,h}\!\left(
  -\boldsymbol\epsilon_k+s[\mathbf b_{k;u}]_+
  +(r_k-s)[-\mathbf b_{k;u}]_+\right).
\]
Under \eqref{eq:BCDX-coefficient-dictionary}, this is the exchange rule of \cite[Definition~2.3]{BCDX23}.
Induction from the common initial torus identifies the cluster variables, normalized monomials, and rational pullbacks. The generating sets identify the two prescribed-coefficient algebras. Theorem~\ref{thm:coefficient-evaluation} then identifies their common algebra with the evaluation image.
\end{proof}

The recursive-pattern identification \(\mathcal A_{q,h}^{\mathrm{ord}}(S)
=\mathcal A_{q,h}^{\mathrm{BCDX}}(S)\) uses only the finite exchange rules. Its realization as the image of the universal algebra above uses the universal Laurent theorem. In particular, \eqref{eq:BCDX-evaluation-image} induces a surjection
\[
 \mathcal A_q^{\mathrm{sc}}(S)
 \otimes_{R_{\mathbf z},\sigma_h}k_t
 \twoheadrightarrow\mathcal A_{q,h}^{\mathrm{BCDX}}(S).
\]
Its injectivity is not needed. All subsequent conclusions apply to this common prescribed-coefficient algebra.

\subsection{Mutation and reachable chambers after evaluation}

We interpret a transformed diagram through its normalized local wall products on finite rational subdivisions. This permits comparison of the \(R_{\mathbf z}\)-forms with the positive representatives used in Section~\ref{sec:mutation}.

\begin{lemma}[Mutation after coefficient evaluation]
\label{lem:mutation-after-coefficient-evaluation}
For \(S'=\mu_k(S),\) the evaluated transform is orderwise finite in the target completion and satisfies
\[
 T_{k,S}(\mathfrak D_{S,h})\sim\mathfrak D_{S',h}.
\]
The reverse transport uses \(T_{k,S}^{-1}.\)
\end{lemma}

\begin{proof}
After the faithful universal root extension, the normalized local products of the \(R_{\mathbf z}\)-form agree with those of \(\mathfrak D_S^+.\) Consequently the target-degree restrictions of Proposition~\ref{prop:mutated-diagram-well-defined} apply to these products: a forbidden degree has zero coefficient after an injective extension, and hence already before extension.

For target order \(N\), Lemma~\ref{lem:tropical-order-control}
gives a finite source bound \(N'(N,0)\), independent of the
coefficient values. The exceptional incoming factor has its own
explicit finite rational formula. Thus the transform to this order
uses finitely many source walls and a finite rational subdivision.

These bounds and subdivisions are independent of the coefficient values. Applying \(\sigma_h\) may delete terms, and commutes with branch reindexing and the finite normalized group operations. It therefore produces an orderwise finite target diagram. The universal nonexceptional path identities specialize on these finite quotients. At the exceptional wall, Lemma~\ref{lem:exceptional-conjugation} is an identity of finite rational actions whose denominators are shifted exchange polynomials. Their evaluations are nonzero, so this identity and its inverse also specialize. The crossing exponents remain \(-1\) for \(v_k\) and \(+1\) for \(-v_k\) along a crossing from \(H_{k,-}\) to \(H_{k,+}.\)

It follows that the local path identities of Theorem~\ref{thm:mutation} hold after evaluation at every target order. The proof of seed faithfulness applies over \(k_t,\) since the lattice data and the indeterminate \(t\) are unchanged. Hence the identities of actions give the required path-group identities. Passing to the inverse limit proves the equivalence. Reversing these identities gives the actual inverse tropical transport of Lemma~\ref{lem:T-seed-inverse}.
\end{proof}

\begin{corollary}[Specialized reachable cluster fan]
\label{cor:specialized-cluster-fan}
The cones \(C_u\) are reachable chamber closures of \(\mathfrak D_{S,h}.\) Coefficient evaluation preserves their fan, lattice markings, and face relations.
\end{corollary}

\begin{proof}
Evaluation cannot create a nonidentity local product from an identity product. Thus the interior of every universal cluster chamber remains free of essential walls. We check that its boundary facets persist.

In the seed-\(u\) coordinates, a general point of the \(k\)-th facet of \(C_u^+\) has the normalized local product of the \(k\)-th incoming wall. Indeed, compatibility identifies the effective covector of \(\sum_i a_iv_{i;u}\) with \(\sum_i a_i e_{i;u}.\) A covector proportional to \(e_{k;u}\) therefore forces all \(a_i=0\) for \(i\ne k.\) Quotienting by degrees involving these other generators leaves the single commuting incoming wall, which is already consistent. The normalized local product is consequently that incoming factor.

Its unsigned adjoint action after evaluation is
\[
 \Ad_{\Phi_{k,h}(x_h^{v_{k;u}})}
       (X_{u,h}(\boldsymbol\epsilon_k))
 =X_{u,h}(\boldsymbol\epsilon_k)
   P_{k,h}\!\left(tX_{u,h}(\mathbf b_{k;u})\right).
\]
This action is nontrivial because \(\mathbf b_{k;u}\ne0\) and \(P_{k,h}\) has constant and leading coefficients \(1.\) The crossing from \(H_{k,-}\) to \(H_{k,+}\) uses its inverse, which is also nontrivial. The leading term of \(P_{k,h}\) detects the wall by order \(r_k\); finite-order closed supports therefore contain the whole facet by density of its general points.

By Lemma~\ref{lem:mutation-after-coefficient-evaluation}, the same piecewise-linear maps transport the evaluated scattering classes. Nonexceptional local products are transported by invertible branches; an exceptional incoming wall is replaced by another nontrivial incoming wall. Applying these comparisons along the tree and using
\(\mathscr T_u^{-1}\) transports the boundary of \(C_u^+\) to the boundary of \(C_u.\) Thus the full boundary of \(C_u\) remains in the scattering support. Its connected interior cannot join another component of the wall complement, so \(C_u\) is still a chamber closure. The cones and their lattice maps are unchanged, proving the claim.
\end{proof}

This statement concerns reachable chambers. Other wall actions and individual theta coefficients may vanish under coefficient evaluation.

\subsection{Root expressions after evaluation}

\begin{corollary}[Specialized root factorization]
\label{cor:specialized-root-factorization}
Let \(K_h/\QQ(t)\) split the \(P_{i,h},\) and write
\[
 P_{i,h}(Y)=\prod_{j=1}^{r_i}(1+q_{ij}^{(h)}Y).
\]
Over \(K_h,\) the evaluated scattering class has a representative with the evaluated wall-factor expressions of Theorem~\ref{thm:seed-positive}. Its theta coefficients, including the Laurent coefficients of reachable cluster monomials, belong to the image semiring
\[
 \ZZ_{\ge0}[t^{\pm1},q_{ij}^{(h)}].
\]
\end{corollary}

\begin{proof}
Form the positive products and broken-line expressions over the free-root ring
\[
 A_{\mathrm{fr}}:=\QQ(t)[y_{ij}:i,j],
\]
with the uniform quantum parameter $t.$ Then apply
\[
 \eta_h:A_{\mathrm{fr}}\longrightarrow K_h,
 \qquad t\longmapsto t,\qquad y_{ij}\longmapsto q_{ij}^{(h)}.
\]
This and the map \(R_{\mathbf z}\xrightarrow{\sigma_h}k_t\hookrightarrow K_h\) give the same incoming coefficients, since
\[
 \eta_h\bigl(e_s(y_{i1},\ldots,y_{ir_i})\bigr)
 =e_s(q_{i1}^{(h)},\ldots,q_{ir_i}^{(h)})
 =h_{i,s}(t)=\sigma_h(z_{i,s}).
\]
Orderwise scattering uniqueness identifies the two evaluated constructions up to equivalence. Evaluating the free-root wall products gives the stated representative. At each order, evaluate the finitely many positive broken-line contributions of Proposition~\ref{prop:theta-positivity}; summing gives the assertion about theta coefficients. The cluster-monomial case follows from \eqref{eq:specialized-cluster-theta}.

Every plethystic factor is formed before evaluation. Relations among the specialized roots, including repeated roots, may identify or cancel terms; the conclusion asserts existence of the indicated expressions. No map \(K_{\mathbf z}\to K_h,\) or intrinsic plethystic operation on \(K_h,\) is required.
\end{proof}

The image-semiring statement does not assert positivity in the
elementary coefficients \(h_{i,s}(t).\) If the chosen root evaluation takes values in a Laurent polynomial semiring with nonnegative coefficients, composition gives Laurent positivity in that semiring. For arbitrary algebraic roots, the conclusion is the evaluated positive expression.

 The universal polynomial separation formula of Section~\ref{sec:fg-wall} also evaluates to a separation formula for the generalized quantum cluster algebra \cite{FPY24}; see Corollary~\ref{cor:BCDX-separation-after-evaluation}. Only reciprocity and the endpoint conditions on \(h_{i,s}(t)\) are required. In particular, this deduction does not require \(h_{i,s}(1)>0.\)
\section{Fock--Goncharov decomposition and wall-crossing automorphisms}
\label{sec:fg-wall}

We first work over the formal coefficient ring \(R_{\mathbf z},\) using the charts \(X_u,\)fraction skew fields \(\mathcal F_u,\) and embeddings \(\iota_u\) of Section~\ref{sec:scattering-algebra}. Coefficient evaluation gives the charts \(X_{u,h}\) and embeddings \(\iota_u^{(h)}\) of Section~\ref{sec:coefficient-evaluation}. The fixed lattices \(N,M\), their seed bases, and the coordinate maps \(\beta_u\) retain their previous meanings.

We identify the non-tropical factor of a mutation path with its
wall-crossing action, construct its polynomial coefficient factor
over \(R_{\mathbf z},\) and then evaluate the coefficients. This gives the separation formula for the generzlied quantum cluster algebras without a condition on \(h_{i,s}(1).\) Our convention \(t=q^{1/2}\) agrees with \cite{FPY24}; their pairing matrix and parameters are \(D=I_n\) and \(d_i=1.\) The exchange degrees remain \(r_i.\)

\subsection{The algebraic Fock--Goncharov factor}

For \(\alpha\in\ZZ^n,\) put
\[
 \widehat Y_u^\alpha:=X_u(\widetilde B_u\alpha),
 \qquad \widehat Y_{i;u}:=X_u(\mathbf b_{i;u}).
\]
Here \(\widehat Y_u^\alpha\) is a based monomial, rather than an unadjusted ordered product of the \(\widehat Y_{i;u}.\) If \(\overline{\mathbf b}\) is the mutable part of \(\mathbf b\in\ZZ^m,\) compatibility gives
\begin{equation}\label{eq:FG-basic-commutation}
 \Lambda_u(\widetilde B_u\alpha,\mathbf b)
 =\alpha^T\overline{\mathbf b},
 \qquad
 \widehat Y_u^\alpha X_u(\mathbf b)
 =t^{2\alpha^T\overline{\mathbf b}}
       X_u(\mathbf b)\widehat Y_u^\alpha.
\end{equation}
In particular, the exponent
\begin{equation}\label{eq:FG-integral-exponent}
 a_\alpha(\mathbf b):=\alpha^T\overline{\mathbf b}
\end{equation}
is integral.

For \(a\in\ZZ\), define the brace power
\begin{equation}\label{eq:FPY-brace-power}
 \mathcal P_i^{\{a\}}(Z):=
 \begin{cases}
 \displaystyle\prod_{\ell=1}^{a}P_i(t^{2\ell-1}Z),&a>0,\\[5pt]
 1,&a=0,\\[3pt]
 \displaystyle\prod_{\ell=a}^{-1}
       P_i(t^{2\ell+1}Z)^{-1},&a<0.
 \end{cases}
\end{equation}
All factors are functions of the same \(Z\), so they commute. We write \(\mathcal P_i^{-\{a\}}=(\mathcal P_i^{\{a\}})^{-1};\) this generally differs from \(\mathcal P_i^{\{-a\}}.\) After coefficient evaluation, \(\mathcal P_{i,h}^{\{a\}}\)denotes the same formula with \(P_i\) replaced by \(P_{i,h}.\)

\begin{definition}[Fock--Goncharov automorphism]
\label{def:FG-automorphism}
For \(i\in I_{\uf}\) and \(0\ne\alpha\in\ZZ^n,\) set
\begin{equation}\label{eq:FG-automorphism-exact}
 \Psi_{i;u}(\widehat Y_u^\alpha)(X_u(\mathbf b))
 :=X_u(\mathbf b)\,
       \mathcal P_i^{-\{a_\alpha(\mathbf b)\}}
                    (\widehat Y_u^\alpha).
\end{equation}
\end{definition}

\begin{lemma}\label{lem:FG-automorphism-well-defined}
The assignment \eqref{eq:FG-automorphism-exact} extends to an
automorphism of \(\mathcal F_u.\) Its inverse replaces \(-\{a_\alpha(\mathbf b)\}\) by \(\{a_\alpha(\mathbf b)\}.\) The same statements hold in the evaluated charts.
\end{lemma}

\begin{proof}
The brace powers satisfy
\[
 \mathcal P_i^{\{a+a'\}}(Z)
 =\mathcal P_i^{\{a\}}(t^{2a'}Z)\,
       \mathcal P_i^{\{a'\}}(Z).
\]
Together with \eqref{eq:FG-basic-commutation}, this identity proves multiplicativity on torus monomials. The monomial \(\widehat Y_u^\alpha\) is fixed, since \(\alpha^TB_u\alpha=0.\) The formula with the inverse brace factor is therefore an inverse homomorphism. Both maps extend to the fraction skew field.

The matrix \(\widetilde B_u\) is injective, so \(\widetilde B_u\alpha\ne0.\) Every denominator after evaluation is consequently a nonzero polynomial in a nonconstant torus monomial, with constant and leading coefficients \(1.\) The same calculation applies over \(k_t.\)
\end{proof}

\begin{proposition}[Fock--Goncharov decomposition]
\label{prop:FG-decomposition}
For an edge \(u\xrightarrow{i}u'\) and either \(\varepsilon\in\{+1,-1\},\)
\[
 \mu_{i;u}^{\mathrm{sc}}
 =\Psi_{i;u}(\widehat Y_{i;u}^{\,\varepsilon})^\varepsilon
       \circ\phi_{i;u,\varepsilon}.
\]
The same identity holds for the evaluated pullbacks.
\end{proposition}

\begin{proof}
Write \(E=E_{i,\varepsilon}^{\widetilde B_uR}.\) The non-tropical factor fixes \(X_u(\boldsymbol\epsilon_j)\) for \(j\ne i.\) On the toric image of the remaining coordinate, its action is
\[
 X_u(E\boldsymbol\epsilon_i)\,
       P_i(t^{-\varepsilon}\widehat Y_{i;u}^{\,\varepsilon}),
\]
because \(a_{\varepsilon\mathbf f_i}(E\boldsymbol\epsilon_i)
=-\varepsilon,\), where \(\mathbf f_i\) is the \(i\)-th standard basis vector of \(\ZZ^n.\) Quantum-torus multiplication gives
\[
 \sum_{s=0}^{r_i}z_{i,s}
 X_u\!\left(-\boldsymbol\epsilon_i
       +s[\varepsilon\mathbf b_{i;u}]_+
       +(r_i-s)[-\varepsilon\mathbf b_{i;u}]_+\right).
\]
Reciprocity makes this independent of \(\varepsilon\), and Proposition~\ref{prop:scattering-exchange-relation} identifies it with the mutation pullback. The calculation is finite and also holds after replacing \(z_{i,s}\) by \(h_{i,s}(t).\) This is the normalization of \cite[Proposition~3.3]{FPY24}.
\end{proof}

\subsection{Signed toric transport and wall crossing}

Write \(X=X_{t_0}\) and \(\widehat Y^\alpha=\widehat Y_{t_0}^\alpha.\) For \(\alpha\in\ZZ^n,\) the intrinsic scattering degree is \(v_\alpha=\beta_{t_0}(\widetilde B_{t_0}\alpha),\) so
\begin{equation}\label{eq:Y-scattering-dictionary}
 x^{v_\alpha}
 =X(\widetilde B_{t_0}\alpha)=\widehat Y^\alpha.
\end{equation}
For \(0\ne\alpha\in\ZZ_{\ge0}^n,\) let \(\Phi_{i,\alpha}:=\Phi_i(x^{v_\alpha}).\) Its finite adjoint formula is
\begin{equation}\label{eq:normalized-wall-action-cvector}
 \Ad_{\Phi_{i,\alpha}^{\,\eta}}(X(\mathbf b))
 =X(\mathbf b)\,
   \mathcal P_i^{\eta\{a_\alpha(\mathbf b)\}}
                (\widehat Y^\alpha),
 \qquad \eta\in\{+1,-1\}.
\end{equation}
We use these factors for incoming walls and their transports to
reachable cluster facets.

\begin{lemma}[Single-wall comparison]
\label{lem:single-wall-FG}
Under \eqref{eq:Y-scattering-dictionary},
\[
 \Ad_{\Phi_{i,\alpha}^{\,\eta}}
 =\Psi_{i;t_0}(\widehat Y^\alpha)^{-\eta}.
\]
\end{lemma}

\begin{proof}
The two sides agree on every based monomial by
\eqref{eq:FG-automorphism-exact} and
\eqref{eq:normalized-wall-action-cvector}.
\end{proof}

Fix a path \(\mathbf i=(i_1,\ldots,i_s)\) from \(t_0=u_0\) to \(u_s=u.\) Write \(\mathsf C_v\) for its principal \(C\)-matrix, reserving \(C_v\) for the chamber of Section~\ref{sec:scattering-algebra}.
Put
\[
 \mathbf c_j:=\mathsf C_{u_{j-1}}\mathbf f_{i_j},
 \qquad
 \varepsilon_j:=\operatorname{sgn}(\mathbf c_j),
 \qquad
 \mathbf c_j^+:=\varepsilon_j\mathbf c_j
       \in\ZZ_{\ge0}^n\setminus\{0\}.
\]
These signs are those before the \(j\)-th mutation. Define the signed toric product
\[
 \phi_{u_j}^{t_0}
 :=\phi_{i_1;u_0,\varepsilon_1}\circ\cdots\circ
      \phi_{i_j;u_{j-1},\varepsilon_j}.
\]
The extended \(g\)-vector convention \eqref{eq:extended-g-vector-coordinates} gives
\begin{equation}\label{eq:FG-signed-toric-transport}
 \phi_{u_j}^{t_0}(X_{u_j}(\mathbf b))
 =X(\widetilde G_{u_j}\mathbf b),\qquad
 \widetilde G_{u_j}\widetilde B_{u_j}
 =\widetilde B_{t_0}\mathsf C_{u_j}.
\end{equation}
Indeed, the signed recurrence is \(\widetilde G_{u_j}=\widetilde G_{u_{j-1}}
E_{i_j,\varepsilon_j}^{\widetilde B_{u_{j-1}}R};\) these are \cite[Proposition~2.17 and (2.14)]{FPY24}.

Choose an admissible path \(\gamma_{\mathbf i}\) from the interior of \(C_u\) to the initial chamber, crossing the successive common facets in the order \(s,s-1,\ldots,1.\) On the \(j\)-th facet the positive effective degree is
\(v_{\mathbf c_j^+}.\) Tropical duality \(\mathsf C_{u_{j-1}}^T G_{u_{j-1}}^{\mathrm{cl}}=I_n\) shows that \(\mathbf c_j\) points toward \(C_{u_{j-1}}.\) Consequently,
\begin{equation}\label{eq:FG-DM-crossing-sign}
 \operatorname{sgn}
   \Lambda(v_{\mathbf c_j^+},\dot\gamma_{\mathbf i})
 =\varepsilon_j,\qquad
 \eta_j=\operatorname{sgn}
   \Lambda(v_{\mathbf c_j^+},-\dot\gamma_{\mathbf i})
 =-\varepsilon_j.
\end{equation}
This is the crossing convention used throughout the paper.

\begin{theorem}[Fock--Goncharov product as wall crossing]
\label{thm:FG-wall-crossing}
For the above path, put
\[
 \Psi(\mathbf i)
 :=\Psi_{i_1;t_0}(\widehat Y^{\mathbf c_1^+})^{\varepsilon_1}
       \circ\cdots\circ
   \Psi_{i_s;t_0}(\widehat Y^{\mathbf c_s^+})^{\varepsilon_s}.
\]
Then
\[
 \Psi(\mathbf i)=\WC_{\gamma_{\mathbf i},
                            \mathfrak D_{S,R_{\mathbf z}}},
 \qquad
 \iota_u=\Psi(\mathbf i)\circ\phi_u^{t_0}.
\]
The first identity compares the rational action on based monomials
with its pointed formal expansion; the second is an identity of fraction-field maps. Both identities hold after coefficient evaluation, with \(\mathfrak D_{S,h}\) and \(\iota_u^{(h)}.\)
\end{theorem}

\begin{proof}
Mutation compatibility identifies the normalized product at each reachable facet with the transported incoming factor \(\Phi_{i_j,\mathbf c_j^+}.\) This statement concerns the local wall product and is unchanged by reduction or subdivision. Lemma~\ref{lem:single-wall-FG} and
\eqref{eq:FG-DM-crossing-sign} identify its crossing action with \(\Psi_{i_j;t_0}(\widehat Y^{\mathbf c_j^+})^{\varepsilon_j}.\) The first crossed wall acts first, so the order \(s,\ldots,1\)
gives exactly the displayed composition.

For each \(j\), the degree \(v_{\mathbf c_j^+}\) belongs to
\(P_S^+\). The proof of
Lemma~\ref{lem:rational-formal-wall-comparison} therefore applies
in the initial chart to the denominators
\[
 P_{i_j}(t^c x^{v_{\mathbf c_j^+}}),\qquad c\in\ZZ.
\]
Their constant terms are \(1\), and their remaining terms have
positive order. The resulting injective pointed expansion identifies
the rational factor with the formal crossing action on based
monomials. These identifications are compatible with finite
composition and every order truncation, proving the first identity.
No identification of completions attached to different seeds is
required.

For the second identity, let \(H=\phi_{u_{j-1}}^{t_0}.\) Equation~\eqref{eq:FG-signed-toric-transport} sends \(\widehat Y_{i_j;u_{j-1}}^{\,\varepsilon_j}\) to \(\widehat Y^{\mathbf c_j^+}.\)
The toric map \(H\) preserves the skew pairing and hence the
integer exponents in \eqref{eq:FG-automorphism-exact}.
It follows that
\[
 H\circ
 \Psi_{i_j;u_{j-1}}
       (\widehat Y_{i_j;u_{j-1}}^{\,\varepsilon_j})^{\varepsilon_j}
 \circ H^{-1}
 =\Psi_{i_j;t_0}(\widehat Y^{\mathbf c_j^+})^{\varepsilon_j}.
\]
Apply this identity successively to Proposition~\ref{prop:FG-decomposition} to move all toric factors to the right.

For evaluated coefficients, the same calculation uses the automorphisms of Lemma~\ref{lem:FG-automorphism-well-defined}. The geometric identification follows from Lemma~\ref{lem:mutation-after-coefficient-evaluation} and Corollary~\ref{cor:specialized-cluster-fan}. Thus no evaluation map on the entire fraction skew field is needed.
\end{proof}

\subsection{Universal \texorpdfstring{\(F\)}{F}-polynomials and separation}

Put \(B=B_{t_0},\) and consider the based quantum polynomial algebra
\[
 \mathcal P_{B,R_{\mathbf z}}
 :=\bigoplus_{\alpha\in\ZZ_{\ge0}^n}R_{\mathbf z}Z^\alpha,
 \qquad
 Z^\alpha Z^\beta=t^{\alpha^TB\beta}Z^{\alpha+\beta}.
\]
Set \(Z_i=Z^{\mathbf f_i}.\) Since \(\widetilde B_{t_0}^T\Lambda_{t_0}\widetilde B_{t_0}=B\) and \(\widetilde B_{t_0}\) is injective, there is an embedding
\begin{equation}\label{eq:FG-quantum-polynomial-embedding}
 \jmath_Y:\mathcal P_{B,R_{\mathbf z}}\hookrightarrow
                  \mathcal T_{t_0,R_{\mathbf z}},
 \qquad Z^\alpha\longmapsto\widehat Y^\alpha.
\end{equation}
For a polynomial \(F,\) write \(F(\widehat Y):=\jmath_Y(F).\) This convention includes the Weyl factors in an ordered expression in \(Z_1,\ldots,Z_n.\)

\begin{definition}[Wall-crossing coefficient]
\label{def:wall-crossing-F}
For a reachable cluster variable with extended \(g\)-vector \(\widetilde{\mathbf g}_{i;u}\in\ZZ^m,\) define
\begin{equation}\label{eq:wall-F-definition}
 F^{\mathrm{wc}}_{i;u}(\widehat Y)
 :=X(-\widetilde{\mathbf g}_{i;u})\,
   \WC_{\gamma_{\mathbf i},\mathfrak D_{S,R_{\mathbf z}}}
        \bigl(X(\widetilde{\mathbf g}_{i;u})\bigr).
\end{equation}
\end{definition}

\begin{theorem}[Universal polynomial separation]
\label{cor:F-polynomial-wall-crossing}
For each reachable cluster variable, there is a unique \(F_{i;u}^{\mathrm{univ}}\in\mathcal P_{B,R_{\mathbf z}}\) with constant term \(1\) such that
\begin{equation}\label{eq:F-wall-normal-order}
 X_{i;u}^{\mathrm{sc}}
 =X(\widetilde{\mathbf g}_{i;u})\,
       F_{i;u}^{\mathrm{univ}}(\widehat Y).
\end{equation}
Moreover, \(F_{i;u}^{\mathrm{univ}}(\widehat Y)
=F^{\mathrm{wc}}_{i;u}(\widehat Y).\)
\end{theorem}

\begin{proof}
Theorem~\ref{thm:universal-cluster-theta} identifies \(X_{i;u}^{\mathrm{sc}}\) with a theta function whose initial exponent is \(g_{i;u}.\) Its straight broken line has coefficient \(1,\) and every genuine bend adds a nonzero element of the sharp monoid \(P_S.\) The universal Laurent theorem, Theorem~\ref{thm:universal-Laurent-positivity}, makes its expansion finite. In initial coordinates it therefore has the form
\[
 X_{i;u}^{\mathrm{sc}}
 =\sum_{\alpha\in\ZZ_{\ge0}^n}
   c_{i;u,\alpha}\,
   X(\widetilde{\mathbf g}_{i;u}
                       +\widetilde B_{t_0}\alpha),
 \qquad c_{i;u,\alpha}\in R_{\mathbf z},\quad c_{i;u,0}=1,
\] with finitely many nonzero coefficients. Quantum-torus multiplication now gives the polynomial
\begin{equation}\label{eq:FG-universal-F-coefficients}
 F_{i;u}^{\mathrm{univ}}
 =\sum_{\alpha\in\ZZ_{\ge0}^n}
 c_{i;u,\alpha}\,
 t^{-\Lambda_{t_0}(\widetilde{\mathbf g}_{i;u},
                         \widetilde B_{t_0}\alpha)}Z^\alpha.
\end{equation} 
It has constant term \(1,\) and satisfies \eqref{eq:F-wall-normal-order}. Injectivity of \(\jmath_Y\) proves uniqueness.

Finally, \(\phi_u^{t_0}(X_u(\boldsymbol\epsilon_i))
=X(\widetilde{\mathbf g}_{i;u}).\) Theorem~\ref{thm:FG-wall-crossing} identifies the wall-crossing image of this monomial with \(X_{i;u}^{\mathrm{sc}}.\) Left multiplication by \(X(-\widetilde{\mathbf g}_{i;u})\) proves the wall-crossing assertion.
\end{proof}

\subsection{Separation after coefficient evaluation}

Let \(h_{i,s}(t)\in\ZZ[t^{\pm1}]\) satisfy
\[
 h_{i,0}(t)=h_{i,r_i}(t)=1,\qquad
 h_{i,s}(t)=h_{i,r_i-s}(t),
\]
and retain \(\sigma_h:R_{\mathbf z}\to k_t.\) Define
\[
 \mathcal P_{B,k_t}
 :=\bigoplus_{\alpha\in\ZZ_{\ge0}^n}k_tZ^\alpha,\qquad
 F_{i;u}^{(h)}:=\sigma_h(F_{i;u}^{\mathrm{univ}}),
 \qquad
 \widehat Y_h^\alpha:=X_{t_0,h}(\widetilde B_{t_0}\alpha).
\]
The map \(\sigma_h\) acts only on the coefficients of the finite polynomial. The embedding \(\jmath_{Y,h}(Z^\alpha)=\widehat Y_h^\alpha\) remains injective.

\begin{corollary}[Separation for prescribed coefficients]
\label{cor:BCDX-separation-after-evaluation}
The cluster variable in \cite{BCDX18} corresponding to \(i,u\) satisfies
\begin{equation}\label{eq:BCDX-separation-no-positivity}
 X_{i;u}^{(h)}
 =X_{t_0,h}(\widetilde{\mathbf g}_{i;u})\,
       F_{i;u}^{(h)}(\widehat Y_h),
 \qquad F_{i;u}^{(h)}\in\mathcal P_{B,k_t}.
\end{equation}
This polynomial has constant term \(1\) and is the right wall-crossing coefficient for \(\mathfrak D_{S,h}.\) No condition on \(h_{i,s}(1)\) is required.
\end{corollary}

\begin{proof}
Both sides of \eqref{eq:F-wall-normal-order} are finite Laurent expressions over \(R_{\mathbf z}.\) Apply \(\operatorname{ev}_{h,t_0}\) to this identity. Theorem~\ref{thm:coefficient-evaluation} identifies the evaluated cluster variable, and Corollary~\ref{cor:BCDX-after-evaluation} identifies its algebra with the generalized quantum cluster algebra. The coefficient map fixes \(t\) and the based monomial exponents, so its right-hand side is exactly that of \eqref{eq:BCDX-separation-no-positivity}. Evaluation preserves polynomiality and the coefficient of \(Z^0.\) Corollary~\ref{cor:evaluation-extended-g-vectors} also preserves the displayed extended \(g\)-vector. The evaluated identity of Theorem~\ref{thm:FG-wall-crossing} then yields
\[
 F_{i;u}^{(h)}(\widehat Y_h)
 =X_{t_0,h}(-\widetilde{\mathbf g}_{i;u})\,
   \WC_{\gamma_{\mathbf i},\mathfrak D_{S,h}}
      \bigl(X_{t_0,h}(\widetilde{\mathbf g}_{i;u})\bigr).
\]
Injectivity of \(\jmath_{Y,h}\) gives uniqueness.
\end{proof}

\begin{proposition}[Identification with the Fu--Peng--Ye \(F\)-element]
\label{prop:FPY-F-identification}
Use the same initial chart and mutation path in \cite{FPY24}, with
\(q^{1/2}=t\), \(D=I_n\), and the same \(R\),
\(\widetilde B_u\), \(\Lambda_u\), and \(h_{i,s}(t)\).
Their rational element of Definition~3.9 equals
\(F_{i;u}^{(h)}\) in the quantum \(Y\)-skew field. In particular,
\[
 \jmath_{Y,h}(F_{i;u}^{\mathrm{FPY}})
 =X_{t_0,h}(-\widetilde{\mathbf g}_{i;u})X_{i;u}^{(h)}
 =F_{i;u}^{(h)}(\widehat Y_h).
\]
\end{proposition}

\begin{proof}
The normalization and symbol correspondence are
\[
\begin{array}{c|c}
 \text{notation of \cite{FPY24}}&\text{notation here}\\ \hline
 q^{1/2},\ D,\ d(j)&t,\ I_n,\ 1\\
 r(j),\ \varepsilon_j,\ c_j^+&r_{i_j},\ \varepsilon_j,\ \mathbf c_j^+\\
 \widehat Y_{t_0}^{\alpha}&X_{t_0,h}(\widetilde B_{t_0}\alpha)\\
 (a,b)_D&a^Tb\\
 \widetilde g_{i;u}&\widetilde{\mathbf g}_{i;u}.
\end{array}
\]
Here \(d(j)\) is the inverse compatibility parameter of FPY;
it is \(1\), whereas the exchange degree is \(r_{i_j}\).
Lemma~\ref{lem:intrinsic-extended-g-vectors} verifies the last
row by the defining recurrence. The brace powers in
\eqref{eq:FPY-brace-power} are exactly FPY's shifted powers after
\(q^{1/2}=t\); they are not ordinary integer powers.

For the first \(j\) edges let \(\Psi_{\le j}\) denote the
evaluated product in Theorem~\ref{thm:FG-wall-crossing}, and put
\(\Psi_{<j}=\Psi_{\le j-1}\). FPY's intermediate factors are
\[
 L_j=\Psi_{<j}\bigl(P_{i_j,h}
                      (t\widehat Y_h^{\mathbf c_j^+})\bigr).
\]
More precisely, their shifted powers mean the transported
expressions
\(\Psi_{<j}(\mathcal P_{i_j,h}^{\{a\}}
 (\widehat Y_h^{\mathbf c_j^+}))\); this also specifies all
Weyl shifts. Fix \(g=\widetilde{\mathbf g}_{i;u}\) and set
\(a_j=(\mathbf c_j^+)^T\overline g\). Inserting the adjacent
inverse monomials gives the finite skew-field identity
\begin{equation}\label{eq:FPY-telescoping-identification}
 \begin{split}
 X_{t_0,h}(-g)\Psi_{\le s}(X_{t_0,h}(g))
 &=\prod_{j=1}^{s}
   \bigl[\Psi_{<j}(X_{t_0,h}(-g))
         \Psi_{\le j}(X_{t_0,h}(g))\bigr]\\
 &=\prod_{j=1}^{s}
   \Psi_{<j}\!\left(
      \mathcal P_{i_j,h}^{-\varepsilon_j\{a_j\}}
                   (\widehat Y_h^{\mathbf c_j^+})\right).
 \end{split}
\end{equation}
Products have increasing \(j\) from left to right. This is the
ordered product in \cite[Lemma~3.8 and Definition~3.9]{FPY24},
including \(L_j\), \(\varepsilon_j\), \(c_j^+\), \(d(j)\),
and the pairing exponent. These algebraic results require no
positivity of \(h_{i,s}(1)\).

Theorem~\ref{thm:FG-wall-crossing} identifies
\(\Psi_{\le s}(X_{t_0,h}(g))=X_{i;u}^{(h)}\).
Corollary~\ref{cor:BCDX-separation-after-evaluation} identifies
the same left quotient with \(F_{i;u}^{(h)}(\widehat Y_h)\).
The torus embedding \(\jmath_{Y,h}\) extends injectively to its
fraction skew field, so equality of the images proves equality
of the rational \(F\)-elements.
\end{proof}

\paragraph{Comparison with Fu--Peng--Ye.}
We use the numbering of the published SIGMA version of
\cite{FPY24}. Proposition~\ref{prop:FPY-F-identification} identifies
the same rational \(F\)-element, rather than just a polynomial
satisfying a similar formula. Polynomiality in
\cite[Theorem~3.11(2)]{FPY24} and the resulting separation formula
in \cite[Corollary~3.12]{FPY24} assume \(h_{i,s}(1)>0\).
Here polynomiality and the constant term \(1\) are established
over \(R_{\mathbf z}\) before evaluation. Thus, for the normalized
seeds considered here, the identified FPY element is polynomial
and its separation identity holds for arbitrary prescribed
reciprocal Laurent coefficients. This does not assert coefficientwise
positivity after an arbitrary evaluation.

\begin{remark}[Noncommutative order]
The coefficient factor is on the right. If \(X_{i;u}^{\mathrm{sc}}=X(\widetilde{\mathbf g}_{i;u})F(\widehat Y),\) then the left quotient is \(F(\widehat Y),\) whereas the right quotient is
\[
 X_{i;u}^{\mathrm{sc}}X(-\widetilde{\mathbf g}_{i;u})
 =X(\widetilde{\mathbf g}_{i;u})F(\widehat Y)
       X(-\widetilde{\mathbf g}_{i;u}).
\]
The same order is used in the evaluated formulas.
\end{remark}

\subsection{A worked \texorpdfstring{\(r=(2,1)\)}{r=(2,1)} example}

\begin{example}[A separation formula without coefficient positivity]
\label{ex:r21-worked}
Take \(m=n=2,\) with initial coordinate matrices
\[
 \Lambda_{t_0}=\widetilde B_{t_0}
 =\begin{pmatrix}0&1\\-1&0\end{pmatrix},
 \qquad (r_1,r_2)=(2,1).
\]
Then \(\widetilde B_{t_0}^T\Lambda_{t_0}=I_2,\) and the intrinsic degrees are \(v_1=-e_{2;t_0}^*,\) \(v_2=e_{1;t_0}^*.\)
Write
\[
 P_1(Y)=1+z_{1,1}Y+Y^2,\qquad P_2(Y)=1+Y.
\]
After evaluation \(z_{1,1}\mapsto h(t),\) the polynomial \(h(t)\in\ZZ[t^{\pm1}]\) is arbitrary. Choose \(\zeta\) in a splitting field with \(\zeta+\zeta^{-1}=h(t).\) Then
\[
 P_{1,h}(Y)=(1+\zeta Y)(1+\zeta^{-1}Y).
\]

\smallskip
\noindent\emph{Lift and projection.}
Adjoin central directions \(u_+,u_-,u_2\) and take the lifted degrees \(f_\pm=v_1+u_\pm,\) \(f_2=v_2+u_2.\) The three ordinary incoming factors project to
\[
 \Psi_t(y_+x^{v_1})\Psi_t(y_-x^{v_1}),
 \qquad \Psi_t(y_2x^{v_2}).
\]
The coefficient evaluation \((y_+,y_-,y_2)\mapsto(\zeta,\zeta^{-1},1)\) recovers \(\Phi_{1,h}(x^{v_1})\) and \(\Phi_{2,h}(x^{v_2}).\) Since \(\Lambda(v_1,v_2)=1,\) the order-two factors on the outgoing base ray \(-\RR_{\ge0}(v_1+v_2)\) are
\begin{equation}\label{eq:r21-projected-outgoing}
 \Psi_t(y_+y_2x^{v_1+v_2})\,
 \Psi_t(y_-y_2x^{v_1+v_2})
 \pmod{\text{order }>2}.
\end{equation}
This is the projection of the two lifted corrections of degrees \(f_++f_2\) and \(f_-+f_2.\)

\smallskip
\noindent\emph{Mutation and the right coefficient factor.}
For the first mutation, the signed toric matrix is
\[
 E_{1,+}^{\widetilde B_{t_0}R}
 =\begin{pmatrix}-1&0\\2&1\end{pmatrix},
 \qquad
 \widetilde{\mathbf g}_{1;u_1}
 =\begin{pmatrix}-1\\2\end{pmatrix}.
\]
The mutated exchange matrix is \(\widetilde B_{u_1}=-\widetilde B_{t_0}.\) In the fixed lattice \(M,\) its degrees satisfy
\[
 v_{1;u_1}=-v_1,\qquad v_{2;u_1}=v_2+2v_1.
\]
Write \(X=X_{t_0}\) and \(\widehat Y_1=X(-\boldsymbol\epsilon_2).\) The universal exchange relation is
\begin{equation}\label{eq:r21-exchange-relation}
 X_{1;u_1}^{\mathrm{sc}}
 =X(-\boldsymbol\epsilon_1+2\boldsymbol\epsilon_2)
   +z_{1,1}X(-\boldsymbol\epsilon_1+\boldsymbol\epsilon_2)
   +X(-\boldsymbol\epsilon_1).
\end{equation}
The path from \(C_{u_1}\) to the initial chamber crosses the first incoming wall from \(H_{1,-}\) to \(H_{1,+}.\) Its DM crossing exponent is \(-1,\), so
\begin{align}
 \WC_{\gamma,\mathfrak D_{S,R_{\mathbf z}}}
        (X(\widetilde{\mathbf g}_{1;u_1}))
 &=
 X(\widetilde{\mathbf g}_{1;u_1})
       P_1(t^{-1}\widehat Y_1)
 \notag\\
 &=
 X(\widetilde{\mathbf g}_{1;u_1})
       (1+z_{1,1}t^{-1}\widehat Y_1+t^{-2}\widehat Y_1^2).
 \label{eq:r21-wall-crossing-F}
\end{align}
Indeed, \(\Lambda_{t_0}(\widetilde{\mathbf g}_{1;u_1},
-\boldsymbol\epsilon_2)=1,\) so the Weyl factor \(t^s\) in \(X(\widetilde{\mathbf g}_{1;u_1})\widehat Y_1^s\) cancels \(t^{-s}.\)
Thus
\[
 F_{1;u_1}^{\mathrm{univ}}=1+z_{1,1}t^{-1}Z_1+t^{-2}Z_1^2,
 \qquad
 F_{1;u_1}^{(h)}=1+h(t)t^{-1}Z_1+t^{-2}Z_1^2.
\]
For example, \(h(t)=-1\) gives
\[
 F_{1;u_1}^{(h)}=1-t^{-1}Z_1+t^{-2}Z_1^2,
\]
although \(h(1)=-1.\) For \(h(t)=0,\) the middle term disappears and the formula remains polynomial with constant term \(1.\) Both choices satisfy reciprocity and the endpoint conditions.
\end{example}

\appendix
\section{Hypotheses and proof dependencies}
\label{app:hypothesis-audit}

This appendix records the external inputs and the order of the main
arguments. We retain the Laurent coefficient forms
\[
 k_t=\QQ[t^{\pm1}],
 \qquad
 R_{\mathbf z}=k_t\otimes_{\QQ}R_{\mathbf z}^{\mathrm{cl}}.
\]

\subsection{Davison--Mandel positivity and the central lift}

The external positivity input is Theorem~\ref{thm:DM}. It is applied to the central lattice extension, where the root variables are recorded by actual lifted degrees. The hypotheses are verified in Proposition~\ref{prop:DM-hypotheses-central-lift}; coefficient descent is proved separately.

{\small
\renewcommand{\arraystretch}{1.12}
\begin{xltabular}{\textwidth}
{@{}>{\raggedright\arraybackslash}p{.24\textwidth}
    >{\raggedright\arraybackslash}X@{}}
\toprule
\textbf{Step} & \textbf{Verification} \\
\midrule
\endfirsthead
\toprule
\textbf{Step} & \textbf{Verification} \\
\midrule
\endhead
\midrule
\multicolumn{2}{r}{\emph{(continued on next page)}} \\
\endfoot
\bottomrule
\endlastfoot

Cone and filtration &
The lifted generators span a sharp rational cone. Their coordinate
order gives the pronilpotent filtration
(Proposition~\ref{prop:DM-hypotheses-central-lift}).
For balanced degrees, \(\lvert\alpha\rvert
=\ord^\dagger(w(\alpha)+u(\alpha))=\ord(w(\alpha)).\) \\

Effective incoming directions &
Incoming wall generators have nonzero effective covectors, although
the ambient skew form may be degenerate. In the seed setting this
follows from Lemma~\ref{lem:compatible-seed-indices-active}.
General wall degrees are positive combinations of the lifted
generators and retain their full multidegrees. \\

Positive cylindrical representative &
At each order \(N,\) the finitely many relevant supporting
hyperplanes and cylindrical supports give finite rational
subdivisions in the base.
Proposition~\ref{prop:global-positive-cylinder-gluing} inserts full
positive factors of birth order \(N\) and retains all earlier factors.
Their coefficients are identified by
Lemmas~\ref{lem:compatible-generic-germs}
and~\ref{lem:raywise-bps-uniqueness}. \\

Projection and consistency &
Lemma~\ref{lem:cylinder-pullback},
Proposition~\ref{prop:cylinder-projection-formal}, and
Theorem~\ref{thm:cylinder-projection} show that cylinder projection
preserves order, local finiteness, and path products.
Thus it preserves consistency and equivalence. \\

Uniqueness and reduction &
Local BPS coefficients are unique by actual lifted degree.
Orderwise reduction preserves the local products and their
positivity on each degree ray.
Proposition~\ref{prop:positive-cylinder-choice} gives independence
up to equivalence; the polyhedral representative need not be unique. \\


\end{xltabular}
}

The central identification \eqref{eq:central-Lie-isomorphism}
retains the lifted degree and its normalization before projection.
The normalized Laurent Lie form
\eqref{eq:normalized-quantum-torus-Lie} and the seed pairing index
\eqref{eq:seed-pairing-index} keep the adjoint coefficients in the
specified Laurent rings. Plethystic factors are formed in the free
root variables before ordinary orderwise coefficient evaluation
(Lemma~\ref{lem:plethystic-base-change}); no intrinsic plethystic
operation on a splitting field is used.

\subsection{Fixed-lattice mutation and classical support detection}

All seeds use the same lattices \(N,M\) and the same skew form \(\Lambda;\) their bases and coordinate matrices change. Our conventions are
\[
 \Lambda(v_i,m)=\langle e_i,m\rangle,
 \qquad
 \eta_{\gamma,v}
 =\operatorname{sgn}\Lambda(v,-\dot\gamma).
\]
Thus a crossing from \(H_{k,-}\) to \(H_{k,+}\) has exponent \(-1\) for the old degree \(v_k\) and \(+1\) for its replacement \(-v_k.\) Reverse tropical transport uses \(T_{k,S}^{-1},\) as in Lemma~\ref{lem:T-seed-inverse}.

\begin{center}
\small
\renewcommand{\arraystretch}{1.12}
\begin{tabularx}{\textwidth}
{@{}>{\raggedright\arraybackslash}p{.24\textwidth}
    >{\raggedright\arraybackslash}X@{}}
\toprule
\textbf{Step} & \textbf{Verification} \\
\midrule
Reciprocal roots &
The triangular symmetric-function calculation in
Lemma~\ref{lem:classical-support-detection} gives a finite faithful
root-field extension. Lemma~\ref{lem:positive-reciprocal-product-nonvanishing}
proves nonvanishing of the positive reciprocal products. \\

Classical detection &
Specialization at \(t=1\) is taken on regular adjoint actions, or on
the corresponding normalized Lie forms, as in
Remarks~\ref{rem:no-specialization-wall-element}
and~\ref{rem:normalized-classical-limit}.
Lemmas~\ref{lem:oriented-support-invariant}
and~\ref{lem:classical-support-detection} identify essential oriented
support. This does not assert preservation of every BPS degree. \\

Target completion &
The classical generalized mutation theorem, used in
Proposition~\ref{prop:mutated-diagram-well-defined}, identifies the
target oriented rays. Injectivity of the seed pairing on the
scattering lattice and saturation of the target monoid give target
degree membership. The two integral branches then give a finite
source-order bound for each target-order bound; the exceptional
incoming wall is handled separately. \\

Quantum consistency &
Lemmas~\ref{lem:linear-halfspace-transport}
and~\ref{lem:exceptional-conjugation} give the quantum action identity
\eqref{eq:mutation-loop-intertwining}.
Lemma~\ref{lem:seed-adjoint-faithful} upgrades it to a path-group
identity. \\

Mutation equivalence &
Lemma~\ref{lem:incoming-data-mutation} identifies the full quantum
incoming data and verifies that the other walls are outgoing.
Theorem~\ref{thm:KS-uniqueness} then proves
Theorem~\ref{thm:mutation} in the fixed lattice \(M.\) \\
\bottomrule
\end{tabularx}
\end{center}

Here classical mutation controls the target rays; the finite-order
argument establishes quantum well-definedness, and the quantum path
calculation establishes consistency. For broken lines,
Lemma~\ref{lem:exceptional-bending-correspondence} checks complementary polynomial choices and Weyl powers in both directions. Proposition~\ref{prop:broken-lines-mutation} uses this termwise calculation for the actual transformed representative. Equivalence then identifies theta sums,
without requiring a canonical bijection between broken lines of
arbitrary equivalent representatives.

\subsection{Cluster theta functions and the Laurent phenomenon}

Proposition~\ref{prop:scattering-exchange-relation} defines the finite
exchange rule, and Theorem~\ref{thm:scattering-involutivity} verifies
the inverse chart transitions. In the principal case,
Proposition~\ref{prop:principal-classical-seed-pattern} identifies
their classical limits with the generalized principal-coefficient
pattern.
Under the standing normalization \(\widetilde B^T\Lambda=[I_n\mid0],\) Proposition~\ref{prop:universal-cluster-fan} identifies the reachable chamber fan.

The vectors \(g_{i;u}\in M\) have initial coordinates
\[
 \widetilde{\mathbf g}_{i;u}
 =\beta_{t_0}^{-1}(g_{i;u})\in\ZZ^m,
\]
the extended \(g\)-vectors of
\eqref{eq:extended-g-vector-coordinates}.

Proposition~\ref{prop:theta-positivity} is formulated on the independently defined pointed sectors of Definition~\ref{def:pointed-seed-completion}. Its triangular transition is a module automorphism, not an assertion about an unspecified global theta-algebra completion.

Lemma~\ref{lem:theta-classical-detection} proves nonvanishing of
ordinary theta coefficients at \(t=1\) in independent universal
reciprocal-root coordinates.
Lemma~\ref{lem:reachable-theta-finiteness} therefore lifts classical
finite support and the no-bending identity before rational chart
pullbacks are applied.
Lemma~\ref{lem:rational-formal-wall-comparison} proves the Ore property and injectivity of the expansion map. Using it, Lemma~\ref{lem:theta-chart-compatibility} then gives
Theorem~\ref{thm:universal-cluster-theta}, followed by
Theorem~\ref{thm:universal-Laurent-positivity}.
The coefficients already lie in the retained \(R_{\mathbf z}\)-form;
the root extension detects vanishing.
Thus the quantum Laurent phenomenon is a conclusion of this chain.

\subsection{Coefficient evaluation and the generalized quantum cluster algebra comparison}

For reciprocal \(h_{i,s}(t)\in\ZZ[t^{\pm1}]\) with \(h_{i,0}=h_{i,r_i}=1,\) the map
\[
 \sigma_h:R_{\mathbf z}\longrightarrow k_t,\qquad
 z_{i,s}\longmapsto h_{i,s}(t)
\]
is applied to the finite-order coefficient forms and finite Laurent
expressions. Theorem~\ref{thm:coefficient-evaluation} identifies the evaluated cluster variables and theta functions.
Corollary~\ref{cor:specialized-Laurent} gives Laurentness after
evaluation; it requires no nonvanishing test on the evaluated
coefficients.

The comparison with \cite{BCDX18,BCDX23} is made at this stage.
With the parameter dictionary \eqref{eq:BCDX-coefficient-dictionary},
Corollary~\ref{cor:BCDX-after-evaluation} gives
\[
 \operatorname{ev}_h\bigl(\mathcal A_q^{\mathrm{sc}}(S)\bigr)
 =\mathcal A_{q,h}^{\mathrm{ord}}(S)
 =\mathcal A_{q,h}^{\mathrm{BCDX}}(S).
\]
The generalized quantum cluster algebra is taken over \(k_t,\) with frozen variables inverted. Equality of the two prescribed-coefficient recursive patterns uses only the finite exchange rules.
Its realization as this evaluation image uses the universal Laurent
theorem.
The induced tensor base-change map is surjective; injectivity is
not asserted.

Lemma~\ref{lem:mutation-after-coefficient-evaluation} inherits the target-order bounds before evaluation. Corollary~\ref{cor:specialized-cluster-fan} preserves the reachable
fan by verifying persistence of its facets. Other walls and theta
coefficients may disappear. Corollary~\ref{cor:specialized-root-factorization} gives evaluated positive root expressions, without identifying the universal and specialized splitting fields or asserting positivity in the coefficients \(h_{i,s}(t).\)

\subsection{Fock--Goncharov decomposition and separation}

Throughout Section~\ref{sec:fg-wall},
\[
 \widetilde B_u^T\Lambda_u=[I_n\mid0].
\]
Proposition~\ref{prop:FG-decomposition} gives the single-step decomposition. The signed toric transport \eqref{eq:FG-signed-toric-transport} uses the extended \(g\)-matrix. For the reverse chamber path, the crossing exponent is
\(-\varepsilon_j,\) as in \eqref{eq:FG-DM-crossing-sign}.
These identities yield Theorem~\ref{thm:FG-wall-crossing}.

Theorem~\ref{cor:F-polynomial-wall-crossing} constructs a finite
polynomial from the pointed universal theta expansion:
\[
 X_{i;u}^{\mathrm{sc}}
 =X_{t_0}(\widetilde{\mathbf g}_{i;u})\,
      F_{i;u}^{\mathrm{univ}}(\widehat Y),
 \qquad
 F_{i;u}^{\mathrm{univ}}\in\mathcal P_{B,R_{\mathbf z}},
 \qquad
 F_{i;u}^{\mathrm{univ}}(0)=1.
\]
Here \(B=B_{t_0}\) and \(\widehat Y^\alpha=X_{t_0}(\widetilde B_{t_0}\alpha).\) The coefficient factor is on the right and is recovered by multiplication by \(X_{t_0}(-\widetilde{\mathbf g}_{i;u})\) on the left.

Applying \(\sigma_h\) gives Corollary~\ref{cor:BCDX-separation-after-evaluation}:
\[
 X_{i;u}^{(h)}
 =X_{t_0,h}(\widetilde{\mathbf g}_{i;u})\,
      F_{i;u}^{(h)}(\widehat Y_h),
 \qquad
 F_{i;u}^{(h)}
 =\sigma_h(F_{i;u}^{\mathrm{univ}})
 \in\mathcal P_{B,k_t},
 \qquad F_{i;u}^{(h)}(0)=1.
\]
Polynomiality and the constant term are established before
evaluation. Consequently this separation formula requires only
reciprocity and the endpoint conditions on \(h_{i,s}(t).\)
Proposition~\ref{prop:FPY-F-identification} verifies that the polynomial is exactly the rational element of the published \cite[Definition~3.9]{FPY24}. The hypothesis \(h_{i,s}(1)>0\) in \cite[Theorem~3.11(2) and Corollary~3.12]{FPY24} is not an input to this argument. It makes no coefficientwise positivity claim for arbitrary
prescribed coefficients.

\section*{Acknowledgements}

The author would like to thank Professor Changjian Fu for helpful
suggestions and discussions. The author also acknowledges the use of
GPT 5.6 and DeepSeek V4.1 Flash for language polishing and for
explaining certain concepts during the preparation of this manuscript.
These tools were not used to generate mathematical proofs or to
formulate the main results. The author is solely responsible for the
mathematical content and for all statements in this paper.

\bibliographystyle{plain}
\bibliography{reference}
\end{document}